\documentclass[preprint,11pt]{elsarticle}
\biboptions{compress} 

\usepackage[a4paper]{geometry}
\usepackage{lmodern} 

\usepackage{amssymb}
\usepackage{amsfonts}
\usepackage{amsmath}
\usepackage{mathtools} 
\usepackage{xfrac} 

\usepackage{graphicx}
\usepackage{caption}
\usepackage{subfig}

\usepackage{acro} 
\usepackage{comment} 
\usepackage{lineno} 
\usepackage[ruled, vlined, english, onelanguage, algo2e]{algorithm2e}

\usepackage{listings}
\usepackage{xcolor}

\definecolor{codegray}{rgb}{0.5,0.5,0.5}
\definecolor{codegraydark}{rgb}{0.3,0.3,0.3}
\definecolor{codepurple}{rgb}{0.58,0,0.82}
\definecolor{backcolour}{rgb}{0.95,0.95,0.92}

\lstdefinestyle{mystyle}{
    backgroundcolor=\color{backcolour},   
    commentstyle=\color{codegraydark},
    keywordstyle=\color{blue},
    numberstyle=\tiny\color{codegray},
    stringstyle=\color{codepurple},
    basicstyle=\ttfamily\footnotesize,
    breakatwhitespace=false,         
    breaklines=true,                 
    captionpos=b,                    
    keepspaces=true,                 
    numbers=none,                    
    showspaces=false,                
    showstringspaces=false,
    showtabs=false,                  
    tabsize=2
}

\usepackage{makecell} 
\usepackage{hyperref}
\usepackage[capitalise]{cleveref}
\Crefname{equation}{}{}

\newcommand{\orcid}[1]{\href{https://orcid.org/#1}{\includegraphics[width=10pt]{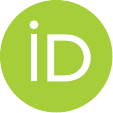}}}

\numberwithin{equation}{section}

\usepackage{fancyhdr}
\newcommand*\nid{\text{d} }   
\newcommand{\pd}[2]{\frac{\partial #1}{\partial #2} } 
\newcolumntype{?}[1]{!{\vrule width #1pt}} 

\Crefname{algocf}{Algorithm}{Algorithms}

\DeclareAcronym{ode}{
	short = ODE,
	long  = ordinary differential equation
}

\DeclareAcronym{RHS}{
	short = RHS,
	long = right-hand-side
}

\DeclareAcronym{PDE}{
	short = PDE,
	long  = partial differential equation
}

\DeclareAcronym{MoL}{
	short = MoL,
	long  = method of lines
}

\DeclareAcronym{ssp}{
	short = SSP,
	long  = strong stability preserving
}

\DeclareAcronym{tvd}{
	short = TVD,
	long  = total varitation diminishing
}

\DeclareAcronym{DG}{
	short = DG,
	long  = Discontinuous Galerkin
}

\DeclareAcronym{dgsem}{
	short = DGSEM,
	long  = discontinuous Galerkin spectral element method
}

\DeclareAcronym{hllc}{
	short = HLLC,
	long  = Harten-Lax-Van Leer-Contact
}

\DeclareAcronym{hlle}{
	short = HLLE,
	long  = Harten-Lax-Van Leer-Einfeldt
}

\DeclareAcronym{hll}{
	short = HLL,
	long  = Harten-Lax-Van Leer
}

\DeclareAcronym{mhd}{
	short = MHD,
	long  = magnetohydrodynamics
}

\DeclareAcronym{Mhd}{
	short = MHD,
	long  = Magnetohydrodynamics
}

\DeclareAcronym{vrMHD}{
	short = VRMHD,
	long  = visco-resistive Magnetohydrodynamics
}

\DeclareAcronym{GLM}{
	short = GLM,
	long  = generalized Lagrangian multiplier
}

\DeclareAcronym{glm-mhd}{
	short = GLM-MHD,
	long  = generalized Lagrangian multiplier magnetohydrodynamics
}

\DeclareAcronym{PRKM}{
	short = PRKM,
	long = Partitioned Runge-Kutta method
}
\DeclareAcronym{ARKMs}{
	short = ARKMs,
	long = Additive Runge-Kutta methods
}
\DeclareAcronym{PERK}{
	short = P-ERK,
	long = Paired-Explicit Runge-Kutta
}

\DeclareAcronym{IMEX}{
	short = IMEX,
	long = implicit-explicit
}

\DeclareAcronym{AMR}{
	short = AMR,
	long = adaptive mesh refinement
}

\DeclareAcronym{CFL}{
	short = CFL,
	long = Courant-Friedrichs-Lewy
}

\DeclareAcronym{DoF}{
	short = DoF,
	long = degree of freedom
}

\DeclareAcronym{DoFs}{
	short = DoFs,
	long = degrees of freedom
}

\journal{Journal of Computational Physics}

\begin{document}
	
	\begin{frontmatter}
		
		\title{Multirate Time-Integration based on Dynamic ODE Partitioning through Adaptively Refined Meshes for Compressible Fluid Dynamics}
		
		\cortext[cor1]{Corresponding author. E-Mail: doehring@acom.rwth-aachen.de}
		
		\author[1]{Daniel Doehring\corref{cor1} \orcid{0009-0005-4260-0332}}
		\author[1,2,3]{Michael Schlottke-Lakemper \orcid{0000-0002-3195-2536}}
		\author[4]{Gregor J. Gassner \orcid{0000-0002-1752-1158}}
		\author[1]{and\\ Manuel Torrilhon \orcid{0000-0003-0008-2061}}
		
		\affiliation[1]{
			organization={Applied~and~Computational~Mathematics,~RWTH~Aachen~University},
			addressline={\allowbreak Schinkelstrasse~2}, 
			city={Aachen},
			postcode={52062}, 
			state={North~Rhine-Westphalia},
			country={Germany.}
		}
		\affiliation[2]{
			organization={High-Performance Computing Center Stuttgart (HLRS), University of Stuttgart},
			country={Germany.}
		}
		\affiliation[3]{
			organization={High-Performance Scientific Computing, Center for Advanced Analytics and Predictive Sciences, \\ University of Augsburg},
			country = {Germany}
		}
		\affiliation[4]{
			organization={Department of Mathematics and Computer Science, Center for Data and Simulation Science, \\University of Cologne},
			country={Germany.}
		}

		\begin{abstract}
			In this paper, we apply the Paired-Explicit Runge-Kutta (P-ERK) schemes by Vermeire et. al. \cite{vermeire2019paired,nasab2022third} to dynamically partitioned systems arising from adaptive mesh refinement.
			The P-ERK schemes enable multirate time-integration with no changes in the spatial discretization methodology, making them readily implementable in existing codes that employ a method-of-lines approach.

			We show that speedup compared to a range of state of the art Runge-Kutta methods can be realized, despite additional overhead due to the dynamic re-assignment of flagging variables and restricting nonlinear stability properties.
			The effectiveness of the approach is demonstrated for a range of simulation setups for viscous and inviscid convection-dominated compressible flows for which we provide a reproducibility repository.

			In addition, we perform a thorough investigation of the nonlinear stability properties of the Paired-Explicit Runge-Kutta schemes regarding limitations due to the violation of monotonicity properties of the underlying spatial discretization.
			Furthermore, we present a novel approach for estimating the relevant eigenvalues of large Jacobians required for the optimization of stability polynomials.
		\end{abstract}
		
		%
		%
		
		%
		%
		
		
	\end{frontmatter}
	
	
	
	%
	\section{Introduction}
	Unsteady convection-dominated flows are omnipresent in nature and engineering.
	Examples thereof include turbulent Navier-Stokes equations, magnetohydrodynamics, Euler equations of gas dynamics, and aeroacoustics besides many others.
	Despite describing notably different physical phenomena, these \acp{PDE} arise from conservation laws describing the balance of mass, momentum, and energy \cite{godlewski2013numerical, leveque2002finite}.
	This allows for a unified treatment of these equations in the form of hyperbolic-parabolic balance laws. 
	When written in conserved variables $\boldsymbol u$, the governing \acp{PDE} considered in this work are of the form
	\begin{equation}
		\label{eq:ConservationLaw}
		\partial_t \boldsymbol u(t, \boldsymbol x) + \nabla \cdot \boldsymbol f\big(\boldsymbol u(t, \boldsymbol x), \nabla \boldsymbol u(t, \boldsymbol x)\big) = 
		\boldsymbol 0,
	\end{equation}
	with flux function $\boldsymbol f$, which may contain viscous terms function of $\nabla \boldsymbol u$.
	
	A common feature of the aforementioned \acp{PDE} is that over the course of a numerical simulation, solution structures of significantly different spatial scales can be present simultaneously.
	Prominent examples thereof are hydrodynamic instabilities, shock waves and turbulent eddies.
	To resolve these structures efficiently it is for practical simulations desirable to use \ac{AMR} \cite{berger1984adaptive} where the computational grid is refined precisely where a high resolution is required to resolve the relevant physical phenomena.

	Given the inherently different nature of time (directed, one dimension) and space (non-directed, multiple dimensions) the temporal and spatial derivatives in \cref{eq:ConservationLaw} are conventionally treated separately.
  This approach is commonly referred to as the \ac{MoL} \cite{cockburn2001runge}, where the application of spatial discretization techniques such as Finite Differences/\allowbreak{}Volumes/\allowbreak{}Elements or \ac{DG} leads to a typically large system of \acp{ode}
	\begin{subequations}
		\label{eq:Semidiscretization}
		\begin{align}
			\boldsymbol U(t_0) &= \boldsymbol U_0 \\
			\label{eq:SemidiscretizationODE}
			\boldsymbol U'(t) &= \boldsymbol F\big(t, \boldsymbol U(t) \big).
		\end{align}
	\end{subequations}
	%
	
	A benefit of the \ac{MoL} approach is that the \ac{PDE} can be solved in a modular fashion, i.e., the spatial discretization and the time-integration can be performed independently.
	Besides the ease of implementation, the \ac{MoL} approach benefits from the availability of a large body of literature on both the spatial discretization techniques and time-integration methods.
	This enables an independent analysis of the building blocks of the overall fully discrete system, see for instance \cite{cockburn2001runge} for the \ac{DG} method.

	When using explicit time integration methods to solve the \ac{ode} system \eqref{eq:Semidiscretization}, the maximum stable timestep is restricted by some form of \ac{CFL} condition \cite{courant1967partial} which depends on the employed spatial discretization technique, the \acp{PDE} which are solved, and the smallest grid size in the domain.
	Although implicit time-integration methods remedy this problem, they require a careful implementation to solve the potentially huge \ac{ode} system \eqref{eq:SemidiscretizationODE} in each Runge-Kutta stage efficiently.
	Additionally, if a Newton-based method is used to solve the nonlinear system, the Jacobian of the \ac{RHS} $\boldsymbol F$ needs to be computed, stored, and inverted in each stage.
	To make this process feasible in large-scale simulations, storage of the Jacobian in sparse or matrix-free format is a must, increasing the complexity of the implementation.
	Nevertheless, we remark that efficient fully implicit time-integration methods for the simulation of compressible flows are availabe, see for instance \cite{bijl2002implicit, yang2014scalable}.
	Additionally, implicit-explicit \ac{IMEX} methods \cite{ascher1995implicit} may be a viable middle ground balancing benefits and drawbacks of both implicit and explicit methods.
	
	The presence of locally refined cells results for explicit methods in a reduction of the overall timestep, causing in principle unnecessary effort in the non-refined parts of the domain, thereby reducing the computational efficiency.
	This is a classic problem which has been tackled already in the 1980s when techniques for stiff \acp{ode}, i.e., \ac{ode} problems with vastly varying time-scales where studied \cite{rentrop1985partitioned, bruder1988partitioned}.
	Since then, plenty of multirate time-integration techniques have been proposed, see for instance \cite{gunther2001multirate, constantinescu2007multirate, grote2015runge, krivodonova2010efficient, schlegel2009multirate, gunther2016multirate, cooper1983additive, jenny2020time} and references therein.
	A common feature of such techniques is that they are in general non-trivial to implement and often require changes in the spatial discretization methodology.
	In particular, for some methods intermediate values need to be stored \cite{grote2015runge}, interpolation of fine to coarse data may be required \cite{grote2015runge, krivodonova2010efficient} or nontrivial flux exchange routines need to be devised \cite{schlegel2009multirate, jenny2020time}.
	In addition, such hybrid schemes combining spatial and temporal discretization require analysis starting from scratch, with little transferability of results from the individual spatial and temporal discretization techniques \cite{constantinescu2007multirate, grote2015runge}.
	
\begin{figure}[!t]
		\centering
		\includegraphics[height=.45\textwidth]{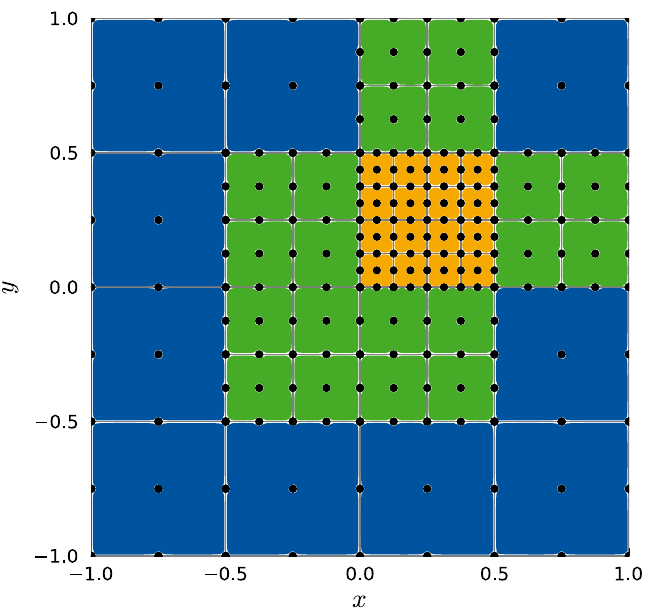}
		\caption[Illustration of the assignment of grid cells and DG coefficients to partitions $r= 1, 2, 3$ based on the minimum edge length $h$ of the grid cells.]
		{Illustration of the assignment of grid cells (grey boundaries) to partitions $r= 1, 2, 3$ based on the minimum edge length $h$ of the grid cells.
		Given a base number of stage-evaluations $E^{(1)}$ that is used to integrate the coarse blue cells, a $2E^{(1)} = E^{(2)}$ and a $4E^{(1)} = E^{(3)}$ stage-evaluation method for the green (medium) and orange (fine) cells, respectively, are required to perform the time-integration without reduction in $\text{CFL}$ number (assuming constant wave speeds $\rho_i$ across the domain).\\
		The Legendre-Gauss-Lobatto integration points of a second-order ($k=2$) \ac{DG} method are indicated by the black dots.
		Note that on each edge there are actually two points stacked on each other, one for each adjacent cell to allow for discontinuous solutions.}
		\label{fig:Grid_Schematic}
	\end{figure}

	Recently, a class of stabilized \acp{PRKM} have been proposed which achieve multirate effects while being easy to implement - in particular, there are no changes in the spatial discretization required.
	Termed \ac{PERK} \cite{vermeire2019paired}, these schemes amount only to application of different Runge-Kutta methods to the semidiscretization \eqref{eq:Semidiscretization} with no need for specialized inter/extrapolation or flux updating strategies.
	In fact, they can be implemented similar to the well-known \ac{IMEX} schemes \cite{ascher1995implicit}, which apply different Runge-Kutta methods for the different operators present in the numerical scheme.
	A schematic illustration of the application of the \ac{PERK} schemes to an adaptively refined mesh is provided in \cref{fig:Grid_Schematic}.

	Furthermore, as the \ac{PERK} schemes can be cast into the framework of \acp{PRKM}, one can utilize results on order conditions \cite{hairer1981order}, conservation of linear invariants \cite{hundsdorfer2015error} and nonlinear stability properties \cite{higueras2006strong, higueras2009characterizing, hundsdorfer2007analysis, ketcheson2013spatially}.
	In \cite{vermeire2019paired, nasab2022third} the \ac{PERK} schemes have been applied to static meshes for the simulation of the compressible Navier-Stokes equations.
	In this work, we apply the \ac{PERK} schemes to dynamically refined meshes based on \ac{AMR} which requires a dynamic re-assignment of the unknowns $\boldsymbol u$ to the corresponding integrators.
	We investigate whether speedup compared to other Runge-Kutta methods tailored to wave propagation and \ac{DG} methods can be realized, despite the additional overhead due to the re-assignment.
	Besides considering viscous flows, we also apply the \ac{PERK} schemes to purely hyperbolic problems.
	
	This paper is structured as follows: First, we introduce the class of \acp{PRKM} by revisiting central properties such as order of consistency, linear/absolute stability, internal consistency, conservation properties and criteria for nonlinear stability.
	Next, we briefly introduce the \ac{PERK} methods devised by Vermeire with a special focus on showcasing issues related to nonlinear stability.
	As the \ac{PERK} schemes can be constructed from optimized stability polynomials, we discuss the construction of the stability polynomials in \cref{sec:ConstructionStabPoly} with an emphasis on the construction of the spectrum for which we propose an efficient algorithm.
	In \cref{sec:Validation} we briefly compare the \ac{PERK} schemes to the composing standalone methods.
	Next, we turn our attention in \cref{sec:Methodology} onto describing the implementation aspects of the partitioning and devising a measure for estimation of the effectiveness of the \ac{PERK} schemes.
	Then, we present a range of classical testcases for both viscous and inviscid flows in \cref{sec:Applications} for which we provide a reproducibility repository \cite{doehring2024multirateRepro}.
	Finally, \cref{sec:Conclusions} concludes the paper.
	\section{Partitioned Runge-Kutta Methods}
	In this work, we consider coefficient-based partitioned \acp{ode} 
	of the form 
	\begin{subequations}
		\label{eq:PartitionedODESys}
		\begin{align}
			\def\arraystretch{1.4}
			\boldsymbol U(t_0) &= \boldsymbol U_0 \\
			\label{eq:PartitionedODESys2}
			\boldsymbol U'(t) &= \begin{pmatrix} \boldsymbol U^{(1)}(t) \\ \vdots \\ \boldsymbol U^{(R)}(t) \end{pmatrix}' = \begin{pmatrix} \boldsymbol F^{(1)}\big(t, \boldsymbol U^{(1)}(t), \dots,  \boldsymbol U^{(R)}(t) \big) \\ \vdots \\ \boldsymbol F^{(R)}\big(t, \boldsymbol U^{(1)}(t), \dots, \boldsymbol U^{(R)}(t) \big) \end{pmatrix} = \boldsymbol F\big(t, \boldsymbol U(t) \big)
		\end{align}
	\end{subequations}
	where $R \in \mathbb N$ denotes the number of partitions which correspond to the levels of refinement in the context of quad/octree \ac{AMR}, cf. \cref{fig:Grid_Schematic}.
	In the case of non-uniform meshes the partitioning of the semidiscretization \eqref{eq:Semidiscretization} is based on the minimum edge length $h$ of the grid cells.
	For the \ac{DG} method employed in this work, the coefficients of the local solution polynomials are the unknowns of the \ac{ode} system \eqref{eq:PartitionedODESys} and can be uniquely assigned to a partition based on the minimum edge length $h$ of the associated cell.
	In the following, the superscript $(\cdot)^{(r)}$ indicates quantities corresponding to the $r$'th partition.
	
	Systems of the form \eqref{eq:PartitionedODESys} motivate the usage of \acp{PRKM}, which are given in Butcher form given by \cite[Chapter~II.15]{HairerWanner1}
	\begin{subequations}
		\label{eq:PartitionedRKSystem}
		\begin{align}
			\label{eq:PartitionedRKFirstEq}
			\boldsymbol U_0 &= \boldsymbol U(t_0) \\
			\label{eq:PartitionedRKSecondEq}
			\boldsymbol K_i^{(r) } &= \boldsymbol F^{(r) }\left( t_n + c_i^{(r) } \Delta t, \boldsymbol U_n + \Delta t \sum_{j=1}^S \sum_{k=1}^R  a_{i,j}^{(k) } \boldsymbol K_j^{(k) }\right), \quad r = 1, \dots, R, \:\: i = 1, \dots, S \\
			\boldsymbol U_{n+1} & = \boldsymbol U_n + \Delta t  \sum_{i=1}^S  \sum_{r=1}^R   b_i^{(r) }  \boldsymbol K_i^{(r)} \, .
		\end{align}
	\end{subequations}
	Here, $S$ denotes the number of stages $\boldsymbol K_i$ and $a_{i,j}^{(r)}$ form the Butcher array matrix $A^{(r)}$.

	\acp{PRKM} have been introduced in the 1970s \cite{doi:10.1137/0713054, griepentroggemischte} and have received notable usage for the integration of Hamiltonian systems \cite{doi:10.1137/0733019, abia1993partitioned}.
	An essentially equivalent class of Runge-Kutta methods are \ac{ARKMs} \cite{KENNEDY2003139, cooper1983additive, sandu2015generalized}, which have been developed in order to achieve efficient integration of multiscale and stiff systems \cite{gunther2016multirate, cooper1983additive}.
	Similar works have been performed with \acp{PRKM} \cite{gunther2001multirate, rentrop1985partitioned, bruder1988partitioned}.
	Since the \ac{PERK} methods we apply in this work are a special case of \acp{PRKM} \cite{nasab2022third}, we briefly revisit conditions for convergence order, absolute stability, conservation of linear invariants, internal consistency, and nonlinear stability.
	\subsection{Order Conditions}
	\label{sec:OrderConditions}
	In addition to the classical conditions for convergence of order $p$, \acp{PRKM} need to satisfy additional conditions which have been derived by multiple authors, see \cite{griepentroggemischte, hairer1981order, albrecht1987new} and for a compact summary \cite{jackiewicz2000order}.
	For identical abscissae $\boldsymbol c^{(r)} \equiv \boldsymbol c$ and weights $\boldsymbol b^{(r)} \equiv \boldsymbol b$ (more on this restriction below), these in principle additional order constraints reduce for $p= 1, 2, 3$ conveniently to the classical ones:
	\begin{subequations}
		\label{eq:PRKOrderCondsSimplified}
		\begin{align}
			p=1:& &\boldsymbol b^T \boldsymbol 1 				 &\overset{!}{=} 1           \\
			p=2:& &\boldsymbol b^T \boldsymbol c 				 &\overset{!}{=} \frac{1}{2} \\
			p=3:& &\boldsymbol b^T \boldsymbol c^2 			 &\overset{!}{=} \frac{1}{3} \\ 
			&			&\boldsymbol b^T A^{(r)} \boldsymbol c &\overset{!}{=} \frac{1}{6} 
			\quad \forall \: r = 1, \dots, R \, .
		\end{align}	
	\end{subequations}
	Here, $\boldsymbol 1 \in \mathbb R^S$ denotes the vector of ones and the exponentiation of vectors is to be understood element-wise.

	For order $p=4$, however, even for identical abscissae $\boldsymbol c$ and weights $\boldsymbol b$ non-trivial coupling conditions between the Butcher arrays $A^{(r)}$ of the different methods arise \cite{jackiewicz2000order}, rendering the construction of an optimized fourth-order \ac{PERK} scheme significantly more difficult.
	In the original publication \cite{vermeire2019paired} second-order \ac{PERK} schemes are proposed which have then be extended to third-order in a subsequent work \cite{nasab2022third}.
	\subsection{Absolute Stability}
	\label{subsec:AbsoluteStability}
	Absolute stability, or linear stability, examines the asymptotic behaviour of a time integration method applied to the constant-coefficient linear test system \cite[Chapter~IV.2]{HairerWanner2}
	\begin{equation}
		\label{eq:TestEquation}
		\boldsymbol U'(t) = J \, \boldsymbol U(t).
	\end{equation}
	In this work, $J$ corresponds to the Jacobian of the \ac{RHS} $\boldsymbol F\big(t, \boldsymbol U(t)\big)$ \eqref{eq:Semidiscretization}
	\begin{equation}
		\label{eq:Jacobian}
		J(t, \boldsymbol U) \coloneqq \pd{\boldsymbol F\big(t, \boldsymbol U(t)\big)}{\boldsymbol U}.
	\end{equation}
	A time-integration method is called absolutely stable if its associated stability function $P(z)$ (which is a polynomial for explicit methods) is in magnitude less or equal to one when evaluated at the scaled eigenvalues with non-positive real part
	\begin{equation}
		\label{eq:AStab}
		\vert P( \Delta t \lambda) \vert \overset{!}{\leq} 1 \quad \forall \: \lambda \in \boldsymbol \sigma(J): \text{Re} (\lambda) \leq 0 \, .
	\end{equation}
	In other words, we demand that the scaled spectrum lies in the region of absolute stability of the method with stability polynomial $P(z)$.
	Since the \ac{RHS} $\boldsymbol F$ 
	depends on the solution $\boldsymbol U$ itself, we require \eqref{eq:AStab} to hold for all (unknown) states $\boldsymbol U(t)$ reached during the simulation.
	This is in principle difficult to ensure, but in practice the spectra $ \boldsymbol \sigma\big(J(\boldsymbol U)\big)$ are observed to be relatively similar to the spectrum of the initial condition $\boldsymbol \sigma\big(J(\boldsymbol U_0)\big)$.
	By 'similar' we are referring to the characteristic distribution of the eigenvalues in the complex plane, which remains (up to scaling due to mesh refinement) robust over the course of a simulation.
	
	The linear stability properties of \acp{PRKM} have been studied in \cite{mclachlan2011linear} with special focus on separable Hamiltonian systems.
	For \ac{ARKMs}, a linear stability analysis is performed in \cite{sandu2015generalized} based on the test problem $U' = \sum_{r=1}^R \lambda^{(r)} U$.

	For \acp{PRKM} it is customary to investigate the test equation
	\begin{equation}
		\boldsymbol U' = \Lambda  \boldsymbol U, \quad \Lambda \in \mathbb C^{N\times N}
	\end{equation}
	which can be naturally partitioned according to a set of mask-matrices $I^{(r)} \in \{0, 1\}^{N \times N}$:
	\begin{equation}
		\Lambda = \sum_{r=1}^R \Lambda^{(r)} = \sum_{r=1}^R I^{(r)} \Lambda \, .
	\end{equation}
	The mask-matrices fulfill the property $\sum_{r=1}^R I^{(r)} = I$, with $I$ being the identity matrix of dimensions $N \times N$.
	The mask matrices $I^{(r)}$ are required to specify the exact partitioning of \eqref{eq:PartitionedODESys2}:
	\begin{align}
		\boldsymbol U (t) &= \sum_{r=1}^R \underbrace{I^{(r)} \boldsymbol U (t)}_{\eqqcolon \boldsymbol U^{(r)}(t)} \, ,
		&\boldsymbol F\big(t, \boldsymbol U(t)\big) = \sum_{r=1}^R \underbrace{I^{(r)} \boldsymbol F\big(t, \boldsymbol U(t)\big)}_{\eqqcolon \boldsymbol F^{(r)}\big(t, \boldsymbol U(t) \big)}
		\, .
	\end{align}
	The matrix-valued stability function $P(Z) \in \mathbb C^{N\times N}$ is with $Z^{(r)} \coloneqq \Delta t \Lambda^{(r)}$ given by \cite{hundsdorfer2015error}
	\begin{equation}
		\label{eq:StabFuncPRK}
		P(Z) = I + \left( \sum_{r=1}^R \left(\boldsymbol b^{(r)} \otimes I \right)^T \left(Z^{(r)} \otimes I\right) \right) \left( I_{NS}  - \sum_{r=1}^R \left(A^{(r)} \otimes I\right) \left( Z^{(r)}\otimes I\right)\right)^{-1} (\boldsymbol 1 \otimes I).
	\end{equation}
	where $\otimes$ denotes the Kronecker product and $I_{NS}$ the identity matrix of dimension $N \cdot S$.
	As discussed in \cite{mclachlan2011linear}, it is difficult to infer information on the linear stability of the overall partitioned scheme from \eqref{eq:StabFuncPRK} based on the stability properties of the individual schemes $(A^{(r)}, b^{(r)})$ for arbitrary mask matrices $I^{(r)}$.
	In contrast to classic Runge-Kutta methods, one can for \acp{PRKM} no longer solely infer the action of the time-integration scheme from the stability functions of the individual methods \cite{mclachlan2011linear}.
	This issue returns in \cref{subsec:NonlinearStabPERK} when the nonlinear stability properties of a special class of \acp{PRKM} are discussed.
	
	In practice, however, it is observed that absolute stability of each method $(A^{(r)}, \boldsymbol b^{(r)})$ for its associated spectrum $\boldsymbol \lambda^{(r)}$ suffices for overall linear stability, irrespective of the concrete realization of the mask matrices $I^{(r)}$ \cite{vermeire2019paired,nasab2022third}.
	To check this, one can either compute $P(Z)$ via \eqref{eq:StabFuncPRK} or set up the fully discrete system directly, i.e., apply the Runge-Kutta method to \eqref{eq:TestEquation}.
	Then, linear stability follows from applying stability criteria for discrete-time linear time-invariant systems \cite[Chapter~8.6]{hespanha2018linear}.
	This is illustrated for a concrete example in \cref{subsubsec:LinearStabFullyDiscrete}.
	\subsection{Internal Consistency}
	To ensure that the approximations of the partitioned stages $\boldsymbol K_i^{(r)}$ approximate indeed the same timesteps $c_i$ across partitions $r$, we require internal consistency \cite{hundsdorfer2015error, hundsdorfer2007analysis, sandu2015generalized}, i.e.,
	\begin{equation}
		\label{eq:InternalConsistency}
		\sum_{j=1}^{i-1} a^{(r)}_{i,j} \overset{!}{=} c_i \quad \forall \: i = 1, \dots , S \: \: \forall \: r = 1, \dots , R \, ,
	\end{equation}
	which is equivalent to all $R$ Runge-Kutta methods having stage order one \cite{hundsdorfer2015error}.
	It has been demonstrated in \cite{hundsdorfer2007analysis, hundsdorfer2015error} that internal consistency is required to avoid spurious oscillations at the interfaces between partitions.
	\subsection{Conservation of Linear Invariants}
	Since the underlying \acp{PDE} from which the semidiscretization \eqref{eq:Semidiscretization}	is constructed correspond oftentimes to conservation laws like \eqref{eq:ConservationLaw} prescribing the conservation of mass, momentum, and energy it is natural to require that the application of a \ac{PRKM} preserves this property for conservative spatial discretizations such as finite volume and \ac{DG} schemes.
	For equation based partitioning according to \eqref{eq:PartitionedODESys2} this is ensured if the weights $\boldsymbol b$ are identical across partitions \cite{constantinescu2007multirate, hundsdorfer2015error}, i.e.,
	\begin{equation}
		\label{eq:Conservation}
		b_i^{(r_1)}	 = b_i^{(r_2)} = b_i, \: \forall \: i = 1, \dots, S, \: \: \forall \: r_1, r_2 = 1, \dots, R \, .
	\end{equation}
	For flux-based partitioning this restriction may be relaxed, see \cite{ketcheson2013spatially} for a detailed discussion.
	\subsection{Nonlinear Stability}
	Nonlinear stability includes (among stronger results for scalar equations) positivity preservation of physical quantities such as pressure and density and the suppression of spurious oscillations around discontinuities.
	Time integration methods guaranteeing these properties are termed \ac{ssp} and have been well-studied over the past years 
	\cite{shu1988efficient, gottlieb2011strong}.
	The nonlinear stability properties of partitioned/additive Runge-Kutta methods have been investigated in \cite{higueras2006strong, higueras2009characterizing, hundsdorfer2007analysis, ketcheson2013spatially}.
	In the aforementioned works it was shown that for an overall \ac{ssp} \ac{PRKM} the individual methods have to be \ac{ssp} and the timestep $\Delta t$ is furthermore restricted by 
	\begin{equation}
		\label{eq:SSPPERK}
		\Delta t = \min_{r = 1, \dots, R} \left \{c^\text{SSP}_r\right \} \Delta t_\text{Forward Euler}.
	\end{equation}
	Here, $c^\text{SSP}_r$ denotes the \ac{ssp} coefficient \cite[Chapter~3]{gottlieb2011strong} (radius of absolute monotonicity \cite{higueras2006strong, kraaijevanger1991contractivity}) of the $r$'th method.
	Except for some special cases where the timestep is restricted due to absolute stability (see, e.g., \cite{kubatko2014optimal}) \cref{eq:SSPPERK} renders the application of \ac{ssp}-capable \acp{PRKM} in many cases ineffective, as the method with the smallest \ac{ssp} coefficient determines the overall admissible timestep which guarantees nonlinear stability.

	\vspace{0.5\baselineskip} 
	Having discussed the properties of \acp{PRKM}, we now turn our attention to a special case thereof, namely the Paired-Explicit Runge-Kutta (P-ERK) methods proposed by Vermeire \cite{vermeire2019paired}.
	\section{Paired-Explicit Runge-Kutta (P-ERK) Methods}
	\label{sec:PERK}
	Paired-Explicit Runge-Kutta (P-ERK) methods have been originally proposed in \cite{vermeire2019paired} and were extended to third-order in \cite{nasab2022third}.
	These schemes target systems which exhibit locally varying characteristic speeds, most prominently due to non-uniform grids.
	In the context of the \ac{DG} schemes considered in this work, the local time scales $\beta_i$ depend both on the spectral radii $\rho_i, i = 1, \dots, N_D$ of the directional Jacobians $J_i = \partial_{\boldsymbol u} \boldsymbol f_i, i = 1, \dots, N_D$ of the flux function \eqref{eq:ConservationLaw} and the local grid size $h_i$:
	\begin{equation}
		\label{eq:CFL_DG}
		\beta \sim \min_{i = 1, \dots, N_D}\frac{h_i}{(k+1) \cdot \rho_i}.
	\end{equation}
	Here, $N_D \in \{1, 2, 3\}$ denotes the number of spatial dimensions and $k$ the local polynomial degree of the \ac{DG} approximation.
	For a more detailed discussion on the particular form of the CFL constraint for \ac{DG} methods the interested reader is referred to \cite{cockburn2001runge, krais2021flexi, warburton2008taming}.
	When applying an explicit time-integration method to the semidiscretization \eqref{eq:SemidiscretizationODE} the maximum timestep $\Delta t$ is restricted according to 
	\begin{equation}
		\label{eq:TimestepConstraint}
		\Delta t \overset{!}{\leq} \text{CFL} \cdot \beta
	\end{equation}
	where $\text{CFL} \in \mathbb R_+$ depends on the employed time-integration method.
	As discussed in \cref{sec:ConstructionStabPoly} we seek to maximize this factor by optimizing the domain of absolute stability for the different methods.

	Note that the estimate \eqref{eq:CFL_DG} is only a valid bound for the maximum timestep $\Delta t$ in the convection-dominated case, i.e., when the spectral radius $\rho$ of the Jacobian of the fully discrete system \cref{eq:SemidiscretizationODE} $J_{\boldsymbol F}(\boldsymbol U)$ scales as 
	\begin{equation}
		\label{eq:Convection_Diffusion_Scaling}
		\rho \sim \frac{\vert a \vert}{h} \gg \frac{\vert d \vert}{h^2}
	\end{equation}
	where $h$ denotes the smallest characteristic cell size and $a, d$ are suitable quantifiers for the influence of convection and diffusion, respectively.

	For such problems, one can partition the semidiscretization \eqref{eq:Semidiscretization} according to the local CFL stability constraint.
	For each partition $r$, an optimized stability polynomial $P_{p; E^{(r)}}^{(r)}(z)$ of order $p$ and degree $E^{(r)}$ is constructed which is then used to derive a specific Runge-Kutta method $(A^{(r)}, \boldsymbol b^{(r)})$.
	In particular, methods with more \ac{RHS} evaluations $E$ and larger $\text{CFL}_{E}$ are used in regions with smaller CFL number, while methods with smaller $E$ and smaller region of absolute stability (smaller $\text{CFL}_E$) are used for the elements of larger size and higher CFL.
	This is illustrated in \cref{fig:Grid_Schematic} where a three times refined grid using a quadtree refinement strategy is shown.
	This enables a more efficient time integration since the evaluations $E$ of the right-hand-side of \eqref{eq:SemidiscretizationODE} are performed only in the regions where they are actually required, and saved in other parts of the domain.
	Most important, efficiency is gained only in the coarse regions of the mesh, where cheaper integrators are used instead of the expensive one used in the finest parts.
	This idea is conceptually similar to local time-stepping methods which perform additional operations only in the fine regions of the mesh by reducing the timestep locally \cite{constantinescu2007multirate, grote2015runge}.
	\subsection{Introduction}
	\label{sec:PERK_Intro}
	To illustrate the concept of the \ac{PERK} schemes, we present the Butcher Tableaus of the second-order accurate \ac{PERK} method obtained from a four-stage and a six-stage method:
	\begin{equation}
		\label{eq:PERK_ButcherTableauClassic}
		\renewcommand\arraystretch{1.3}
		\begin{array}
			{c|c|c c c c c c|c c c c c c c c}
			i & \boldsymbol c & & & A^{(1)} & & &  & & & A^{(2)} & & & & & \\
			\hline
			1 & 0 & & & & & &             										 & & & & & & & & \\
			2 & \sfrac{1}{10} & \sfrac{1}{10} & & & & &        & \sfrac{1}{10} & & & & & \\
			3 & \sfrac{2}{10} & \sfrac{2}{10} & 0 & & & &      & \sfrac{2}{10} - a_{3,2}^{(2)} & a_{3,2}^{(2)} & & & &  \\ 
			4 & \sfrac{3}{10} & \sfrac{3}{10} & 0 & 0 & & &    & \sfrac{3}{10} -a_{4,3}^{(2)}& 0 & a_{4,3}^{(2)} & \\ 
			5 & \sfrac{4}{10} & \sfrac{4}{10} - a_{5,4}^{(1)}  & 0  & 0 & a_{5,4}^{(1)}  & & & \sfrac{4}{10} - a_{5,4}^{(2)} & 0 & 0 & a_{5,4}^{(2)} &  & \\
			6 & \sfrac{5}{10} & \sfrac{5}{10} - a_{6,5}^{(1)}  & 0 & 0 & 0 & a_{6,5}^{(1)} & & \sfrac{5}{10} - a_{6,4}^{(2)} & 0 & 0 & 0 & a_{6,4}^{(2)} & \\
			\hline
			& \boldsymbol b^T & 0 & 0 & 0 & 0 & 0 & 1          & 0 & 0 & 0 & 0             & 0 & 1
		\end{array}
	\end{equation}
	The first method $A^{(1)}$ requires only the computation of stages $\boldsymbol K_1^{(1)}, \boldsymbol K_4^{(1)}, \boldsymbol K_5^{(1)}, \boldsymbol K_6^{(1)}$, i.e., four evaluations of $\boldsymbol F^{(1)}$ while the second method $A^{(2)}$ requires computation of all stages $\boldsymbol K_i^{(2)}, i = 1, \dots, 6$, i.e., six evaluations of $\boldsymbol F^{(2)}$.
	When looking at the first method in isolation, it becomes a reducible \cite[Chapter~IV.12]{HairerWanner2} method, i.e., the stages $\boldsymbol K_2^{(1)}, \boldsymbol K_3^{(1)}$ do not influence the final solution $\boldsymbol U_{n+1}^{(1)}$ and the Butcher tableau could be truncated to a four-stage method.
	In the context of \acp{PRKM}, however, the second and third stage are required for an internally consistent cf., \cref{eq:InternalConsistency}, update of the intermediate state 
	$ \boldsymbol U_n + \Delta t \sum_{j=1}^S \sum_{r=1}^R  a_{i,j}^{(r)} \boldsymbol K_j^{(r)}$, see \eqref{eq:PartitionedRKSecondEq}.

	This particular form of the Butcher tableaus \eqref{eq:PERK_ButcherTableauClassic} comes with three advantages: First, the \ac{PERK} methods are low-storage, requiring only storage of two Runge-Kutta stages at a time.
	Second, the computational costs per stage are not increasing if higher stage methods are used, as always at most two previous stages are used for computing the next stage.
	Third, the sparse structure of both $A^{(r)}$ and $\boldsymbol b^{(r)}$ allows for a simplified computation of the coefficients from a certain stability polynomial.
	The \ac{PERK} methods are clearly internally consistent and conservative, i.e., methods which satisfy \eqref{eq:InternalConsistency} and \eqref{eq:Conservation}, respectively.
	The restriction to shared timesteps $c_i$ and weights $b_i$ comes also with the advantage of significantly simplified order constraints \cite{jackiewicz2000order},\cite[Chapter~II.15]{HairerWanner1}, cf. \eqref{eq:PRKOrderCondsSimplified}.
	
	We remark that a different choices of abscissae $\boldsymbol c$ and weights $\boldsymbol b$ than presented in \eqref{eq:PERK_ButcherTableauClassic} following \cite{vermeire2019paired} are possible.
	One options is to use $c_S = 1$ and $b_1 = b_S = \sfrac{1}{2}$ which in the $S=2$ case boils down to Heun's method.
	Another compelling option is given by $c_S = \sfrac{2}{3}$ and $b_1 = \sfrac{1}{4}, b_S = \sfrac{3}{4}$ which are known in the $S=2$ case as Ralston's method which minimizes the truncation error \cite{ralston1962runge}.
	The choice proposed in \cite{vermeire2019paired} corresponds for $S=2$ to the explicit midpoint method.
	
	The stability polynomials corresponding to \eqref{eq:PERK_ButcherTableauClassic} are given by 
	\begin{subequations}
		\label{eq:PERK_StabPolys}
		\begin{align}
			P_{2;4}^{(1)}(z) &= 1 + z +\frac{1}{2}z^2 + \alpha_3^{(1)} z^3 + \alpha_4^{(1)} z^4\\
			P_{2;6}^{(2)}(z) &= 1 + z +\frac{1}{2}z^2 + \alpha_3^{(2)} z^3 + \alpha_4^{(2)} z^4 + \alpha_5^{(2)} z^5 + \alpha_6^{(2)} z^6.
		\end{align}
	\end{subequations}
	The first subscript of $P^{(r)}$ denotes the linear order of consistency $p$ of the method, while the second subscript denotes the degree of the polynomial which is equal to the number of \ac{RHS} evaluations $E^{(r)}$.

	The free coefficients $\boldsymbol \alpha^{(1)} \in \mathbb R^{2}, \boldsymbol \alpha^{(2)} \in \mathbb R^{4}$ can be optimized for a certain spectrum $\boldsymbol \sigma\big(J_{\boldsymbol F}(\boldsymbol U)\big)$ of the Jacobian $J_{\boldsymbol F}(\boldsymbol U)$ \eqref{eq:Jacobian}.
	Using optimized stabilized Runge-Kutta methods further increases the efficiency of the scheme as the region of absolute stability can be closely adapted to the spectrum $\boldsymbol \sigma\big(J_{\boldsymbol F}(\boldsymbol U)\big)$ \cite{ketcheson2013optimal}.
	For convection-dominated problems, cf. \eqref{eq:Convection_Diffusion_Scaling}, the spectra $\boldsymbol \sigma\big(J_{\boldsymbol F}(\boldsymbol U)\big)$ of the semidiscretization increases linearly with the smallest inverse grid size $\sfrac{1}{h}$.
	At the same time, the largest possible timestep $\Delta t_{p; E}$ (in the sense of absolute stability) of an optimized $p$'th order, $E$ degree stability polynomial scales also (asymptotically) linearly with the degree = number of stage evaluations $E$ \cite{jeltsch1978largest, owren1990some, ketcheson2013optimal}.
	Consequently, the usage of an optimized $2E$ stage-evaluation Runge-Kutta method for the simulation of a convection-dominated \ac{PDE} with smallest grid size $\frac{1}{2} h$ allows using the same timestep $\Delta t$ as the simulation of the same problem on grid size $h$ with an $E$ stage-evaluation scheme.
	
	In contrast, for diffusion-dominated \acp{PDE}, the spectral radius increases quadratically with the inverse of the grid size, cf. \eqref{eq:Convection_Diffusion_Scaling}.
	Thus, for an efficient treatment of such problems the region of stability needs also to increase quadratically with the number of stages $S$, which is the case for Runge-Kutta-Chebyshev methods \cite{abdulle2015explicit, van1996development} targeting exactly such problems.
	These methods, however, do not extend without special treatment (see for instance \cite{verwer2004rkc, torrilhon2007essentially}) significantly into the complex plane, thus bearing issues for the convection-induced complex eigenvalues.
	Consequently, we restrict ourselves in this work to convection-dominated problems for which the required number of stages increases only linearly while refining the mesh.
	\subsection{Nonlinear Stability Properties}
	\label{subsec:NonlinearStabPERK}
	The key element to the effectiveness of the \ac{PERK} schemes is the fact that evaluations of the right-hand-side can be avoided in regions of the domain with accordingly increased characteristic time $\beta$.
	This, in turn, requires that the weight vector $\boldsymbol b$ is relatively sparse in order to be able to omit as many stage computations as possible.
	In the context of the previous example \cref{eq:PERK_ButcherTableauClassic}, it is clear that the computation of $\boldsymbol K^{(1)}_2$ and $\boldsymbol K^{(1)}_3$ can only be avoided if $b_2 = b_3 = 0$.

	The presence of zeros in the shared weight vector, however, renders the method non \ac{ssp} as the requirement $\boldsymbol b^{(r)} > \boldsymbol 0, \, r = 1, \dots R$ \cite{higueras2009characterizing, ferracina2005extension, kraaijevanger1991contractivity} (the inequality is to be understood component-wise) is violated for the \ac{PERK} schemes.
	In particular, this \ac{ssp} requirement is always violated for the highest stage method with $E=S$, as for this scheme there are stage evaluations $K_i$ with $b_i = 0$.
	This extends to the other (reducible) methods with $E < S$ as well, with only exception of the base method where the evaluated stages coincide with the non-zeros in the weight vector (which could be \ac{ssp} as a standalone scheme).
	Consequently, the \ac{PERK} schemes offer no guarantees in terms of nonlinear stability properties, and, in fact, we show in the remainder of this section that this leads to problems even in the simplest examples.
	\subsubsection{Linear Advection with Godunov's Scheme}
	\label{subsubsec:LinearAdvectionGodunov}
	To illustrate potential issues with nonlinear stability, we consider the canonical first-order Finite Volume Godunov scheme applied to the 1D advection equation with smooth initial data
	\begin{subequations}
		\label{eq:1DAdvection}
		\begin{align}
			u_t + u_x &= 0, \\
			u_0(t_0 = 0, x) &= 1 + \frac{1}{2} \sin(\pi x)
		\end{align}
	\end{subequations}
	on $\Omega = (-1, 1)$ equipped with periodic boundaries.
	For a potentially non-uniform spatial discretization the semidiscretization \eqref{eq:Semidiscretization} reads for Upwind/Godunov flux
	\begin{subequations}
		\label{eq:UpwindFullSystem}
		\begin{align}
			\renewcommand\arraystretch{1.3}
			\boldsymbol U_0 &= \boldsymbol U(t_0) \\
			\label{eq:UpwindFullSystemODE}
			\frac{\nid }{\nid t} \boldsymbol U(t) &= \underbrace{\begin{pmatrix} -\frac{1}{\Delta x_1} & 0 & \dots & 0 & \frac{1}{\Delta x_1} \\ \frac{1}{\Delta x_2} & -\frac{1}{\Delta x_2} & 0  & \dots & 0\\
					& \ddots & \ddots &  \\ 0 & \dots & 0 & \frac{1}{\Delta x_N} & -\frac{1}{\Delta x_N}\end{pmatrix}}_{\eqqcolon L} \boldsymbol U(t) \, ,
		\end{align}
	\end{subequations}
	where 
	\begin{equation}
		U_i(t) \coloneqq \frac{1}{\Delta x_i}\int_{x_{i-\sfrac{1}{2}}}^{x_{i+\sfrac{1}{2}}} u(t, x) \, \nid x, \quad i = 1, \dots, N
	\end{equation}
	denote the cell-averages.
	Applying a two-level ($R=2$) \ac{PRKM} to \cref{eq:UpwindFullSystemODE} leads to the fully discrete system 
	\begin{equation}
		\label{eq:PartitionedSystem}
		\boldsymbol U_{n+1} = \underbrace{\begin{pmatrix} D^{(1)}(\Delta t, L) \\ D^{(2)}(\Delta t, L) \end{pmatrix}}_{
		\eqqcolon D } 
	\boldsymbol U_{n}
	\end{equation}
	which is equivalent to the matrix-valued stability function \eqref{eq:StabFuncPRK} with $R=2$ and specified mask matrices $I^{(1)}, I^{(2)}$.
	\subsubsection{Linear Stability of the Fully Discrete System}
	\label{subsubsec:LinearStabFullyDiscrete}
	First, we consider an equidistant discretization with $N=64, \Delta x = \Delta x^{(1)} = \Delta x^{(2)} = \sfrac{2}{N}$ where we use the optimized sixteen stage-evaluation method in the center of the domain, i.e., within $[-0.5, 0.5]$, and the optimized eight stage-evaluation method in the remainder of the domain, see \cref{fig:u816_64}.
	For the Godunov scheme, the spectrum $\boldsymbol \sigma\big(J(\boldsymbol U) \big)$ is circular with center at $-a / \Delta x$ and radius $\rho = a / \Delta x$.
	For this spectrum the optimal stability polynomial for second-order accurate methods is known analytically for all polynomial degrees \cite{owren1990some}.
	The Butcher array defining values $a_{i,i-1}^{(r)}$ of the employed $E = \{8, 16\}$ stage methods are provided in compact form in \cref{eq:PERK_8_16_CompactForm}.
	The other non-zero entry $a_{i,1}^{(r)}$ is determined by the internal consistency condition \eqref{eq:InternalConsistency}, i.e., $a_{i,1}^{(r)} = c_i - a_{i,i-1}^{(r)}$.
	\begin{equation}
		\label{eq:PERK_8_16_CompactForm}
		\begin{array}
			{c|c|c|c}
			i & \boldsymbol c &  a_{i, i-1}^{(1)}, E^{(1)} = 8 & a_{i, i-1}^{(2)}, E^{(2)} = 16 \\
			\hline 
			1 & 0 & 0 & 0  \\
			2 & \sfrac{1}{30} & 0 & 0 \\
			3 & \sfrac{2}{30} & 0 & 0.008333333333333335 \\
			4 & \sfrac{3}{30} & 0 & 0.01333333333333334 \\
			5 & \sfrac{4}{30} & 0 & 0.019047619047619042 \\
			6 & \sfrac{5}{30} & 0 & 0.025641025641025637 \\
			7 & \sfrac{6}{30} & 0 & 0.033333333333333354 \\
			8 & \sfrac{7}{30} & 0 & 0.042424242424242434 \\
			9 & \sfrac{8}{30} & 0 & 0.053333333333333295 \\
			10 & \sfrac{9}{30} & 0 & 0.06666666666666667 \\
			11 & \sfrac{10}{30} & 0.019841269841269837 & 0.08333333333333337 \\
			12 & \sfrac{11}{30} & 0.04489795918367346 & 0.10476190476190472 \\
			13 & \sfrac{12}{30} & 0.07792207792207795 & 0.13333333333333336 \\
			14 & \sfrac{13}{30} & 0.12380952380952381 & 0.17333333333333337 \\
			15 & \sfrac{14}{30} & 0.19230769230769232 & 0.23333333333333323 \\
			16 & \sfrac{15}{30} & 0.3061224489795918 & 0.3333333333333334 \\
		\end{array}
	\end{equation}

\begin{figure}[ht]
		\centering
		\includegraphics[width=.45\textwidth]{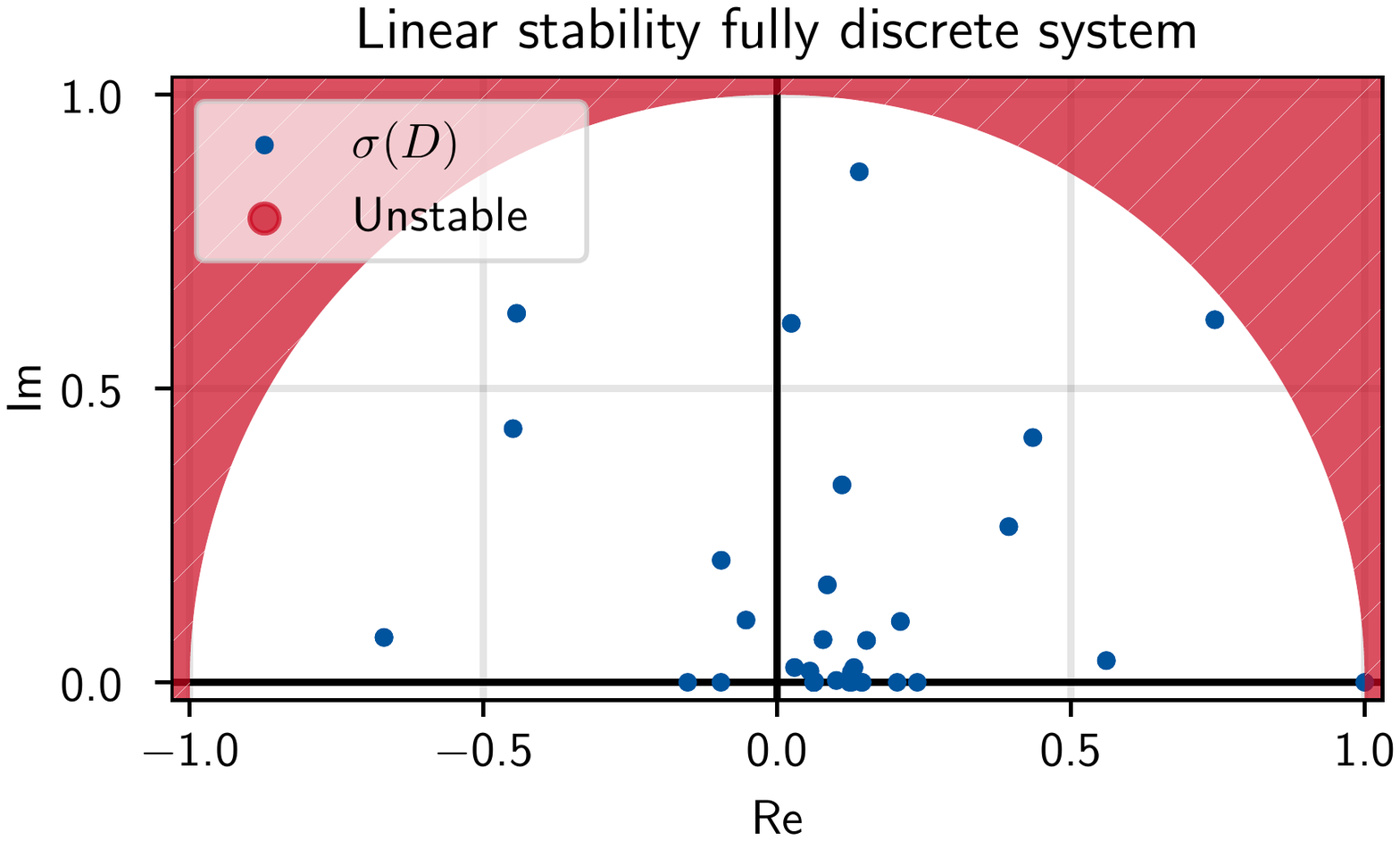}
		\caption[Spectrum $\boldsymbol \sigma(D)$ of the fully discrete matrix $D$ of the $\text{P-ERK}_{2; \{8, 16\}}$ scheme applied to \eqref{eq:UpwindFullSystemODE}.]
		{Spectrum $\boldsymbol \sigma(D)$ of the fully discrete matrix $D$ of the $\text{P-ERK}_{2; \{8, 16\}}$ scheme applied to \eqref{eq:UpwindFullSystemODE}.
		The unstable region with $\vert \lambda \vert > 1$ is shaded in red.
		Note that the purely real eigenvalue around $(1, 0)$ is indeed stable with numeric value $1.0 - 28 \varepsilon$ and machine epsilon $\varepsilon = \mathcal{O}\left(10^{-16}\right)$.
		}
		\label{fig:LinearStability_FullyDiscreteSystemUniform}
	\end{figure}

	The solution $\boldsymbol U_1 = \boldsymbol U(\Delta t)$ after one step with maximum stable timestep $\Delta t = 0.21875$ is presented in \cref{fig:u816_64}.
	Clearly, there are heavy oscillations at the interface between the two schemes.
	It is natural to ask whether these oscillations increase over time and cause a blow-up of the overall simulation.
	To answer this question, we consider the spectrum of $D$ (or equivalently $P(Z)$) which is also displayed in \cref{fig:LinearStability_FullyDiscreteSystemUniform}.
	In the present case the spectral radius $\rho(D) < 1$ and thus the overall scheme is linearly stable \cite[Chapter~8.6]{hespanha2018linear}, which is confirmed by long time runs.
	This shows that the observed oscillations are no consequence of violated linear/absolute stability, but instead are due to the nonlinear stability properties of the scheme which are explored in the next subsection.
	\subsubsection{Nonlinear Stability of the Fully Discrete System}
	\label{subsubsec:NonlinearStabFullyDiscrete}
\begin{figure}[!t]
		\centering
		\subfloat[{Result $\boldsymbol U_1$ of \eqref{eq:PartitionedSystem} after one step.
		The vertical dashed lines indicate the partitioning into $\boldsymbol U^{(1)}$ and $\boldsymbol U^{(2)}$.}]{
			\label{fig:u816_64}
			\centering
			\includegraphics[width=.45\textwidth]{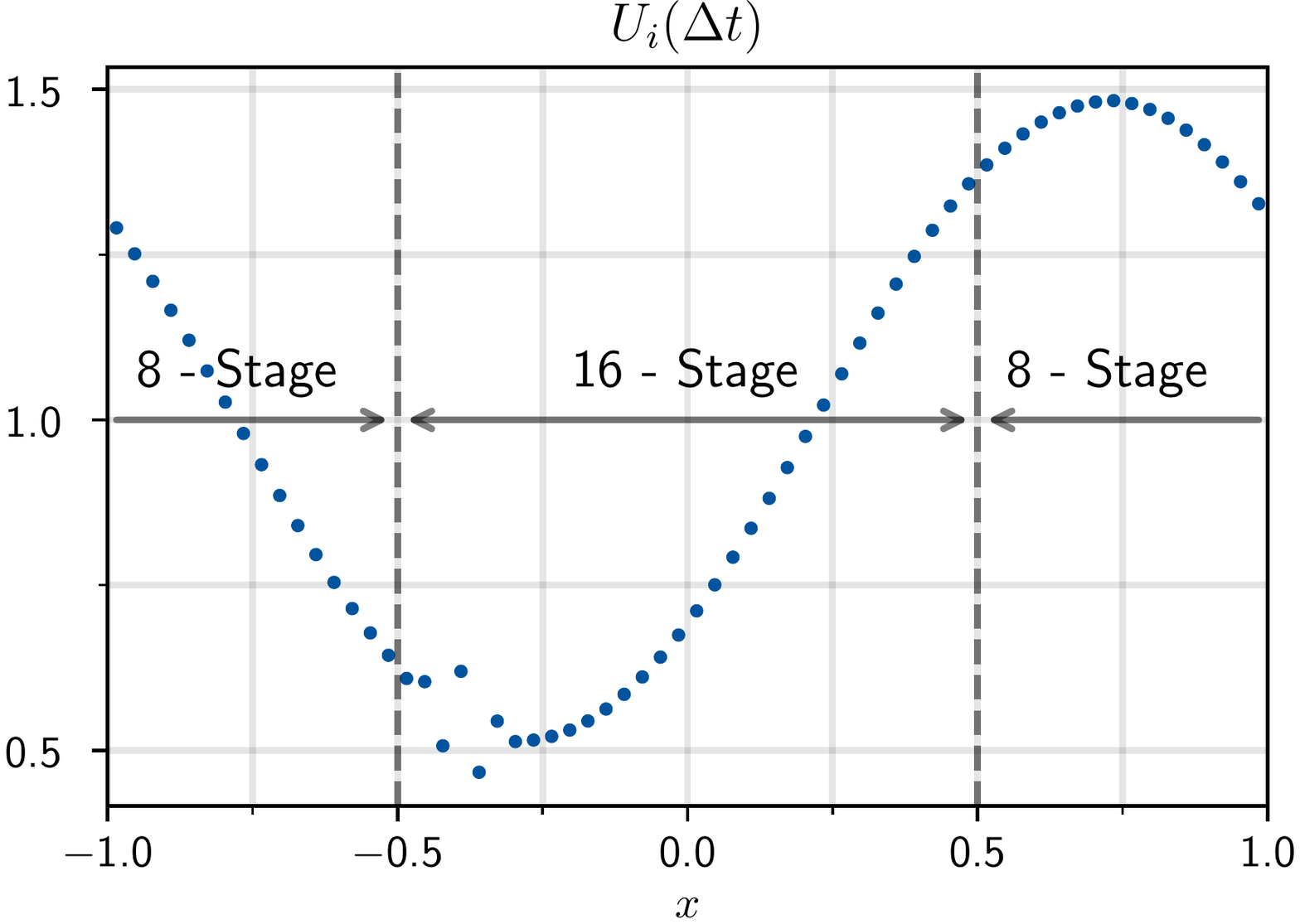}
		}
		\hfill
		\subfloat[{Non-zero entries $D_{ij}$ of fully-discrete system matrix $D$, cf. \eqref{eq:PartitionedSystem}.}]{
			\label{fig:MatrixPlot}
			\centering
			\includegraphics[width=.47\textwidth]{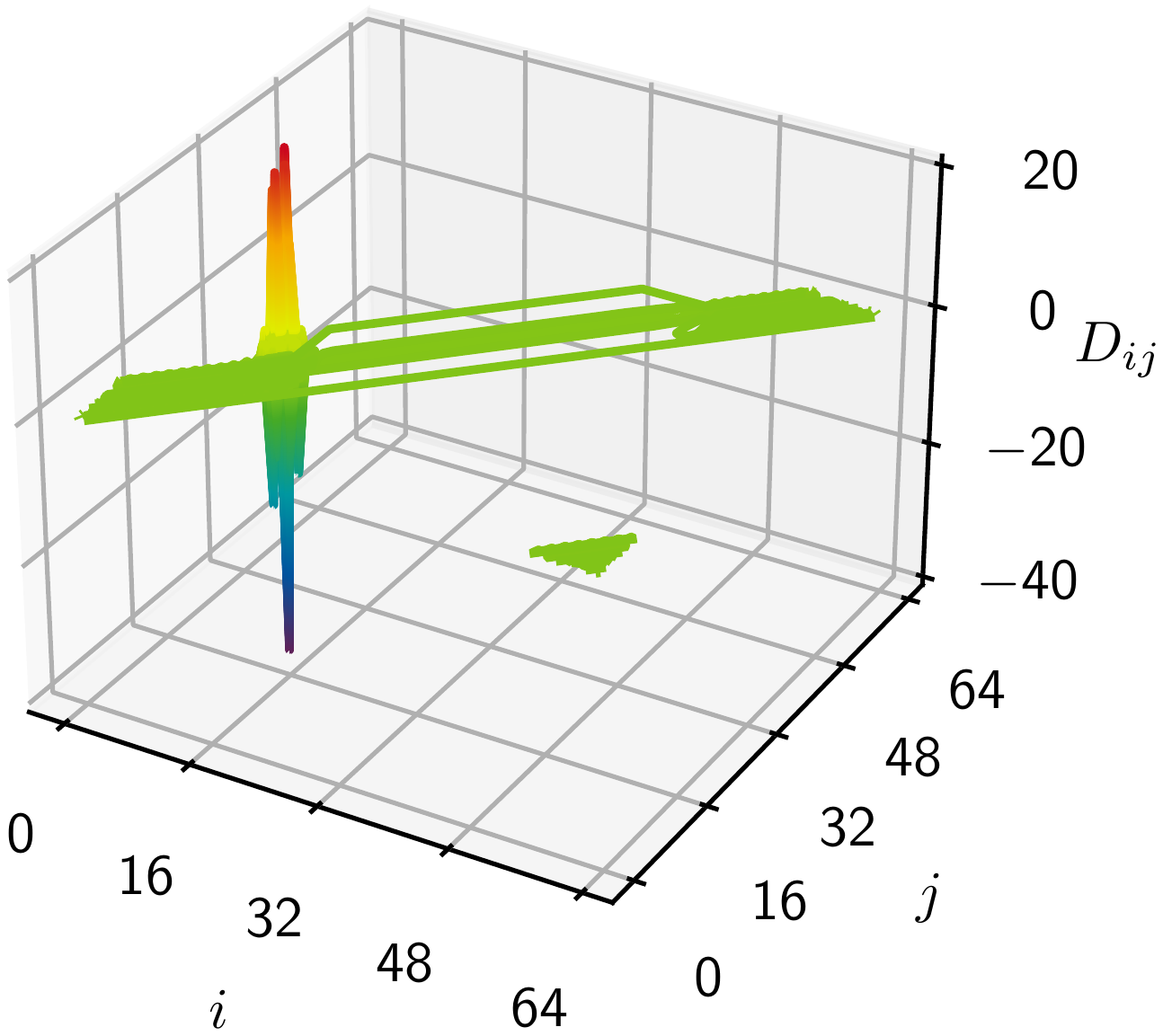}
		}
		\caption[Fully-Discrete Matrix $D$ corresponding to \eqref{eq:UpwindFullSystemODE} with temporal integration through a second-order $8-16$ \ac{PERK} scheme and results $U_i(\Delta t)$ after performing one timestep.]
		{Cell averages $U_i(\Delta t)$ (\cref{fig:u816_64}) after performing one timestep of the second-order 8-16 optimized \ac{PERK} scheme applied to \eqref{eq:UpwindFullSystemODE} on a uniform mesh with $\Delta x = \sfrac{2}{64}$.
		The entries of the fully discrete Matrix $D$ (cf \eqref{eq:PartitionedSystem}) are shown in \cref{fig:MatrixPlot}.}
		\label{fig:Matrix_and_u816_64}
	\end{figure}

	The fully discrete system \eqref{eq:PartitionedSystem} is monotonicity-preserving in the sense of Harten \cite{HARTEN1997260} if and only if
	\begin{equation}
		\label{eq:MonotonicityRequirement}
		D_{ij} \geq 0 \quad \: \forall \: i, j = 1, \dots, N.
	\end{equation}
	The application of the \ac{PERK} schemes leads in general to a violation of \eqref{eq:MonotonicityRequirement}, even for uniform meshes as clearly seen in \cref{fig:u816_64}.
	To shed some more light on this, we present a plot of the matrix entries $D_{ij}$ in \cref{fig:MatrixPlot}.
	It is evident from \cref{fig:MatrixPlot} that monotonicity is lost at the intersection of the schemes, precisely where the higher-stage method retrieves only first-order data from the lower stage scheme.
	first-order information means in this case that the intermediate state $\boldsymbol U_n + \Delta t \sum_{i=1}^S \sum_{r=1}^R b_i^{(r) } \boldsymbol K_i^{(r)}$ (cf. \eqref{eq:PartitionedRKSecondEq}) received on the $r=1$ partition only information from one evaluation of the \ac{RHS}.
	In particular, the rows $i = 17, \dots, 24$ of $D$ have negative entries which originates from the fact that for these points the higher stage method still relies on first-order information only, precisely the previous iterate $\boldsymbol U_n$ plus the first stage $\boldsymbol K_1$.
	In particular, for the $S=16, E = \{8, 16\}$ combination at stage $i=8$ the only update to $\boldsymbol U^{(1)}$ is due to $\left(c_8 - a_{8,1}^{(1)}\right) \boldsymbol K_8^{(1)}$, while $\boldsymbol U^{(2)}$ receives an update including $a_{8,7}^{(2)}\boldsymbol K_7^{(2)}$ which contains seven consecutive evaluations of the \ac{RHS}.
	
	The influence of the negative entries can be clearly seen in \cref{fig:u816_64}, where oscillations arise behind the transition from one scheme to the other.
	Since the initial condition is transported from the left to the right, the usage of the upwind fluxes corresponds to using information from the left.
	This is why there are oscillations in this case at $x=-0.5$, since the higher stage-evaluation method receives only first-order information based upon $\boldsymbol K_1^{(1)}$.
	At $x= 0.5$ there are for positive advection velocity $a=1$ no oscillations, as the lower stage-evaluation method retrieves high order information at any time.
	In turn, for a negative transportation velocity, i.e., $u_t -u_x=0$ the oscillations arise at $x= 0.5$.

	We mention that each row sum $d_i \coloneqq \sum_{j=1}^N D_{ij}$ equals one, reflecting that the scheme is conservative, agreeing with the criteria \eqref{eq:Conservation}.
	Furthermore, the \ac{PERK} scheme converges with second-order once the asymptotic regime is reached. 
	
	This loss of monotonicity is noteworthy as it foreshadows possible problems for more complex scenarios.
	The first-order Godunov finite volume scheme is known to be very diffusive, yet already for this system and uniform meshes, artificial oscillations are introduced.
	There are several ways how the oscillations can be mitigated. 
	First of all, the oscillations become less severe if the timestep is reduced since this leads to a reduction in magnitude of the negative entries $D_{ij}$.
	Next, for a higher spatial resolution $N$ the influence of the negative entries is reduced as neighboring values $U_i$ are of more similar value (note that no oscillations can occur for a constant solution since $d_i \coloneqq \sum_{j=1}^N D_{ij} = 1$).
	Finally, the oscillations can also be reduced if the difference in the number of stage evaluations $E$ is reduced, with monotonicity recovered for the same method in each part of the domain.
	These observations are quantified in \cref{subsubsec:NonlinearStabPERK_NonUniform} for the intended application of the \ac{PERK} schemes to non-uniform grids.
	\subsubsection{Alternative P-ERK Schemes}
	It is also instructive to consider a different stage evaluation pattern of the \ac{PERK} methods.
	The original \ac{PERK} methods perform all shared evaluations at the end, cf. \eqref{eq:PERK_ButcherTableauClassic}.
	In principle, also different stage evaluation patterns are possible while preserving the key features of the \ac{PERK} schemes.
	For instance, one can construct ''alternating'' schemes such as 
	\begin{equation}
		\label{eq:PERK_ButcherTableauAlternating}
		\renewcommand\arraystretch{1.3}
		\begin{array}
			{c|c|c c c c c c}
			i & \boldsymbol c & & & A^{(1)} & & \\
			\hline 
			1 & 0 & & & & & &    \\
			2 & \sfrac{1}{10} & \sfrac{1}{10} & & & & &   \\
			3 & \sfrac{2}{10} & \sfrac{2}{10} & 0 & & & & \\ 
			4 & \sfrac{3}{10} & \sfrac{3}{10} -a_{4,3}^{(1)} & 0 & a_{4,3}^{(1)} & & & \\ 
			5 & \sfrac{4}{10} & \sfrac{4}{10} & 0 & 0  & 0 & & \\
			6 & \sfrac{5}{10} & \sfrac{5}{10} - a_{6,4}^{(1)} & 0 & 0 & a_{6,4}^{(1)} & 0 & \\
			\hline
			& \boldsymbol b^T & 0 & 0 & 0 & 0 & 0 & 1
		\end{array}
	\end{equation}
	from the same stability polynomial $P_{2;4}^{(1)}(z)$ \eqref{eq:PERK_StabPolys} that would be used for the traditional \ac{PERK} scheme \eqref{eq:PERK_ButcherTableauClassic}.
	For the first scheme, specified according to \eqref{eq:PERK_ButcherTableauAlternating}, the stages $\boldsymbol K_1^{(1)}, \boldsymbol K_3^{(1)}, \boldsymbol K_4^{(1)}, \boldsymbol K_6^{(1)}$ are evaluated, while for the classic \ac{PERK} scheme \eqref{eq:PERK_ButcherTableauClassic} $\boldsymbol K_1^{(1)}, \boldsymbol K_4^{(1)}, \boldsymbol K_5^{(1)}, \boldsymbol K_6^{(1)}$ are evaluated.

	Another possible design pattern for \ac{PERK} schemes is to perform the shared evaluations at the early stages of the method, rather than at the end as for the traditional \ac{PERK} schemes.
	The Butcher Tableau for the $p=2$, $E=4, S=6$ scheme with this stage evaluation pattern is given by
	\begin{equation}
		\label{eq:PERK_ButcherTableauEarlyShared}
		\renewcommand\arraystretch{1.3}
		\begin{array}
			{c|c|c c c c c c}
			i & \boldsymbol c & & A^{(1)} & & & \\
			\hline 
			1 & 0 & & & & & &    \\
			2 & \sfrac{1}{10} & \sfrac{1}{10} & & & & &   \\
			3 & \sfrac{2}{10} & \sfrac{2}{10} - a_{3,2}^{(1)} & a_{3,2}^{(1)} & & & & \\ 
			4 & \sfrac{3}{10} & \sfrac{3}{10} & 0 & 0 & & &\\ 
			5 & \sfrac{4}{10} & \sfrac{4}{10} & 0  & 0 & 0 & & \\
			6 & \sfrac{5}{10} & \sfrac{5}{10} - a_{6,3}^{(1)} & 0 & a_{6,3}^{(1)} & 0 & 0 &\\
			\hline
			& \boldsymbol b^T & 0 & 0 & 0 & 0 & 0 & 1
		\end{array}
	\end{equation}
	where the free coefficients $a_{3,2}^{(1)}, a_{6,3}^{(1)}$ are again obtained from the identical stability polynomial $P_{2;4}^{(1)}(z)$ that may be used for the other methods.
	In this case, stages $\boldsymbol K_1^{(1)}, \boldsymbol K_2^{(1)}, \boldsymbol K_3^{(1)},\boldsymbol K_6^{(1)}$ are evaluated.

\begin{figure}[!t]
		\centering{
			\includegraphics[width=.45\textwidth]{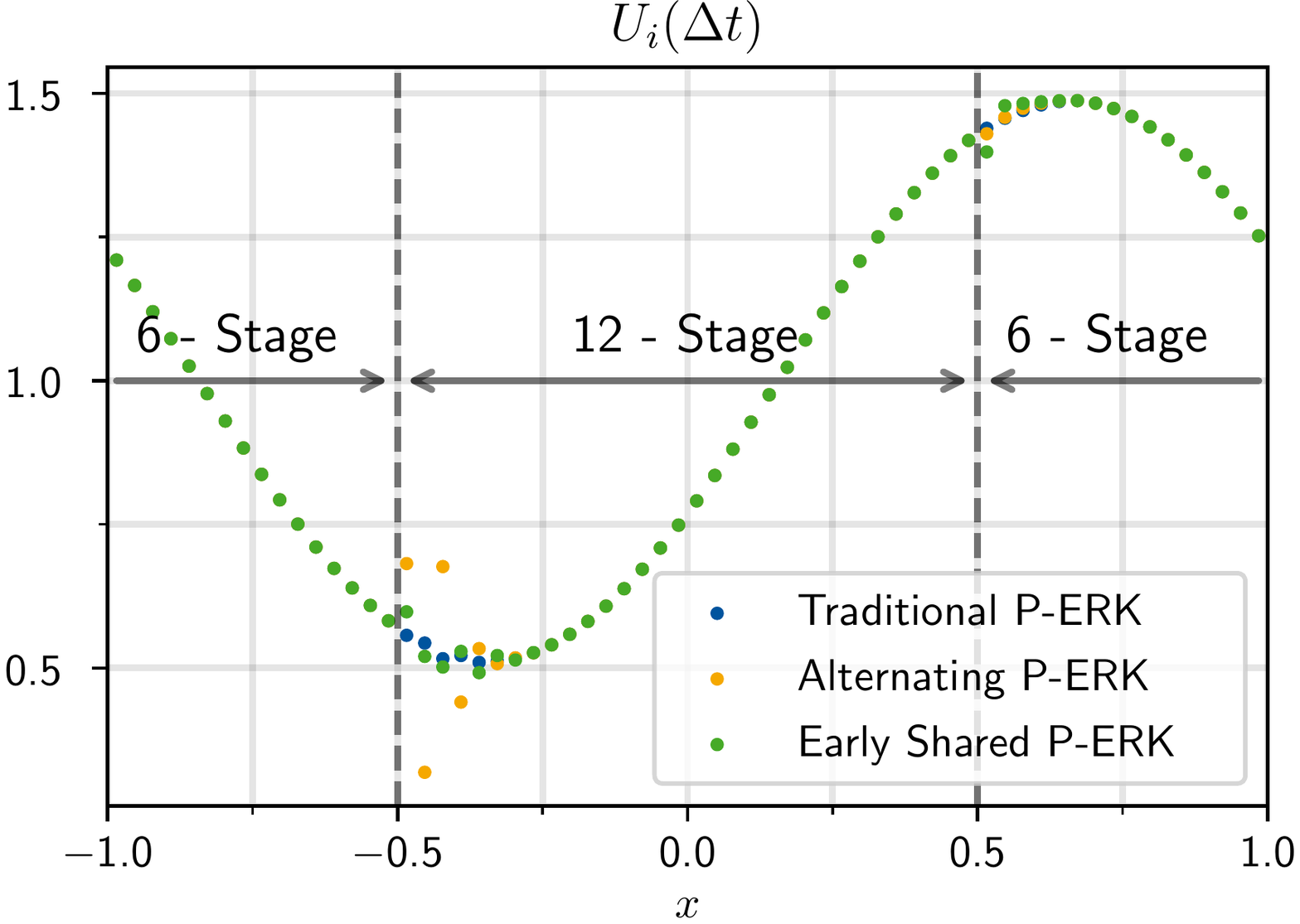}
		}
		\caption[Result $\boldsymbol U^1$ to \cref{eq:PartitionedSystem} after one timestep of two different 6-12 \ac{PERK} schemes with different stage evaluation patterns.]
		{Result $\boldsymbol U^1$ to \cref{eq:PartitionedSystem} after one timestep of two different 6-12 \ac{PERK} schemes with different stage evaluation patterns on a uniform mesh.
		The solution colored in blue correspond to the standard \ac{PERK} scheme \eqref{eq:PERK_ButcherTableauClassic} while the solution obtained with the alternating \ac{PERK} scheme \eqref{eq:PERK_ButcherTableauAlternating} is colored in orange.}
		\label{fig:uPlot_Comparison}
	\end{figure}

	To illustrate the influence of the stage evaluation pattern \cref{fig:uPlot_Comparison} compares the results after one timestep for the standard \cref{eq:PERK_ButcherTableauClassic}, alternating \cref{eq:PERK_ButcherTableauAlternating}, and shared early stage \cref{eq:PERK_ButcherTableauEarlyShared} \ac{PERK} scheme.
	Clearly, the schemes yield different results at the interfaces at the schemes, with the strongest oscillations observed for the alternating scheme.
	Notably, the scheme with shared early stage evaluations causes the largest errors at the fine-coarse interface at $x=0.5$.

	Analogous to the linear stability properties discussed in \cref{subsec:AbsoluteStability}, the action of partitioned schemes is no longer determined by the stability polynomials of the individual schemes only, even for linear \acp{ode}.
	In fact, due to the partitioning there is a discontinuity in the time-integration method introduced and consequently, the final outcome does depend on the design of the stage evaluation pattern.

	This feature is similar in spirit to \ac{ssp} time integration schemes.
	Consider for instance the two-stage second-order stability polynomial, for which infinitely many Runge-Kutta methods exist, for instance the midpoint method and Heun's method mentioned earlier.
	When applied to a linear problem, both give the identical result, as the overall outcome is governed by the same stability polynomial.
	In contrast, when applied to a nonlinear, potentially discontinuous problem, the implementation of the intermediate stage is of crucial importance as this determines the nonlinear stability properties, see the instructive example given in \cite{gottlieb1998total}.
	\subsubsection{Non-uniform Meshes}
	\label{subsubsec:NonlinearStabPERK_NonUniform}
	Coming to the intended application case of the \ac{PERK} schemes to non-uniform meshes, it is observed that the presence of a discontinuity in the system matrix $L$ \eqref{eq:UpwindFullSystemODE} further impairs nonlinear stability.
	To give evidence for this, we consider again the model problem \eqref{eq:1DAdvection} and refine the grid in $\Omega^{(2)} \coloneqq [-0.5, 0.5]$ by a factor $\alpha > 1$, i.e., $\Delta x_2 = \Delta x_1/\alpha$.
	In the non-refined part of the mesh with $\Delta x_1 = \sfrac{2}{64}$ an $E^{(1)} = 8$ stage-evaluation method is used, while in the refined part the number of stage evaluations is determined based on $\alpha$ such that the timestep $\Delta t_8 = 0.21875$ can be kept constant.
	For all configurations presented in this section linear stability was checked by means of the spectral radius of the fully discrete system matrix $D$ \eqref{eq:PartitionedSystem} and long-time simulations.

	To measure the oscillations, we consider the relative total variation (TV) increase after one timestep
	\begin{equation}
		e^{\text{TV}} \coloneqq \frac{\Vert \boldsymbol 
		U_1
		\Vert_\text{TV} - \Vert \boldsymbol U_0 \Vert_\text{TV}}{\Vert \boldsymbol U_0 \Vert_\text{TV}} 
	\end{equation}
	with total variation seminorm $\Vert \boldsymbol U \Vert_\text{TV} \coloneqq \sum_i \vert U_{i+1} - U_i \vert$.
	$e^{\text{TV}}$ is reported in \cref{tab:TotalVariationIncrease_AlphaVariation} for a range of stage evaluations $E^{(2)}$ in the refined region.
	Only for the most modestly refined grid the total variation actually decreases, preserving the diffusive nature of the Godunov scheme.
	For higher refinement ratios, exponential increase in total variation is observed, with dramatic values for the highest refinements.
	\begin{table}[!ht]
		\def\arraystretch{1.3}
		\centering
		\begin{tabular}{c?{1.5}c|c|c|c|c|c|c|c}
			$E^{(2)}$ & 9 & 10 & 11 & 12 & 13 & 14 & 15 & 16 \\
			\hline
			$\alpha$ & 1.125 & 1.25 & 1.375 & 1.5 & 1.625 & 1.75 & 1.875 & 2.0 \\
			\Xhline{3.5\arrayrulewidth}
			$e^{\text{TV}}$ & $
			\boldsymbol{-} 0.03$ & $0.11$ & $0.55$ & $1.65$ & $3.71$ & $7.85$ & $15.4$ & $26.0$
		\end{tabular}
		\caption[Increase in relative total variation $e_{\, \text{TV}}$ for non-uniform meshes with grid sizes $\Delta x, \frac{\Delta x}{\alpha}$.]{Increase in relative total variation $e_{\, \text{TV}}$ for non-uniform meshes with grid sizes $\Delta x, \frac{\Delta x}{\alpha}$.
		The $E^{(1)} = 8$ method is used in the non-refined region for all $\alpha$.
		}
		\label{tab:TotalVariationIncrease_AlphaVariation}
	\end{table}
	\\
	Similar to the case for uniform grid with partitioned schemes, the oscillations decrease with a reduction in timestep, increase in spatial resolution, and smaller difference in the number of stage evaluations $E^{(2)} - E^{(1)}$.
	To illustrate this, we reduce the maximum stable timestep $\Delta t_8$ by a factor $\text{CFL} \leq 1$ and provide the increase in total variation in \cref{tab:TotalVariationIncrease_dtVariation}
	\begin{table}[!ht]
		\def\arraystretch{1.3}
		\centering
		\begin{tabular}{c?{1.5}c|c|c|c|c|c|c}
			CFL & $0.4$ & $0.5$ & $0.6$ & $0.7$ & $0.8$ & $0.9$ & $1.0 $ \\
			\Xhline{3.5\arrayrulewidth}
			$e^{\text{TV}}$ & $-0.01$ & $0.03$ & $0.28$ & $1.21$ & $4.04$ & $11.5$ & $26.0$
		\end{tabular}
		\caption[Increase in relative total variation $e_{\, \text{TV}}$ for non-uniform meshes with varying base grid size.]{Increase in relative total variation $e_{\, \text{TV}}$ for non-uniform meshes with varying base grid size $\Delta x^{(1)} = \sfrac{2}{64}$ and two times ($\alpha = 2$) finer grid in $\Omega^{(2)} = [-0.5, 0.5]$ under reduction of timestep $\Delta t = \text{CFL} \cdot \Delta t_8$.
		The coarse cells are integrated with the $E^{(1)} = 8$ stage evaluation method, while the fine cells are integrated using the $E^{(2)}=16$ method.
		}
		\label{tab:TotalVariationIncrease_dtVariation}
	\end{table}
	\\
	Next, the base resolution $N^{(1)}$, i.e., the resolution in $\Omega^{(1)} \coloneqq (-1, -0.5] \cup [0.5, 1)$ is varied from $N^{(1)}=64$ to $N^{(1)}=4095$ for the $E^{(1)}=8, E^{(2)}=16$ scheme.
	As tabulated in \cref{tab:TotalVariationIncrease_NVariation} the increased resolution causes the increase in total variation to diminish on finer grids.
	\begin{table}[!ht]
		\def\arraystretch{1.3}
		\centering
		\begin{tabular}{c?{1.5}c|c|c|c|c|c|c}
			$N^{(1)}$ & 64 & 128 & 256 & 512 & 1024 & 2048 & 4096 \\
			\Xhline{3.5\arrayrulewidth}
			$e^{\text{TV}}$ & $26.0$ & $7.27$ & $1.82$ & $0.45$ & $0.11$ & $0.02$ & $0.01$ 
		\end{tabular}
		\caption[Increase in relative total variation $e^{\text{TV}}$ for non-uniform meshes with varying base grid size.]{Increase in relative total variation $e_{\, \text{TV}}$ for non-uniform meshes with varying base grid size $\Delta x_\text{Base} = \sfrac{2}{N^{(1)}}$ and two times ($\alpha = 2$) finer grid in $\Omega^{(2)} = [-0.5, 0.5]$.
		The coarse cells are integrated with the $E^{(1)} = 8$ stage evaluation method, while the fine cells are integrated using the $E^{(2)}=16$ method.
		}
		\label{tab:TotalVariationIncrease_NVariation}
	\end{table}
	\\
	Finally, keeping cell size $\Delta x^{(1)} = \frac{2}{64}$ and ratio $\alpha=2$, the difference in stage evaluations
	$\Delta E= E^{(2)} - E^{(1)}$ is varied.
	Here, we always set $E^{(2)} = 2 E^{(1)}$ which implies that $\Delta E= E^{(2)} - E^{(1)} = 2 E^{(1)} - E^{(1)} = E^{(1)}$.
	For smaller $\Delta E$ the increase in total variation is reduced, as shown in \cref{tab:TotalVariationIncrease_EVariation}.
	This can be intuitively explained by the fact that for smaller $\Delta E$ the higher stage method does less stages with only first-order information based on $\boldsymbol K_1^{(1)}$ from the lower stage method.
	\begin{table}[!ht]
		\def\arraystretch{1.3}
		\centering
		\begin{tabular}{c?{1.5}c|c|c|c|c|c|c}
			$\Delta E = E^{(1)}$ & 2 & 3 & 4 & 5 & 6 & 7 & 8 \\
			\Xhline{3.5\arrayrulewidth}
			$e^{\text{TV}}$ & $-0.00$ & $-0.01$ & $-0.00$ & $0.06$ & $0.61$ & $4.25$ & $26.0$ 
		\end{tabular}
		\caption{Increase in relative total variation $e^{\text{TV}}$ for non-uniform meshes with constant base grid size $\Delta x^{(1)} = \sfrac{2}{64}$ and two times ($\alpha = 2$) finer mesh in $\Omega^{(2)} = [-0.5, 0.5]$ under variation of the number of base stage evaluations $E^{(1)} = \{ 2, \allowbreak 3, \allowbreak 4, \allowbreak \dots, \allowbreak 8\}$.
		The fine cells are integrated using the $E^{(2)}=2 E^{(1)}$ stage-evaluation method.
		}
		\label{tab:TotalVariationIncrease_EVariation}
	\end{table}

	Since these issues are encountered for the Godunov scheme which evolves the cell averages, it is immediately clear that conventional limiting procedures are not applicable to remedy the loss of monotonicity.
	In particular, derivates of standard procedures such as MUSCL-type limiters \cite[Chapter~6]{leveque2002finite} or related techniques to enforce, e.g., positivity \cite{doi:10.1098/rspa.2011.0153} do not work as the mean state is already corrupted and therefore cannot be used as a fallback solution.
	In contrast, one would need to devise some form of flux limiting (or in this case stage-limiting) similar to techniques employed for the Lax-Wendroff scheme \cite[Chapter~6.3]{leveque2002finite}[Chapter~III.1]\cite{godlewski2013numerical}, which also evolves cell averages.
	\subsubsection{Relation of Nonlinear Stability Analysis to Published Results}
	We would like to classify the previous observations in the context of the published results on \ac{PERK} schemes.
	In \cite{vermeire2019paired, nasab2022third} the \ac{PERK} schemes have been successfully applied to the simulation of the compressible Navier-Stokes equations for applications such as flows around airfoils.
	The meshes employed in the aforementioned works are characterized by a high resolution and small refinement rations between neighboring elements.
	For instance, the mesh for the simulation of a turbulent flow around the SD7003 airfoil employed in \cite{nasab2022third} has an average refinement ratio of $\bar \alpha \approx 1.15$.
	As illustrated by \cref{tab:TotalVariationIncrease_AlphaVariation} moderate refinement ratios are much less vulnerable to monotonicity violations.
	As the refinement ration of neighboring cells is typically small, the difference in stage evaluations $E^{(r+1)} - E^{(r)}$ may also be relatively small.
	In fact, in \cite{vermeire2019paired, nasab2022third} the authors use \ac{PERK} families with a maximum difference in stage evaluations of $\Delta E = 4$.
	Again, we demonstrated by means of \cref{tab:TotalVariationIncrease_EVariation} that the monotonicity violations are reduced for smaller $\Delta E$.
	In conjunction with the fact that the simulated equations bear a diffusive nature, the potential loss of monotonicity is not observed in the aforementioned works.

	Besides presenting application to the compressible Navier-Stokes equations, the inviscid isentropic vortex advection test case \cite{shu1988efficient, wang2013high} is also studied in \cite{vermeire2019paired, nasab2022third} for the purpose of validating the \ac{PERK} schemes.
	Notably, they employ a random distribution of \ac{PERK} schemes on a uniform grid to test the order of convergence of the schemes.
	For this test case it is in principle possible that the $E^{(R)} = 16$ scheme is placed next to an $E^{(1)} = 2$ (for the second-order accurate case) or $E^{(1)} = 3$ (for the third-order accurate case) scheme.
	Based on the findings in the previous section one might expect that this could reveal the loss of monotonicity described in the previous section.
	This is not the case, however, as the timestep has to be adapted to the scheme with the smallest number of stage evaluations which requires a significant reduction in timestep.
	For instance, the maximum stable timestep for the $p=3, E=3$ scheme is about 8 times smaller than the one for the $p=3, E=16$ scheme.
	As tabulated in \cref{tab:TotalVariationIncrease_dtVariation} there are no artificial oscillations observed if the timestep is chosen small enough.

	To further strengthen the assessment discussed above, we created a test case based on a structured, smoothly refined, high-resolution mesh where we use a second-order, 16-stage \ac{PERK} scheme with stage evaluations $E = \{4, 6, 8, \dots, 16\}$.
	In agreement with the analysis above, this testcase can be run without reduction of the timestep for increased stage evaluations $E^{(R)}$, i.e., the stage normalized timestep $\Delta t / E^{(R)}$ can be kept almost constant for all schemes.
	This application is presented in more detail in \ref{sec:AppSmoothHighResMesh}.
	However, we would like to stress that the use case considered in this work is different from the one considered in \cite{vermeire2019paired, nasab2022third} and in \ref{sec:AppSmoothHighResMesh}, as for quad/octree \ac{AMR} the mesh refinement ratio is always a factor of two and thus the difference in stage evaluations $E^{(r)} - E^{(r-1)}$ grows exponentially with the number of refinement levels.
	\subsection{Higher-Order P-ERK Methods}
	In the original paper \cite{vermeire2019paired} \ac{PERK} schemes of order $p=2$ are presented which are extended to order $p=3$ in \cite{nasab2022third}.
	A particular feature of Runge-Kutta methods is that, given a certain stability polynomial, infinitely many Runge-Kutta methods with corresponding Butcher arrays can be constructed.
	For the third-order \ac{PERK} methods, the proposed archetype of the Butcher tableau is \cite{nasab2022third}
	\begin{equation}
		\label{eq:PERK_ThirdOrderProposed}
		\renewcommand\arraystretch{1.2}
		\begin{array}
			{c|c c c}
			0 & & \\
			\frac{1}{3} & \frac{1}{3} & \\
			1 & -1 & 2 \\
			\hline 
			 & 0 & \frac{3}{4} & \frac{1}{4}
		\end{array}
	\end{equation}
	which corresponds to the third-order Runge-Kutta method presented in \cite{sanderse2019constraint} with parameter $c_2= \frac{1}{3}$.
	The presence of the negative entry $a_{3,1} = -1$ renders the incorporation of the $E=3, p=3$ in a \ac{PERK} ensemble unattractive as the negative entry corresponds essentially to downwinding the numerical fluxes of stage $\boldsymbol K_1$ \cite{gottlieb1998total}.
	Alternatively, we choose the Butcher array of the canonical three-stage, third-order method by Shu and Osher \cite{shu1988efficient} as the base scheme for our $p=3$ \ac{PERK} schemes, i.e., 
	\begin{equation}
		\label{eq:PERK_ThirdOrderShuOsher}
		\renewcommand\arraystretch{1.2}
		\begin{array}
			{c|c c c}
			0 & & \\
			1 & 1 & \\
			\frac{1}{2} & \frac{1}{4} & \frac{1}{4} \\
			\hline 
			& \frac{1}{6} & \frac{1}{6} & \frac{4}{6}
		\end{array}
	\end{equation}
	We emphasize that the usage of the first Runge-Kutta stage $\boldsymbol K_1$ in the final update step does not increase storage requirements of the scheme as $\boldsymbol K_1$ is used anyways to compute the higher stages, cf. \eqref{eq:PERK_ButcherTableauClassic}.
	Also, the computational costs per stage stay at two scalar-vector multiplications and one vector-vector addition.
	
	We have not yet discussed how to choose the free abscissae $c_i, i = 2, \dots, S-2$. 
	In this work, we set
	\begin{equation}
		\label{eq:Free_Timesteps}
		c_i = \frac{i-1}{S-3}, \quad i = 1, \dots, S-2
	\end{equation}
	which corresponds to a linear distribution of timesteps between $c_1=0$ and $c_{S-2} = 1$, similar to the second-order case \cite{vermeire2019paired}.
	We mention that more sophisticated choices of the abscissae are in principle possible, e.g., to improve the internal stability properties of the scheme \cite{ketcheson2014internal}.
	Since this is not easily done and becomes also more relevant for many-stage methods which are not the main focus of this work, we restrict ourselves in this work to the linear distribution according to \eqref{eq:Free_Timesteps}.
	Additionally, in contrast to the second-order accurate method, a nonlinear solver is required to obtain the $a_{i,i-1}$ from the coefficients $\alpha_j$ of the stability polynomial.
	Although in principle not guaranteed, with the choice of abscissae 
	\eqref{eq:Free_Timesteps} we obtain positive Butcher coefficients $a_{i,i-1}$ for all optimized stability polynomials considered in this work.
	
	As mentioned in \cref{sec:OrderConditions} the construction of a fourth order \ac{PERK} scheme is significantly more complex due to the coupling conditions involving the different Butcher arrays $A^{(r)}$ and is left for future research.
	\section{Construction of the Stability Polynomials}
	\label{sec:ConstructionStabPoly}
	Besides their multirate property an additional advantage of the \ac{PERK} schemes is that the stability polynomials of the constituting schemes can be optimized according to the intended application to maximize the admissible timestep $\Delta t$.
	In particular, we can construct stability polynomials that are tailored to the problem at hand, i.e., the spatial discretization, boundary and initial conditions.
	The optimization problem for the $p$'th order, $E^{(r)}$ degree stability polynomial thus reads
	\begin{equation}
		\label{eq:OptimizationProblem}
		\max \Delta t \quad \text{s.t. } \left \vert P^{(r)}_{p;E^{(r)}}(\Delta t \lambda) \right \vert \leq 1, \: \forall \: \lambda \in \boldsymbol \sigma \, .
	\end{equation}
	To solve this problem we resort to the approach developed in \cite{ketcheson2013optimal} which has also been applied in \cite{vermeire2019paired, nasab2022third}.
	In the following, we first describe how the relevant eigenvalues of the spectrum $\boldsymbol \sigma$ are obtained and then discuss the solution of the optimization problem \eqref{eq:OptimizationProblem}.
	\subsection{Generating Jacobian and Spectrum}
	\label{sec:SpectraGeneration}
\begin{figure}[ht]
		\centering
		\subfloat[{Reduced spectrum obtained from the procedure described in \cref{sec:SpectraGeneration}.}]{
			\label{fig:ReducedSpectrum}
			\centering
			\includegraphics[width=.45\textwidth]{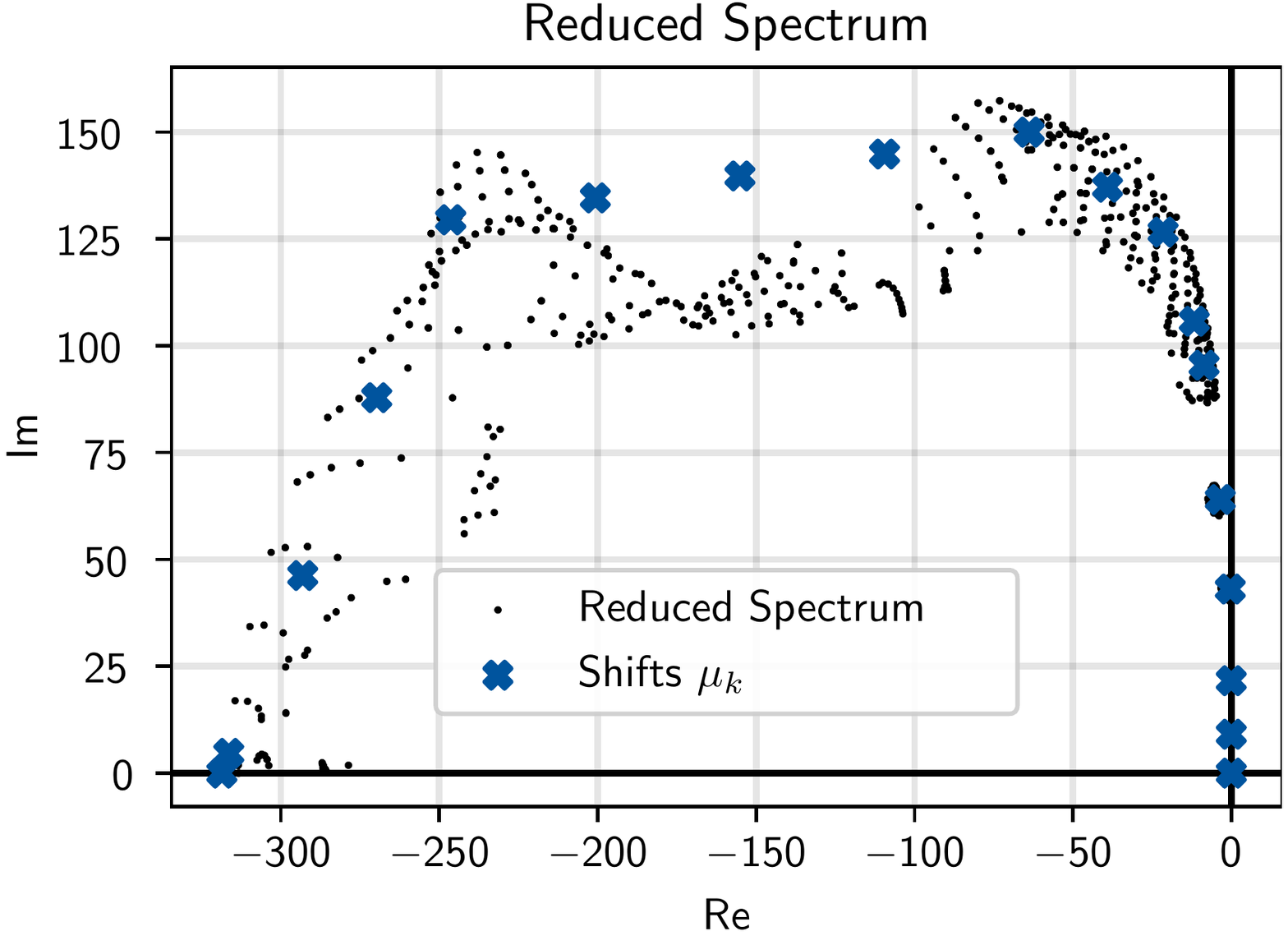}
		}
		\hfill
		\subfloat[{Full spectrum obtained from standard eigendecomposition.}]{
			\label{fig:FullSpectrum}
			\centering
			\includegraphics[width=.45\textwidth]{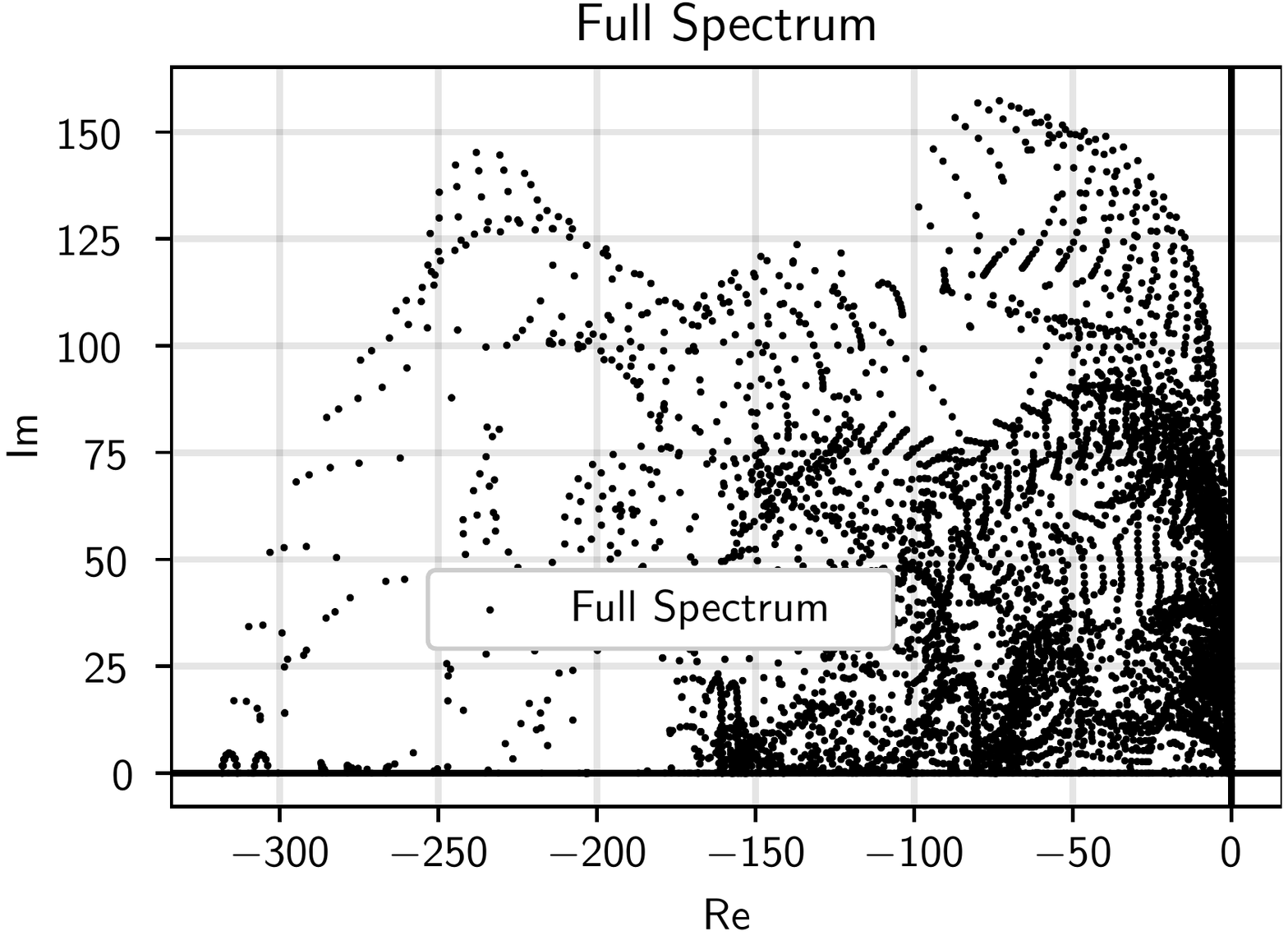}
		}
		\caption[Reduced spectrum obtained from Arnoldi iterations around shifts and full spectrum obtained from standard eigendecomposition.]
		{Reduced spectrum obtained from Arnoldi iterations around shifts (\cref{fig:ReducedSpectrum}) and full spectrum obtained from standard eigendecomposition (\cref{fig:FullSpectrum}).
		The spectrum corresponds to a \ac{dgsem} discretization of the compressible Euler equations with an initial condition leading to a Kelvin-Helmholtz instability.
		}
		\label{fig:SpectraComparisonReducedFull}
	\end{figure}

	Throughout this work we use \texttt{Trixi.jl} \cite{trixi1, trixi2, trixi3}, a high-order \ac{DG} code written in \texttt{Julia} for the numerical simulation of conservation laws.
	\texttt{Trixi.jl} supports algorithmic differentiation (AD) which allows exact and efficient computation of the Jacobian $J(\boldsymbol U)$ \eqref{eq:Jacobian}.
	As long as memory requirements are no concern, the Jacobian can be computed quickly even for relatively large \ac{ode} systems \eqref{eq:Semidiscretization}.
	In contrast, the eigenvalue decomposition to obtain the spectrum $\boldsymbol \sigma\big(J(\boldsymbol U) \big)$ is computationally infeasible for practical simulation setups with millions of \ac{DoFs}.
	At the same time, for the optimization of a stability polynomial the entire spectrum is not required.
	Instead, it suffices to supply a subset of the eigenvalues forming the boundary of the spectrum \cite{trixi2, doehring2024manystage}.
	To obtain this subset of eigenvalues lying on the boundary of the spectrum efficiently, we employ the following routine:
	\begin{enumerate}
		\item Compute $\widetilde{\boldsymbol \sigma}\left( \widetilde{J}(\boldsymbol U_0) \right)$ for a semidiscretization $\boldsymbol U' = \widetilde{\boldsymbol F}(\boldsymbol U)$ with reduced minimal spatial resolution $\widetilde{h}$ such that the complete eigenvalue decomposition is still computational feasible.
		In this work, we use reduced systems with dimensionality $\mathcal O\left(10^3 \times 10^3 \right)$.
		\item Compute the convex hull $\widetilde{\boldsymbol \mu}\left(\widetilde{\boldsymbol \sigma}\right)$ of the reduced spectrum $\widetilde{\boldsymbol \sigma}$.
		As the convex hull points can be quite distant from each other, additional points on the polygons are added such that the distance between two points is below the average distance $\delta$ between hull-points $ \delta \coloneqq \sum_{k=1}^{K-1} \frac{\Vert \tilde{\mu}_{k+1} - \tilde{\mu}_{k} \Vert_2}{K-1}$.
		
		It suffices to consider either the second or third quadrant of the complex plane due to symmetry of $\boldsymbol \sigma$ with respect to the real axis.
		This also assumes that there are no amplifying modes with $\text{Re}(\lambda) >0$ in the spectrum, which is the case for the semidiscretizations arising from conservation laws considered in this work.
		\item Increase the spatial resolution to the actual value used in the simulation, or the maximum value for which $J(\boldsymbol U_0)$ can still be stored in main memory.
		Convert $J(\boldsymbol U_0)$ to sparse format to save memory and computational costs for the subsequent Arnoldi iterations.
		For the \ac{DG} discretizations considered here, $J(\boldsymbol U_0)$ is typically extremely sparse due to the small \ac{DG} stencils.
		As a side-node, ideally one would compute $J(\boldsymbol U_0)$ already in sparse format. 
		This, however, requires special techniques to identify the nonzeros (matrix-coloring) in the Jacobian which are not yet supported in \texttt{Trixi.jl}.
		\item As discussed in \cref{sec:PERK_Intro} the spectrum of convection-dominated \acp{PDE} scales increases linearly with decreasing characteristic grid size $h$.
		Consequently, the convex hull $\widetilde{\boldsymbol \mu} \in \mathbb C^K$ of the reduced spectrum is scaled by $\sfrac{\widetilde{h}}{h}$ to agree with the discretization $\boldsymbol U' = \boldsymbol F(\boldsymbol U)$.
		
		Then, the `outer` eigenvalues of $J(\boldsymbol{U}_0)$ can be estimated by performing Arnoldi iterations with shift $\mu_k \coloneqq \frac{\widetilde{h}}{h} \widetilde{\mu}_k$ to compute estimates to the eigenvalues close to the convex hull point $\mu_k$.
		This can be done efficiently as the Arnoldi iterations for the different hull points $\mu_k, k = 1, \dots, K$ can be parallelized.
		Furthermore, as the Arnoldi iteration involves only matrix-vector products, this approach heavily benefits from storing $J(\boldsymbol{U}_0)$ in sparse format.
		Moreover, the stopping tolerance can be set relatively high as an estimate to the spectrum $\boldsymbol \sigma\left( J(\boldsymbol U_0) \right)$ is not only sufficient, but in fact beneficial since $\widetilde{\boldsymbol \sigma}\left( J(\boldsymbol U) \right)$ changes over the course of the simulation.
		This leads to an overall more robust stability polynomial.
		In terms of software we use \texttt{Arpack.jl}, a \texttt{Julia} wrapper around \texttt{arpack-ng} \cite{lehoucq1998arpack}.
		
		In this work, we try to estimate the outer shape of the spectrum with 1000 eigenvalues, hence supplying to each Arnoldi iteration a target number of $M \coloneqq \sfrac{1000}{K}$ eigenvalues.
		In general, after combining the estimates $\boldsymbol \lambda_k$ one will have less than 1000 distinct eigenvalues as some eigenvalues are found from multiple different shifts $\mu_k$.
		In practice, this still gives a sufficient large set of eigenvalues $\boldsymbol \lambda = \cup_{k=1}^K \boldsymbol \lambda_k$.
		\item If also in step 3 a reduced system had to be considered (for instance due to memory restrictions), scale the obtained eigenvalues $\boldsymbol \lambda$ accordingly.
	\end{enumerate}
	To illustrate the advantages of the presented routine, the reduced spectrum 
	and the complete spectrum computed through a standard eigendecomposition are presented in \cref{fig:SpectraComparisonReducedFull}.
	The spectrum shown in this case corresponds to the initial condition leading to a Kelvin-Helmholtz instability (see \cite{rueda2021subcell} for the precise setup) with the \ac{dgsem} \cite{black1999conservative, kopriva2009implementing} with flux-differencing \cite{gassner2013skew, gassner2016split} using solution polynomials of degree $k=3$, \ac{hlle} surface-flux \cite{einfeldt1988godunov}, subcell shock-capturing \cite{hennemann2021provably} to suppress spurious oscillations, and entropy-conserving volume flux \cite{fisher2013high, ranocha2020entropy}.
	Clearly, the shape of the spectrum is captured very well by the reduced spectrum while being significantly faster to compute.
	Furthermore, by estimating the spectrum using the procedure outlined above, one automatically obtains a reduced spectrum that captures the essential information while speeding up the subsequent optimization of the stability polynomial (discussed in the next section) as there are fewer redundant constraints to be checked.
	In this case, the shifts $\widetilde{\mu}_k$ are obtained from a reduced system with $\widetilde{J} \in \mathbb R^{1024\times 1024}$ for which the eigendecomposition can be performed in $0.4$ seconds.
	The large system is for the purpose of this study chosen such that the full eigendecomposition can still performed in somewhat reasonable time.
	Here, we consider $J \in \mathbb R^{16384\times 16384}$ where the full eigendecomposition takes more than $11$ minutes to complete.
	In contrast, for $K=20$ scaled shifts the Arnoldi iterations require in total only about $30$ seconds to complete.
	Here, $\vert \boldsymbol \lambda \vert = 930$ distinct eigenvalues were found which are displayed in \cref{fig:ReducedSpectrum}.

	We use this spectra estimation for the applications in \cref{sec:Applications} where the problem setup requires already for the initial condition a relatively fine grid resolution.
	Precisely, these are the visco-resistive Orszag-Tang vortex (\ref{subsec:ViscoResistiveOrszagTang}), the Kelvin-Helmholtz instability (\ref{subsec:KelvinHelmholtzInstability}), and Rayleigh-Taylor instability (\ref{subsec:RayleighTaylorInstability}).
	\subsection{Optimization of the Stability Polynomial}
	Equipped with a set of constraining eigenvalues, it remains to find stability polynomials of prescribed degrees $E^{(r)}$ and linear order $p$ 
	which maximize the stable timestep 
	$\Delta t$.
	For this, we use the procedure developed in \cite{ketcheson2013optimal} which allows (in standard double precision) for an efficient and reliable optimization of general stability polynomials up to, say, degree 16.
	In principle, also the optimization of higher-degree polynomials is possible if a suitable basis, i.e., different from the monomials $z^j$ can be found for the spectrum of interest.
	While this is possible for simple cases like parabolic (negative real axis), purely hyperbolic (imaginary axis) and circle/disk spectra, the choice of the basis for the general case is unclear \cite{ketcheson2013optimal}.
	The computational cost of the optimization procedure depends on the number of constraints, i.e., the number of supplied eigenvalues $\boldsymbol \sigma$.
	For moderate numbers of eigenvalues, say, up to 2000 (which are in practice sufficient to estimate the spectrum of the operator) the optimziation can be done in less than a minute.
	
	In the case of quad/octree-based \ac{AMR} the characteristic cell size is halved with every level of refinement.
	This implies in the case of convection-dominated spectra \eqref{eq:Convection_Diffusion_Scaling} that the stable timestep needs to be reduced by a factor two, or, in the case of stabilized Runge-Kutta methods, that the number of stage-evaluation needs to be doubled.
	
	Hence, for an $R$-times refined grid, one needs to construct stability polynomials with degrees up to $S^{(R+1)} = S^{(1)} \cdot 2^R$, i.e., the polynomial degree scales exponentially in the number of grid refinements $R$.
	As mentioned above, the optimization approach developed in \cite{ketcheson2013optimal} is for general spectra (for which no special suited basis can be found) not capable of reliably optimizing such high degree polynomials.
	Through employing a different ansatz for the optimization, one can optimize stability polynomials for convection dominated spectra with degrees larger than 100 \cite{doehring2024manystage}.
	
	Thus, we have in principle the possibility to integrate four to five times refined meshes with the same timestep as used for the base discretization.
	However, as discussed in \cref{subsec:NonlinearStabPERK} the nonlinear stability properties of \ac{PERK} schemes worsen quickly for an increasing difference in stage evaluations $\Delta E^{(r+1)} \coloneqq E^{(r+1)} - E^{(r)}$, cf. \cref{tab:TotalVariationIncrease_EVariation}.
	This renders the application of \ac{PERK} schemes with stage evaluations $E^{(r)} = E^{(1)} \cdot \{1, \allowbreak 2, \allowbreak 4, \allowbreak  8, \allowbreak  16, ...\}$ very inefficient, as the overall timestep needs to be reduced due to the poor nonlinear stability properties.
	As a consequence, we will restrict ourselves for practical cases presented in \cref{sec:Applications} to few-stage \ac{PERK} families.
	\subsection{Constructing P-ERK Methods from Stability Polynomials}
	Equipped with an optimized stability polynomial $P_{p;E^{(r)}}(z)$ it remains to construct the corresponding Runge-Kutta method.
	As the weights $\boldsymbol b^T$ and absciassae $\boldsymbol c$ are shared for all members of the family, these are fixed a priori.
	Thus, it remains to construct the free parameters of the Butcher arrays $A^{(r)}$, i.e., the entries $a_{i,i-1}^{(r)}$.
	For the sake of readability, we will truncate $a_{i, i-1}^{(r)}$ to $a_{i}$ in this section.
	The central relation between the stability polynomial and the Butcher array is given by \cite[Chapter~IV.7]{HairerWanner2}
	\begin{equation}
		P(z) = 1 + z \boldsymbol b^T \left( I - z A \right)^{-1} \boldsymbol 1
	\end{equation}
	where $\boldsymbol 1$ is the column-vector of ones.
	Following \cite{ketcheson2013optimal} the optimal stability polynomials are represented in monomial form, i.e., 
	\begin{equation}
			P_{p;E}(z) = \sum_{i=0}^{p} \frac{z^i}{i!} + \sum_{i=p+1}^{E} \alpha_i z^i, \quad \boldsymbol \alpha \in \mathbb R^{E-p}.
	\end{equation}
	with optimization variables $\boldsymbol \alpha$.
	\subsubsection{Second-Order P-ERK Methods}
	For the second-order \ac{PERK} methods, the general form of the stability polynomial from an $S$-stage method is given by
	\begin{equation}
		\begin{split}
			\label{eq:PERK2_StabPnom}
			P_{2;E}(z) = & \:  1 + z + \frac12 z^2 + c_{S-1} a_{S}z^3 + c_{S-2} a_{S}a_{S-1}z^4 \\
			+ & \dots + c_{S-E+2} z^{S-E+2}\prod_{i=S-E+3}^{S} a_{i}.
		\end{split}
	\end{equation}
	Equipped with the optimized monomial coefficients $\boldsymbol \alpha$ we can now construct the $a_{i} = a_{i, i-1}$ sequentially via
	\begin{equation}
		\label{eq:PERKCoeffsConstruction}
		a_{S-i} = \frac{\alpha_{3+i}}{c_{S-i-1}} \prod_{j=0}^{i-1} \frac{1}{a_{S-j}},
		\quad i = 0, \dots , S - 3.
	\end{equation}
	\subsubsection{Third-Order P-ERK Methods}
	For the third-order \ac{PERK} methods, the general form of the stability polynomial from an $S$-stage method is given by
	\begin{equation}
		\begin{split}
			\label{eq:PERK3_StabPnom}
			P_{3;E}(z) = & \:  1 + z + \frac12 z^2 + \frac16 z^3 \\ 
			+ & 
			\sum_{i = 1}^{S-4} \left( c_{S-1-i} a_{S} b_S \prod_{j=1}^{i} a_{S-j} 
			+ c_{S-2-i} b_{S-1} a_{S-1} \prod_{j=1}^i a_{S-1-j}\right)z^{i+3} \\
			+ & c_{2} b_S a_S \prod_{j=1}^{i} a_{S-j} z^S \: .
		\end{split}
	\end{equation}
	Besides matching the monomial coefficients $\boldsymbol \alpha$ we have a constraint arising from the third-order consistency conditions, cf. \eqref{eq:PRKOrderCondsSimplified}, see also \cite{nasab2022third}.
	\begin{equation}
		\frac{1}{6} \overset{!}{=} b_S a_{S} + b_{S-1} a_{S-1} \: .
	\end{equation}
	In contrast to the second-order case, the coefficients $a_{i}$ cannot be directly computed from the monomial coefficients, but instead a nonlinear system of equations needs to be solved.
	Additionally, we demand that there are no negative entries in the Butcher arrays $A^{(r)}$ corresponding to downwinding \cite{gottlieb1998total} by requiring that the $a_{i}$ and also $a_{i, 1} = c_i - a_i$ are positive.
	In practice, we find such a solution relatively quick using a standard trust-region method with randomly initialized starting points. 
	\section{Validation}
	\label{sec:Validation}
	Before coming to applications we seek to compare the combined \ac{PERK} schemes to standalone optimized schemes with regard to errors and conservation properties.
	For this, we consider the classic isentropic vortex testcase \cite{shu1988efficient, wang2013high} with the same parameters as used in \cite{vermeire2019paired, nasab2022third}, except that here we use a more compact domain $\Omega = [0, 10]^2$.
	In addition to the convergence studies presented in \cite{vermeire2019paired, nasab2022third} we compare the errors of the \ac{PERK} schemes to the errors of an individual scheme.
	This information is essential for assessing the effectiveness of the \ac{PERK} schemes since one is ultimately interested in accuracy per computational cost.
	For this testcase, the domain is initially discretized with 64 elements per direction, on which the solution is reconstructed using polynomials with degrees two (for second-order \ac{PERK} scheme) and three (for third-order \ac{PERK} scheme).
	The mesh is dynamically refined around the vortex based on the distance to the center of the vortex with two additional levels of refinement, i.e., the minimal grid size is in this case $h = \frac{10}{256}$.
	The numerical fluxes are computed using the \ac{hllc} flux \cite{toro1994restoration}.
	For solution polynomials of degree $k=2$ we construct $p=2$ optimal stability polynomials with degrees $3, 6,$ and $12$, while for solution polynomials of degree $k=3$, \allowbreak $p=3$ stability polynomials of degrees $4, 8$, and $16$ are optimized.
	We choose this particular \ac{PERK} families since the presence of the additional free parameter allows the optimization of the lowest-degree stability polynomials, i.e., $P_{2;3}(z), P_{3;4}(z)$, to the spectrum.
	
	The density errors at final time $t_f = 20$
	\begin{equation}
		e_{\rho}(x, y) \coloneqq \rho(t_f, x,y) - \rho_h(t_f, x,y)
	\end{equation}
	are tabulated in \cref{tab:DensityErrors_IsentropicVortex} for each individual optimized Runge-Kutta method and the composed \ac{PERK} scheme.
	\begin{table}
		\def\arraystretch{1.2}
		\centering
		\begin{tabular}{l?{2}c|c}
			Method & $\frac{1}{\vert \Omega \vert }\Vert e_{\rho}(x, y) \Vert_{L^1(\Omega)}$ & $\Vert e_{\rho}(x, y) \Vert_{L^\infty(\Omega)}$ \\
			\Xhline{5\arrayrulewidth}
			$\text{P-ERK}_{2;3}$ & $1.77 \cdot 10^{-5}$ & $6.43\cdot 10^{-4}$ \\
			$\text{P-ERK}_{2;6}$ & $1.77 \cdot 10^{-5}$ & $6.40\cdot 10^{-4}$ \\			
			$\text{P-ERK}_{2;12}$ & $1.74 \cdot 10^{-5}$ & $6.30\cdot 10^{-4}$ \\	
			$\text{P-ERK}_{2;\{3, 6, 12\}}$ & $2.09 \cdot 10^{-5}$ & $5.86 \cdot 10^{-4}$ \\
			\Xhline{2\arrayrulewidth}
			$\text{P-ERK}_{3; 4}$ & $4.05 \cdot 10^{-7}$ & $1.12\cdot 10^{-5}$ \\
			$\text{P-ERK}_{3;  8}$ & $4.05 \cdot 10^{-7}$ & $1.12\cdot 10^{-5}$ \\			
			$\text{P-ERK}_{3; 16}$ & $4.05 \cdot 10^{-7}$ & $1.11\cdot 10^{-5}$ \\	
			$\text{P-ERK}_{3; \{4, 8, 16\}}$ & $3.73 \cdot 10^{-7}$ & $1.36 \cdot 10^{-5}$
		\end{tabular}
		\caption[Density-errors for isentropic vortex advection testcase.]
		{Density-errors for isentropic vortex advection testcase after one pass through the domain $\Omega = [0, 10]^2$ at $t_f = 20$.
		The compressible Euler equations are discretized using the \ac{dgsem} with \ac{hllc} flux \cite{toro1994restoration} on a dynamically refined mesh with cell sizes ranging from $\Delta x^{(1)} = \sfrac{10}{64}$ to $\Delta x^{(3)} = \sfrac{2}{256}$.
		For Runge-Kutta order $p=2$, the local polynomials are of degree $k=2$ while for Runge-Kutta order $p=3$ the solution is constructed using $k=3$ polynomials.}
		\label{tab:DensityErrors_IsentropicVortex}
	\end{table}
	Both $L^\infty$ and domain-normalized $L^1$ errors are for the joint \ac{PERK} schemes of similar size as the errors of the individual schemes.

\begin{figure}[!t]
		\centering
		\subfloat[{Density errors for second-order $E=\{3,6,12\}$ \ac{PERK} family.}]{
			\label{fig:DensityErrorsPERK_P2_CovergenceTest}
			\centering
			\includegraphics[width=.45\textwidth]{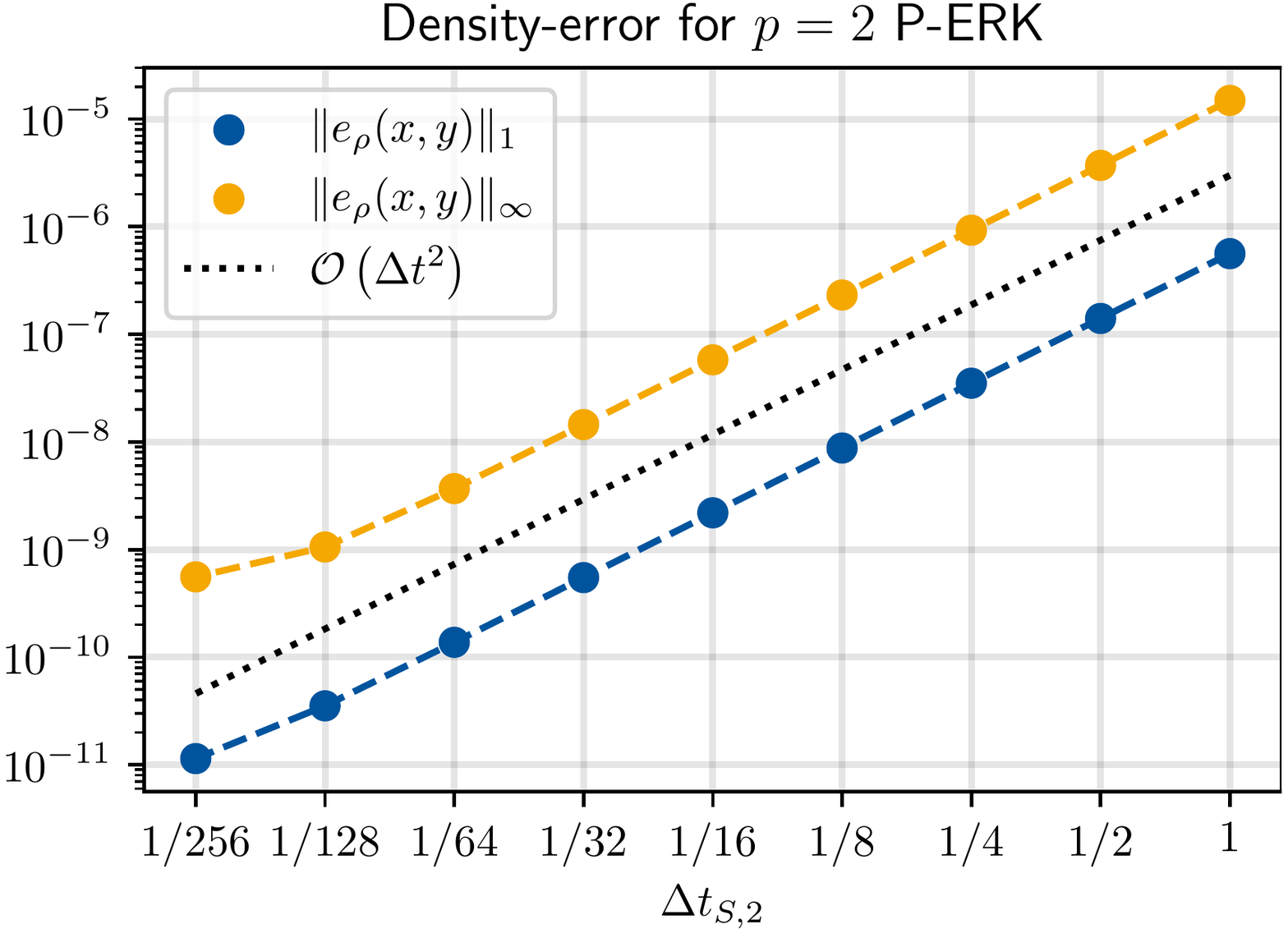}
		}
		\hfill
		\subfloat[{Density errors for third-order $E=\{4,8,16\}$ \ac{PERK} family.}]{
			\label{fig:DensityErrorsPERK_P3_CovergenceTest}
			\centering
			\includegraphics[width=.45\textwidth]{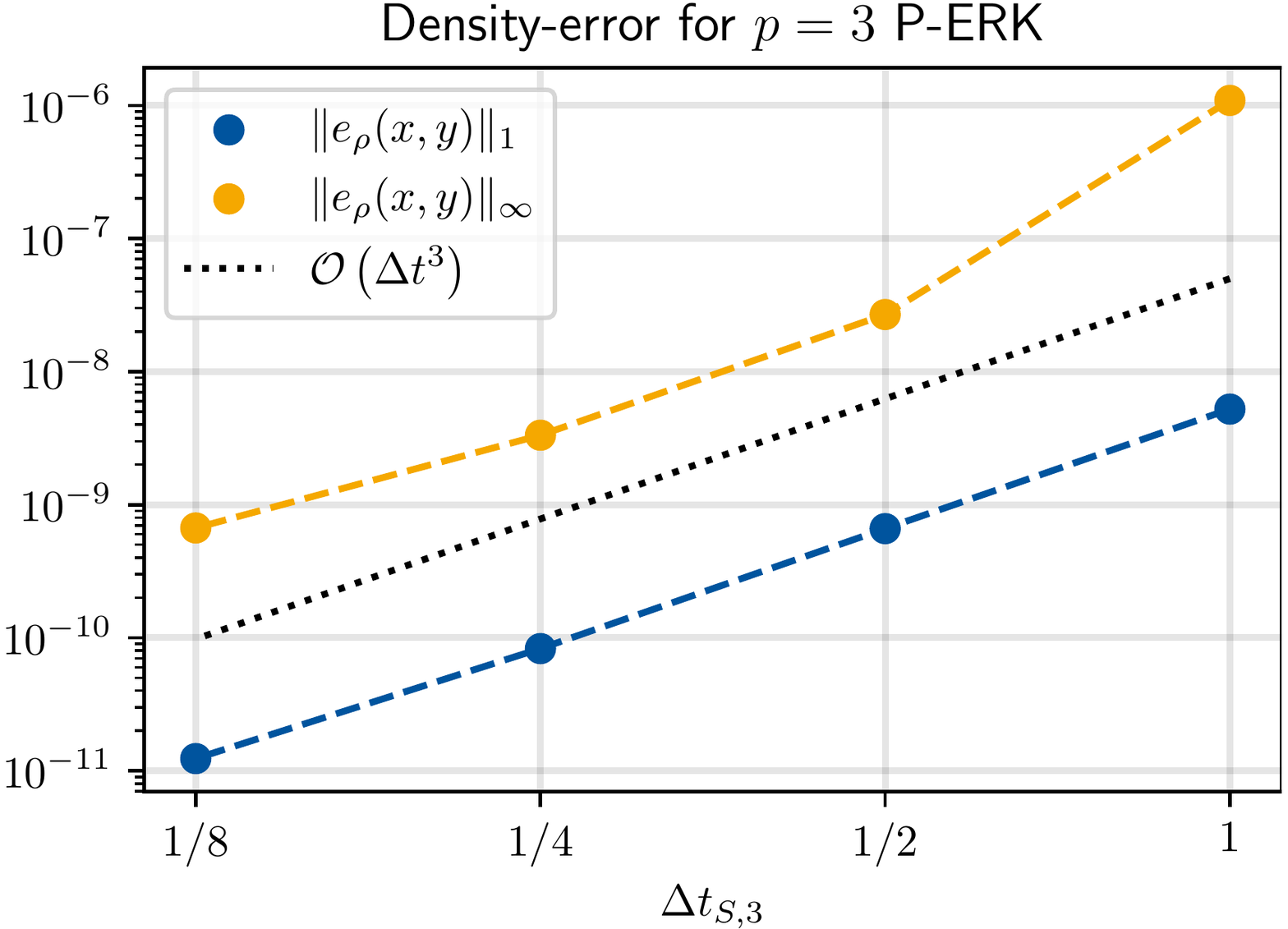}
		}
		\caption[$L^1$ and $L^\infty$ density errors for second and third-order PERK.]
		{$L^1$ and $L^\infty$ density errors for second (\cref{fig:DensityErrorsPERK_P2_CovergenceTest}) and third-order (\cref{fig:DensityErrorsPERK_P3_CovergenceTest}) \ac{PERK} schemes applied to the isentropic vortex advection testcase.
		We observe the expected order of convergence for both second and third-order \ac{PERK} schemes until the spatial discretization error dominates the overall error.
		This becomes especially visible for the $L^\infty$ errors plotted in orange.}
		\label{fig:DensityErrorsPERK_P2P3_CovergenceTest}
	\end{figure}

	In \cite{hundsdorfer2015error} it is shown that at the interfaces of different schemes order reduction \cite[Chapter~II.2]{hundsdorfer2003numerical} may be encountered.
	As argued in \cite{hundsdorfer2015error}, the interfaces between methods act like additional time-dependent internal boundary conditions, thus potentially suffering from effects also encountered for non-periodic problems \cite[Chapter~II.2]{hundsdorfer2003numerical}.
	Nevertheless, for both $L^1$ and $L^\infty$ error, both second and third-order \ac{PERK} schemes display the expected order of convergence, as illustrated in \cref{fig:DensityErrorsPERK_P2P3_CovergenceTest}.
	As in \cite{vermeire2019paired, nasab2022third}, we employ local polynomials of degree $k=6$ for the convergence tests for both second and third-order \ac{PERK} to minimize spatial errors.

	It remains to examine whether the application of \ac{PERK} schemes causes additional violation of the conserved fields compared to a standard Runge-Kutta scheme.
	For this purpose the third-order \ac{PERK} scheme is compared against the three stage, third-order \ac{ssp} method by Shu and Osher.
	As for the comparison of the density-errors of the \ac{PERK} schemes to standalone schemes, the local polynomials are of degree $k=3$.
	The conservation error 
	\begin{equation}
		e_u^\text{Cons}(t_f) \coloneqq \left \vert \int_\Omega u(t_0, \boldsymbol x) \nid \boldsymbol x - \int_\Omega u(t_f, \boldsymbol x) \nid \boldsymbol x \right \vert
	\end{equation}
	is computed for mass, momentum and energy after one pass of the vortex through the domain and tabulated in \cref{tab:ConservationErrors_IsentropicVortex}.
	\begin{table}
		\def\arraystretch{1.2}
		\centering
		\begin{tabular}{l?{2}c|c}
			 & $\text{P-ERK}_{3;\{4, 8, 16\}}$ & $\text{SSP}_{3;3}$ \\
			 	\Xhline{5\arrayrulewidth}
			 $e_\rho^\text{Cons}(t_f)$ & $3.94 \cdot 10^{-13}$ & $7.82 \cdot 10^{-13} $ \\
			 $e_{\rho v_x}^\text{Cons}(t_f)$ & $4.79 \cdot 10^{-14}$ & $4.89 \cdot 10^{-13} $ \\
			 $e_{\rho v_y}^\text{Cons}(t_f)$ & $7.11 \cdot 10^{-14}$ & $5.00 \cdot 10^{-13} $ \\
			 $e_{\rho e}^\text{Cons}(t_f)$ & $1.06 \cdot 10^{-12}$ & $3.64 \cdot 10^{-12} $ \\
		\end{tabular}
		\caption[Conservation-errors for isentropic vortex advection testcase.]
		{Conservation-errors for isentropic vortex advection testcase after one pass through the domain at $t_f = 20$ for the third-order \ac{PERK} scheme and the three-stage, third-order Shu-Osher scheme.}
		\label{tab:ConservationErrors_IsentropicVortex}
	\end{table}
	Clearly, the errors of the $p=3, E=\{4, 8, 16\}$ \ac{PERK} family are on par with the errors for the Shu-Osher scheme.
	In particular, these conservation violations are dominated by the spatial accuracy, i.e., reducing the timestep has very little influence.
	We also compared the $p=2, E=\{3, 6, 12\}$ \ac{PERK} family to the two-stage second-order \ac{ssp} scheme (Heun's method) for a discretization with $k=2$ solution polynomials.
	For this configuration the conservation errors are negligible, both for the \ac{PERK} and $\text{SSP}_{2;2}$ scheme.

	\vspace{0.5\baselineskip} 
	All presented results from this section can be reproduced using the publicly available repository \cite{doehring2024multirateRepro}.
	\section{Methodology}
	\label{sec:Methodology}
	In this section, we outline the procedure for the application of \ac{PERK} schemes to simulations with quad/octree-based \ac{AMR} in \texttt{Trixi.jl}.
	\subsection{Implementation}
	For obtaining the spatial semidiscretization \eqref{eq:Semidiscretization} we use the \ac{dgsem} \cite{black1999conservative, kopriva2009implementing} with flux-differencing \cite{gassner2013skew, gassner2016split} as implemented in \texttt{Trixi.jl}.
	To demonstrate that indeed no changes to the spatial discretization are required, we exemplary provide the standard implementation of the volume integral computation in \texttt{Trixi.jl} in \cref{listing:Trixi} and the implementation used for the \ac{PERK} schemes in \cref{listing:PERK}.
	The only differences is that the \ac{PERK} version takes an additional argument \texttt{elements\_r} which contains the indices of the elements of partition $r$.
	\begin{lstlisting}[language=Python, label=listing:Trixi, caption={Volume integral computation for standard time integration schemes in \texttt{Trixi.jl} \cite{trixi1,trixi2,trixi3}.}, captionpos=b]
function calc_volume_integral!(du, u, mesh, nonconservative_terms, equations,
                               volume_integral::VolumeIntegralWeakForm,
															 dg::DGSEM, cache)
	# Loop over all elements in the cache
	@threaded for element in eachelement(dg, cache)
		weak_form_kernel!(du, u, element, mesh, nonconservative_terms,
	                    equations, dg, cache)
	end
	return nothing
end\end{lstlisting}
	\begin{lstlisting}[language=Python, label=listing:PERK, caption={Volume integral computation for \ac{PERK} schemes in \texttt{Trixi.jl} \cite{doehring2024multirateRepro}.
			In contrast to the standard volume integral computation \cref{listing:Trixi}, the \ac{PERK} version takes an additional argument \texttt{elements\_r} which contains the elements bundling the $k+1$ polynomial coefficients from $\boldsymbol U^{(r)}$, cf. \cref{eq:PartitionedODESys2}.}, captionpos=b]
function calc_volume_integral!(du, u, mesh, nonconservative_terms, equations,
															 volume_integral::VolumeIntegralWeakForm,
															 dg::DGSEM, cache,
															 elements_r::Vector{Int64})
	# Loop over elements of the r'th level
	@threaded for element in elements_r
		weak_form_kernel!(du, u, element, mesh, nonconservative_terms,
											equations, dg, cache)
	end
	return nothing
end\end{lstlisting}
	\subsection{Dynamic Partitioning based on AMR}
	In this work we employ two types of quad-octree based meshes with \ac{AMR}.
	First, we consider a \texttt{TreeMesh} based on a hierarchical tree data structure with uniform cell sizes obtained from subsequent refinements of a square/cubic domain.
	For these types of meshes, the assignment of cells to partitions is trivial as the cellsizes correspond to the position of the nodes in the tree, which makes it very efficient to query the characteristic grid size of each cell.
	Additionally, we also employ in principle unstructured meshes based on the \texttt{p4est} library \cite{BursteddeWilcoxGhattas11}.
	In that case, we cannot simply use the position of a cell in the tree data structure as a measure for the characteristic grid size, since cellsizes can vary even though they are formally on the same level of refinement.
	This implies that the computation of a criterion (e.g., shortest edge length, inner/outer cell diameter, determinant of the coordinate-transformation matrix, ...) is required for the assignment of each cell which introduces additional computational work.
	Most examples considered in this work employ the first mesh type but we also present one application of the latter, where the shortest edge length of each cell is recomputed in the partition assignment step.
	
	\texttt{Trixi.jl} requires four data structures to be assigned to the partitions.
	First, \texttt{elements} and its associated \texttt{boundaries} are are trivially assigned according to the size of the associated cell.
	Second, \texttt{interfaces} are by construction only between conforming cells of thus identical size and are hence also uniquely associated with a level.
	Finally, \texttt{mortars} (which take care of hanging nodes at the non-conforming interface of cells of different sizes \cite{kopriva1996conservative}) are associated to the partition corresponding to the fine neighbor cells in order to receive updates whenever the fine neighboring cell is updated.
	\subsection{Aptitude for Massively-Parallel, Distributed-Memory Architectures}
	We would like to emphasize that the presented \ac{PERK} schemes are particularly well-suited for massively-parallel, distributed-memory architectures.
	This is due to the fact that the partitioning of the cells belonging to a shared-memory unit can be done completely independently of the other units, assuming a unique partitioning criterion is available.
	In the case of distributed storage of the mesh data, a functionality provided by \texttt{p4est}, the index arrays of the \texttt{elements}, \texttt{boundaries}, \texttt{interfaces}, and \texttt{mortars} would then also be stored on each memory unit only.
	In particular, no additional data exchange is needed if the \ac{PERK} method is used over a standard Runge-Kutta method.
	\subsection{Assessing Performance}
	It is clear that the potential speedup of simulation using \ac{PERK} is strongly depending on the ratio of coarse to fine cells.
	In case of two levels of refinement, we can expect almost a speedup of factor two if there is only one refined cell and the remainder of the mesh is coarse.
	Conversely, there is effectively no speedup to be gained if there is only one coarse cell left with the rest of the mesh refined.
	\subsubsection{Actual Number of Right-Hand-Side Evaluations}
	In order to accurately assess the performance of a given \ac{PERK} scheme when applied to a specific problem we require a way to compute the optimal, minimal runtime.
	In the ideal scenario, one would have an optimized scheme for every cell size present in the mesh such that the timestep does not need to be reduced if the cells are further refined throughout the simulation.
	As argued in \cref{subsec:NonlinearStabPERK} this is in practice not possible due to the loss of monotonicity that occurs for the usage of schemes with increasing difference in stage evaluations $\Delta E^{(r+1)} \coloneqq E^{(r+1)} - E^{(r)}$.
	Consequently, we restrict ourselves in practice to, e.g., a three-level ($R=3$) \ac{PERK} family with moderate number of stage evaluations such as $E= \{3, 6, 12\}$.
	If, however, there is now a mesh with more than three levels $L$ of refinement, i.e., $L>R$ optimality is inevitably lost since either the lowest or highest level is integrated non-optimal.
	In this work, we always use the highest stage-evaluation method for the finest level and fill the lower levels with the lowest stage-evaluation method.
	For instance, if we have a mesh with six levels of refinement and a \ac{PERK} scheme with $E= \{3, 6, 12\}$, we would use for the finest cells the $E=12$ scheme, for the second finest cells the $E=6$ scheme, and for the remaining levels the $E=3$ scheme.
	For this case, the number of scalar, i.e., per \ac{DG} solution coefficient, \ac{RHS} evaluations $N_\text{RHS}^\text{Actual}$ can be computed as
	\begin{equation}
		\label{eq:ActualRHSEvals}
		N_\text{RHS}^\text{Actual} = \sum_{i=1}^{ \lceil N_t/ N_\text{AMR} \rceil } N_\text{AMR} \, N \, (k+1)^{N_D} \left  [\sum_{l=1}^R E^{(l)} N_\text{C}(i, l) + \sum_{l=R+1}^L E^{(1)} N_\text{C}(i, l) \right ]
	\end{equation}
	In the formula above, $N_t$ denotes the overall number of timesteps to progress the simulation from initial time $t_0$ to final time $t_f$.
	$N_\text{AMR}$ is the number of timesteps taken until the next (potential) change of the mesh, $N$ denotes the number of fields of the \ac{PDE} \eqref{eq:ConservationLaw}, $k$ is the \ac{DG} solution polynomial degree, and $N_D \in \{1, 2, 3\}$ denotes the number of spatial dimensions.
	$N_\text{C}(i, l)$ denotes the number of cells associated with the $l$'th level assigned at the $i$'th \ac{AMR} call and $E^{(l)}$ denotes the number of stage evaluations for the $l$'th level.
	\subsubsection{Minimal Number of Right-Hand-Side Evaluations}
	For an optimal scheme one would have a method $E^{(l)}$ for every level $l=1, \dots, L$ of refinement.
	In this case, the second sum in \cref{eq:ActualRHSEvals} is not necessary as the first sum would already cover all levels of refinement, i.e., $R=L$.
	Recalling that for the quad/octree meshes considered here the cell sizes vary from level to level by a factor of 2 it is clear that for a coarser cell, in principle, a timestep twice as large can be used.
	For the optimized Runge-Kutta schemes considered here, this translates to employing a method with half the number of stage evaluations $E$.
	Consequently, the total number of additional scalar, i.e., per \ac{DoF}, right-hand-side evaluations $N_\text{RHS+} \coloneqq N_\text{RHS}^\text{Actual} - N_\text{RHS}^{R=L} $ due to the lack of a full range of optimized schemes can be computed as
	\begin{equation}
		\label{eq:AdditionalRHSEvals}
		N_\text{RHS+} 
		= \sum_{i=1}^{\lceil N_t/ N_\text{AMR} \rceil } N_\text{AMR} \, N \, (k+1)^{N_D} \sum_{l=R+1}^L E^{(1)} \left(1 - \frac{1}{2^{l - R}}\right) N_\text{C}(i, l) \, .
	\end{equation}
	The term in parentheses computes the difference in actual stage evaluations $E^{(1)}$ to the number of in principal required stage evaluations $\sfrac{E^{(1)}}{2^{l - R}}$ for the $l=R+1, \dots, L$ coarser levels.

	For most examples considered in this work, a three-member \ac{PERK} family involving the non-optimized $p=3, E^{(1)}=3$ schemes is performing best since this leads to a \ac{PERK} scheme with minimal difference in stage evaluations $\Delta E$ which is advantageous when it comes to avoiding artificial oscillations at the interfaces, see \cref{tab:TotalVariationIncrease_EVariation}.
	Then, the linearly stable timestep is usually doubling already for the $E^{(2)} = 4$ method and doubling a second time for, say, the $E^{(3)} = 7$ scheme.
	While the timestep $\Delta t$ decreases by a factor of four from the $E^{(3)} = 7$ to the $E^{(1)} = 3$ scheme, the number of stage evaluations decreases only by a factor of $\sfrac{3}{7} > \frac{1}{4}$.
	In other words, such a \ac{PERK} scheme is non-optimal on every but the finest level $r=R$.
	Hence, to quantify the additional computational effort due to the lack of a full range of optimized schemes we need to compute the number of additional \ac{RHS} evaluations not only for those levels $l = R+1, \dots, L$ that are not integrated with the optimal scheme but also for the levels $l = 1, \dots, R-1$ for which methods are present.
	Thus, the optimal, gold-standard \ac{PERK} method has stage evaluations $E_\star = E^{(1)}_\star \cdot \left \{1, \allowbreak 2, \allowbreak 4, \allowbreak  \dots, \allowbreak  2^{L-1} \right \}$, or, equivalently, 
	$E_\star= E^{(R)} \cdot \left \{ 2^{-(L-1)}, \allowbreak 2^{-(L-2)}, \allowbreak \dots, \allowbreak 1 \right \} $.
	Consequently, the overall number of additional \ac{RHS} evaluations $N_\text{RHS+}$ compared to a hypothetical optimal \ac{PERK} scheme for which both timestep and number of stage evaluations halve for every member method the number of additional \ac{RHS} evaluations 
	\begin{equation}
		\label{eq:AdditionalRHSEvalsIdealDef}
		N_\text{RHS+}^\text{Hyp} \coloneqq N_\text{RHS}^\text{Actual} - N_\text{RHS}^\text{Opt} 
	\end{equation}
	which can be computed as
	\begin{equation}
		\begin{aligned} 
			\label{eq:AdditionalRHSEvalsIdeal}
			N_\text{RHS+}^\text{Hyp} = \sum_{i=1}^{\lceil N_t/ N_\text{AMR} \rceil } N_\text{AMR} \, N \, (k+1)^{N_D}
				\Bigg[ & \sum_{l=1}^{R-1}  \left( E^{(l)} -\frac{E^{(R)}}{2^{l-1}}\right) N_\text{C}(i, l) \\ 
			 + & \sum_{l=R+1}^L \left( E^{(1)} -\frac{E^{(R)}}{2^{l-1}}\right) N_\text{C}(i, l) \Bigg] \: .
			\end{aligned}
	\end{equation}
	Note that this assumes the same maximal linearly stable timestep is achieved for both the optimal and actual \ac{PERK} method, i.e., $E^{(R)} = E^{(R)}_\star$.

	It is worthwhile to discuss \cref{eq:AdditionalRHSEvalsIdeal} for a concrete example.
	In this work, the method with highest stage-evaluations is often $E^{(3)} = 6$ and it is thus not immediately clear how one would interpret e.g. $\frac{6}{2^2} = 1.5$ stage evaluations for the hypothetical optimal scheme.
	To resolve this, recall that for optimized stabilized Runge-Kutta methods the maximum linearly stable timestep increases linearly with the number of stage evaluations.
	Thus, the computational effort in the ideal case (i.e., without loss of monotonicity) of a \ac{PERK} scheme with base method of $E^{(1)} = 3$ is equivalent to a \ac{PERK} scheme with base method $E^{(1)} = 6$ since (again, in the ideal case) the timestep of the latter family can be doubled, hence cutting the number of timesteps $N_t$ and the \ac{AMR} interval $N_\text{AMR}$ in half.
	Thus, up to rounding, the number of scalar \ac{RHS} evaluations $N_\text{RHS}^\text{Opt}$ for two \ac{PERK} families with different number of base evaluations $E^{(1)}$ is the same:
	\begin{align}
		N_\text{RHS}^\text{Opt} &= \sum_{i=1}^{\lceil N_t/ N_\text{AMR} \rceil } N_\text{AMR} \, N \, (k+1)^{N_D} \sum_{l=1}^L E^{(l)} N_\text{C}(i, l) \\
		&= \sum_{i=1}^{ \lceil (N_t/2)/ (N_\text{AMR}/2) \rceil } \frac{N_\text{AMR}}{2} \, N \, (k+1)^{N_D} \sum_{l=1}^L 2 E^{(l)} N_\text{C}(i, l)
	\end{align}
	Now, one can in principal choose $E^{(R)}$ high enough such that $E^{(R)} / 2^{L - 1}$ is still an integer.

	\subsubsection{Degree of Optimality}
	Equipped with the number of additional \ac{RHS} evaluations $N_\text{RHS+}^\text{Hyp}$ and the average cost of evaluating the \ac{RHS} per cell $\tau_\text{RHS/C}$, which is provided by \texttt{Trixi.jl}, one can compute the potential savings of a fully-optimal scheme.
	Neglecting possible additional overhead due to the usage of more levels, the hypothetical optimal runtime $\tau_\text{Opt}$ can then be computed as
	\begin{subequations}
		\begin{align}
		\label{eq:DefinitionRuntime}
		\tau_\text{Opt} &\coloneqq  N_\text{RHS}^\text{Opt} \cdot \tau_\text{RHS/C} \\
		\label{eq:HypotheticalSpeedUp}
		&= \tau - N_\text{RHS+}^\text{Hyp} \cdot \tau_\text{RHS/C} \, ,
		\end{align}
	\end{subequations}
	where $\tau$ denotes the measured wallclock runtime of an actual simulation with a \ac{PERK} scheme and $\tau_\text{RHS/C}$ quantifies the average cost of evaluating the \ac{RHS} $\boldsymbol F$ \eqref{eq:SemidiscretizationODE} per cell.
	
	We remark that \eqref{eq:DefinitionRuntime} is a simplified model which assumes that the measured wallclock runtime $\tau$ is solely governed by the number of \ac{RHS} evaluations.
	In practice, we see that \texttt{Trixi.jl} may be more a memory- than compute-bound code and thus the actual runtime $\tau$ may not be directly proportional to the number of \ac{RHS} evaluations, but may be affected by other factors such as cache misses due to non-coalesced memory accesses.

	It is important to note that \cref{eq:HypotheticalSpeedUp} is only valid if the \ac{CFL} number of the actual simulation using a \ac{PERK} scheme does not need to be reduced below the value obtained during optimization for absolute stability.
	While this is the case for the hyperbolic-parabolic examples considered in this work, the practical \ac{CFL} number used in the simulation of hyperbolic \acp{PDE} needs to be reduced.
	In this case, the number of additional \ac{RHS} evaluations $N_\text{RHS+}^\text{Hyp}$ is given by 
	\begin{equation}
		\label{eq:AdditionalRHSEvalsIdealDef_CFLReduction}
		N_\text{RHS+}^\text{Hyp} \coloneqq N_\text{RHS}^\text{Actual} - \frac{\text{CFL}_\text{Actual}}{\text{CFL}_\text{Opt}} N_\text{RHS}^\text{Opt} \, .
	\end{equation}
	To examine the performance of the non-optimal schemes we will provide the ratio 
	\begin{equation}
		\label{eq:OptimalityMeasure}
		\kappa \coloneqq \frac{\tau_\text{Opt}}{ \tau} \leq 1
	\end{equation}
	for the applications discussed in \cref{sec:Applications}.
	This is a measure of optimality in the sense that this quantifies how much of the theoretically possible speedup is actually achieved in practice.
	\subsubsection{Memory Footprint}
	The additional data structures due to the partitioning increase the memory footprint of the method.
	This is a classic example of storage vs. computations where we save \ac{RHS} evaluations at the cost of storing additional information.
	In particular, for the herein employed \texttt{TreeMesh} and \texttt{p4est} mesh for each \texttt{element}, \texttt{interface}, \texttt{boundary}, and \texttt{mortar} an integer specifying the corresponding method $r$ is stored.
	For other mesh types, such as the structured mesh used in \ref{sec:AppSmoothHighResMesh} there are only cells required and thus only for those additional data structures indicating the partitioning are required.
	For the considered examples the partitioning variables cause on average a $10-30\%$ increase in memory footprint.
	\section{Applications}
	\label{sec:Applications}
	In this section, we present a number of applications of the \ac{PERK} schemes to classic test problems with quad/octree-based \ac{AMR}.
	All presented cases are made publicly available in the reproducibility repository \cite{doehring2024multirateRepro}.
	We begin by studying hyperbolic-parabolic \acp{PDE}, i.e., systems with viscous fluxes for which we expect in general a better performance of the \ac{PERK} schemes due to inherent damping of potential high-frequency oscillations at the interfaces of methods.
	Additionally, we also consider hyperbolic \acp{PDE} to assess whether the \ac{PERK} schemes can be used efficiently in practice for such systems without the need for additional stabilization.
	In particular, this is a novel contribution since the \ac{PERK} schemes have so far only been applied to hyperbolic-parabolic systems \cite{vermeire2019paired,nasab2022third}.

	For all considered examples, three-level ($R=3$) \ac{PERK} families are observed to be the most efficient \ac{PERK} schemes.
	In particular, the methods are always chosen such that the admissible timestep doubles between levels, yielding values of $\text{CFL} = \{1, 2, 4\}$, cf. \cref{eq:TimestepConstraint}.
	Schemes with more family members, i.e., $R>3$ are observed to suffer from the issues discussed in \cref{subsec:NonlinearStabPERK} which, in turn, demand a reduction in the \ac{CFL} number, rendering them ultimately less efficient than the three-level families.
	In this work, we compare the multi-level \ac{PERK} schemes always to a single optimized \ac{PERK} method and four third-order integrators optimized for \acp{PDE} implemented in \texttt{OrdinaryDiffEq.jl} \cite{DifferentialEquations.jl-2017}.
	These are a seven-stage low-dissipation, low-dispersion scheme optimized for wave propagation problems with \ac{DG} methods \cite{toulorge2012optimal},
	a five-stage method optimized for high-order spectral difference methods applied to wave-propagation problems \cite{parsani2013optimized}, a four stage method optimized for compressible Navier-Stokes equations \cite{kennedy2000low}, and a recently published five-stage method with timestep computation based on error-control \cite{ranocha2022optimized}.
	Furthermore, we compare also to canonical third-order Shu-Osher \ac{ssp} scheme \cite{shu1988efficient}, which is the considered benchmark method in \cite{vermeire2019paired, nasab2022third}.

	For all simulations we aim for cells as coarse as possible that still allow for a meaningful simulation.
	This is important since we could generate artificial speedup by e.g. allowing only four instead of six levels of refinement which would lead to many more lowest-level cells on which the speedup is gained.
	Additionally, we employ robust \ac{AMR} indicator setups to avoid a too strong influence of the inevitably different solutions of the different methods on the \ac{AMR} process and thus on the mesh and computational cost.
	\subsection{Hyperbolic-Parabolic Equations}
	The presence of viscous fluxes makes the \ac{PERK} schemes less vulnerable to oscillations due to monotonicity violations and are thus, at least in the convection-dominated regime, a suitable choice for the integration of such systems.
	In fact, for the hyperbolic-parabolic equations considered here, the $R=3$ \ac{PERK} schemes can be run with the same \ac{CFL} number as the optimized single \ac{PERK} method.
	The hyperbolic-parabolic problems are discretized using the \ac{dgsem}  \cite{black1999conservative, kopriva2009implementing} where the viscous fluxes are treated using the Bassi-Rebay 1 (BR1) scheme \cite{bassi1997high, gassner2018br1}.
	The admissible timestep is in every case computed based on the inviscid \ac{CFL} condition \eqref{eq:CFL_DG}, \eqref{eq:TimestepConstraint} where we ensured that the invisicid terms are indeed limiting the timestep, i.e., not the viscous terms, cf. \eqref{eq:Convection_Diffusion_Scaling}.
	\subsubsection{Periodic Shear-Layer}
	\label{sec:PeriodicShearLayer}
	This example is selected to showcase the potential speed up from using the \ac{PERK} schemes.
	In contrast to the majority of the other examples shown in this work, the finest cells make up only for a small fraction of the entire grid.
	In particular, the mesh is refined only in confined regions where the roll-up of emerging vortices is observed, keeping the majority of the domain relatively coarse.
	
\begin{figure}[!t]
		\centering
		\subfloat[{Velocity $v_x$ in $x$-direction at $t = 0.8$ for the periodic shear layer.}]{
			\label{fig:ShearLayer_vx}
			\centering
			\includegraphics[height=.45\textwidth]{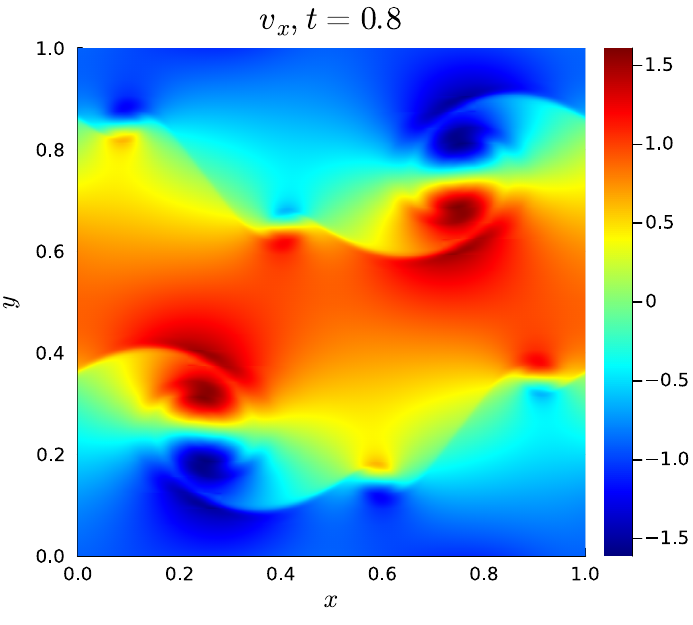}
		}
		\hfill
		\subfloat[{Five times refined mesh at $t = 0.8$.}]{
			\label{fig:ShearLayer_Mesh}
			\centering
			\includegraphics[height=.45\textwidth]{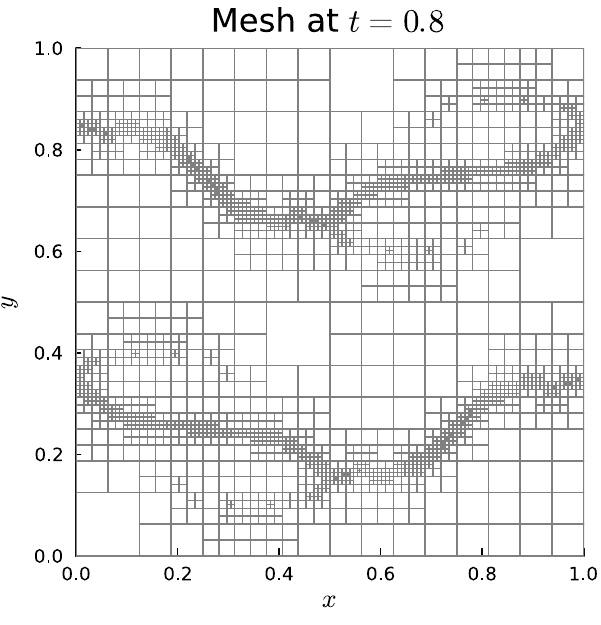}
		}
		\caption[Velocity $v_x$ in $x$-direction at $t = 0.8$ for the doubly periodic double shear layer testcase and corresponding adaptively refined mesh.]
		{Velocity $v_x$ in $x$-direction at $t = 0.8$ for the doubly periodic double shear layer testcase \cref{eq:ShearLayerIC} (\cref{fig:ShearLayer_vx}) and corresponding adaptively refined mesh \cref{fig:ShearLayer_Mesh}.}
		\label{fig:ShearLayer_vx_Mesh}
	\end{figure}

	Following the configuration described in \cite{brown1995performance}, we set
	\begin{equation}
		\label{eq:ShearLayerIC}
		\boldsymbol u_\text{primal}(t_0=0, x, y) = \begin{pmatrix} \rho \\ v_x \\ v_y \\ p \end{pmatrix} = \begin{pmatrix} 1 \\ \begin{cases} \tanh(\kappa (y - 0.25)) & y \leq 0.5 \\ \tanh(\kappa(0.75 - y)) & y > 0.5 \end{cases} \\ \delta \sin\big( 2 \pi (x + 0.25 )\big) \\ c^2 \rho / \gamma \end{pmatrix}
	\end{equation}
	with flow parameters $\kappa = 80, \delta = 0.05, c_0 = 10, \text{Ma}_0 = 0.1, \text{Re}_0 = 2.5 \cdot 10^{4}, \text{Pr} = 0.72$ and run the simulation until $t_f = 1.2$ as in \cite{brown1995performance}.
	In terms of spatial discretization, local polynomials of degree $k=3$ with \ac{hllc} surface-flux and flux differencing for the entropy-stable volume flux by Ranocha \cite{ranocha2018comparison} are employed to discretize the spatial derivatives.
	In contrast to the subsequent examples the subcell shock-capturing \cite{hennemann2021provably} is not required to suppress any spurious oscillations.
	The domain $\Omega = [0, 1]^2$ is initially heavily underresolved with only $2^3=8$ cells per coordinate direction.
	To allow for a physically meaningful simulation that resolves the vortex build-up, \ac{AMR} is employed based on the indicator by Löhner \cite{lohner1987adaptive} with indicator variable $v_x$.
	The implementation in \texttt{Trixi.jl} of the indicator by Löhner follows the adaptations to Legendre-Gauss-Lobatto elements developed in the FLASH code \cite{fryxell2000flash}.
	In this simulation, cells are allowed to be refined up to six times where the finest level corresponds to a discretization with $2^9 = 512$ cells per dimension.
	The velocity component $v_x$ and the computational mesh at $t=0.8$ are displayed in \cref{fig:ShearLayer_vx_Mesh}.

	For this setup, a $p=3$ \ac{PERK} scheme with $E = \{3, 4, 7\}$ stage evaluations is constructed, where the admissible timestep doubles for each family member.
	The mesh is allowed to change every $N_\text{AMR} = 20$ timesteps, which involves the recomputation of the flagging-variables for the \ac{PERK} integrator.
	In this case, the latter consumes $2.5\%$ of the total runtime with is the overhead due to the partitioning.
	Comparing this to the speedup of the stand-alone \ac{PERK} scheme, this overhead is well invested additional computational effort.
	Comparing the simulation with $\text{P-ERK}_{3;\{3,4,7\}}$ to the standalone method $\text{P-ERK}_{3;7}$, the memory footprint increases in the mean by $20 \%$.
	
	The ratio of over ten runs averaged relative runtimes $\sfrac{\bar{\tau}}{\bar{\tau}_{\text{P-ERK}_{3;\{3,4,7\}}}}$ of the different schemes are tabulated in \cref{tab:RunTimes_ShearLayer}.
	Also, we provide the degree of optimality $\kappa = \tau_\text{Opt} / \bar{\tau}$ cf. \cref{eq:OptimalityMeasure} and the number of performed scalar \ac{RHS} evaluations $N_\text{RHS}^\text{Actual}$ \cref{eq:ActualRHSEvals} for the different schemes.
	For each method we compare to, the \ac{AMR} interval is adjusted such that approximately the same number of grid refinements is performed as for the simulation with the \ac{PERK} family.
	This is roughly equivalent to performing the mesh refinement at the same physical timestamps.
	When comparing the different schemes it is important to note that, while the simulations are equipped with identical parameters (except for \ac{CFL} number and \ac{AMR} interval $N_\text{AMR}$), the \ac{AMR} indicator which governs local mesh refinement or coarsening is based on the solution itself, which is inevitably different for the different methods.
	This becomes already apparent for this example as the $\text{DGLDD}_{3;7}$ \cite{toulorge2012optimal} scheme is slightly faster than the optimized $\text{P-ERK}_{3;7}$ scheme which can be traced back to the fact that the finest level is reached later for the simulation with the $\text{DGLDD}_{3;7}$ scheme, hence imposing the timestep restriction later.

	Overall, the $\text{P-ERK}_{3;\{3,4,7\}}$ is the fastest method for this testcase with all other methods being at least $40\%$ slower.
	For this example we have a four times as large overall timestep for the $\text{P-ERK}_{3;\{3,4,7\}}$ scheme compared to the $\text{SSP}_{3;3}$ method while observing a speedup of factor $2.356$.
	We observe that the saved number of \ac{RHS} evaluations due to the optimized \ac{PERK} scheme does not directly translate into a speedup of the same factor.
	For instance, the standalone \ac{PERK} scheme requires roughly $76.2 \%$ more scalar \ac{RHS} evaluations, but is only $46.8\%$ slower.
	This is due to fact that Trixi.jl is for the standard DGSEM with a simple numerical flux memory-bound and thus saved computations do not translate proportional into runtime savings. 
	In particular, data locality is in general lost for the partitioned approach, i.e., the coefficients $U^{(r)}$ are in general not stored in a contiguous memory block. 
	Also, the data structures containing elements and interfaces are accessed in a non-sequential manner which leads to cache misses and thus reduced performance.
	In contrast, for the standard methods everything is processed as present in memory which explains the non-ideal performance of the partitioned approach.

	Comparing this to the examples given in \cite{nasab2022third} where the authors observe speedups of about $4.3$ across examples for an eight times as large timestep (compared to the $\text{SSP}_{3;3}$) we conclude that we have essentially the same ratio of speedup per timestep increase.

	The four-stage method $\text{CKLLS}_{3;4}$ \cite{kennedy2000low} is not significantly faster than the canonical Shu-Osher \ac{ssp} scheme, and is thus not considered in the following examples.
	This is possibly due to optimization for Navier-Stokes equations where diffusion plays a more dominant role than in the convection-dominated flows considered here.
	\begin{table}
		\def\arraystretch{1.2}
		\centering
		\begin{tabular}{l?{2}c|c|c}
			Method & $\sfrac{\bar{\tau}}{\bar{\tau}_{\text{P-ERK}_{3;\{3,4,7\}}}}$ & $\sfrac{\tau_\text{Opt}}{\bar{\tau}}$ & $N_\text{RHS}^\text{Actual}$ \\
			\Xhline{5\arrayrulewidth}
			$\text{P-ERK}_{3;\{3,4,7\}}$ 												& $1.0$   & $68.8 \%$ & $5.26 \cdot 10^{9 \hphantom{0}} $ \\
			$\text{P-ERK}_{3;7}$ 																& $1.468$ & $46.9 \%$ & $9.27 \cdot 10^{9 \hphantom{0}}$ \\
			$\text{SSP}_{3;3}$ 																	& $2.356$ & $29.2 \%$ & $1.42 \cdot 10^{10}$ \\
			$\text{DGLDD}_{3;7}$ \cite{toulorge2012optimal} 		& $1.412$ & $48.8 \%$ & $8.53 \cdot 10^{9 \hphantom{0}}$ \\
			$\text{PKD3S}_{3;5}$ \cite{parsani2013optimized} 		& $1.854$ & $37.1 \%$ & $1.07 \cdot 10^{10}$ \\
			$\text{RDPKFSAL}_{3;5}$ \cite{ranocha2022optimized} & $2.171$ & $31.7 \%$ & $1.30 \cdot 10^{10}$ \\
			$\text{CKLLS}_{3;4}$ \cite{kennedy2000low} 					& $2.350$ & $29.3 \%$ & $1.45 \cdot 10^{10}$
		\end{tabular}
		\caption[Runtimes, optimality fraction \cref{eq:OptimalityMeasure}, and number of scalar \ac{RHS} evaluations \cref{eq:ActualRHSEvals} of different third-order \ac{ode} solvers compared to the optimized $p=3, S= \{3, 4, 7\}$ \ac{PERK} scheme for the periodic shear-layer.]
		{Runtimes, optimality fraction \cref{eq:OptimalityMeasure}, and number of scalar \ac{RHS} evaluations \cref{eq:ActualRHSEvals} of different third-order \ac{ode} solvers compared to the optimized $p=3, S= \{3, 4, 7\}$ \ac{PERK} scheme for the periodic shear-layer.}
		\label{tab:RunTimes_ShearLayer}
	\end{table}
	\subsubsection{Visco-Resistive Orszag-Tang Vortex}
	\label{subsec:ViscoResistiveOrszagTang}
	We consider the classic Orszag-Tang vortex \cite{orszag1979small} in the context of the \ac{vrMHD} equations \cite{warburton1999discontinuous}.
	The ideal \ac{mhd} equations are implemented with divergence cleaning through a \ac{GLM} approach \cite{derigs2018ideal, munz2000divergence} and extended by viscous and resistive terms by means of the BR1 scheme \cite{bassi1997high}.

	As in \cite{warburton1999discontinuous} the domain is $\Omega = [0, 2\pi]^2$ and the gas parameters are set to $\gamma = \frac{5}{3}$, $\text{Pr} = 1$.
	Compared to \cite{warburton1999discontinuous} we reduce viscosity and resistivity by one order of magnitude to $\mu = \nu = 10^{-3}$ to stay convection-dominated (cf. \eqref{eq:Convection_Diffusion_Scaling}) even in the finest parts of the mesh.
	The initial condition is given by \cite{warburton1999discontinuous}
	\begin{equation}
		\label{eq:IC_OrszagTang}
		\boldsymbol u_\text{primal}
		(t_0=0, 
		x, y)	= \begin{pmatrix}
			\rho \\ v_x \\ v_y \\ v_z \\ p \\ B_x \\ B_y \\ B_z
		\end{pmatrix}
		= 
		\begin{pmatrix}
			1.0 \\ -2 \sqrt{\pi} \sin(y) \\
			2 \sqrt{\pi} \sin(x) \\ 
			0.0 \\
			\frac{15}{4} + \frac{1}{4} \cos(4 x) + \frac{4}{5} \cos(2x) \cos(y) - \cos(x)\cos(y) + \frac{1}{4} \cos(2 y)\\
			- \sin(2 y) \\ \sin(2 x) \\ 0.0
		\end{pmatrix}
	\end{equation}
	and the simulation is run until $t_f = 2.0$ as in \cite{warburton1999discontinuous}.
	Note that in the \ac{vrMHD} equations employed in \cite{warburton1999discontinuous} the magnetic terms are nondimensionalized using viscous and resistive Lundquist numbers, which is different from \cite{derigs2018ideal} used in this work.
	The grid cells vary between $h = \sfrac{ 2 \pi}{2^5}$ and $h = \sfrac{2 \pi}{2^9}$, i.e., a mesh with five consecutive levels of refinement is employed.

\begin{figure}[!t]
		\centering
		\subfloat[{Pressure $p$ at $t = 1.0$ for the visco-resistive Orszag-Tang vortex.}]{
			\label{fig:OrszagTang_pressure}
			\centering
			\includegraphics[height=.45\textwidth]{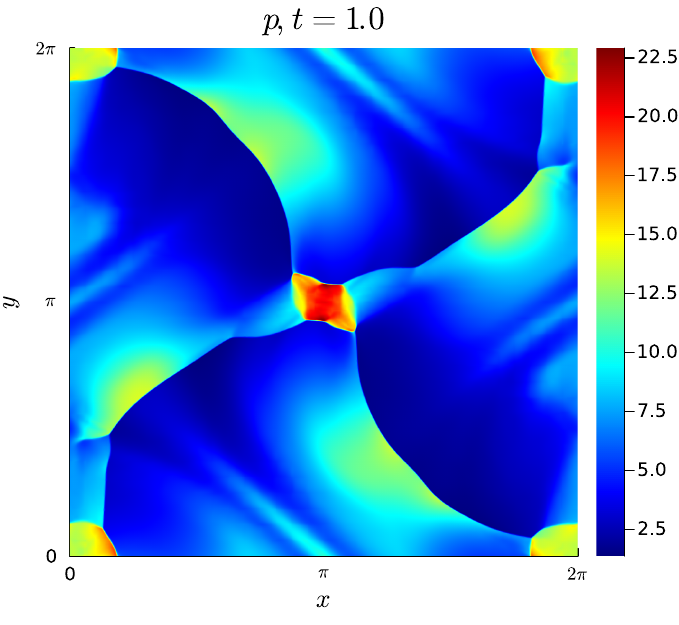}
		}
		\hfill
		\subfloat[{Four times refined mesh at $t = 1.0$.}]{
			\label{fig:OrszagTang_Mesh}
			\centering
			\includegraphics[height=.45\textwidth]{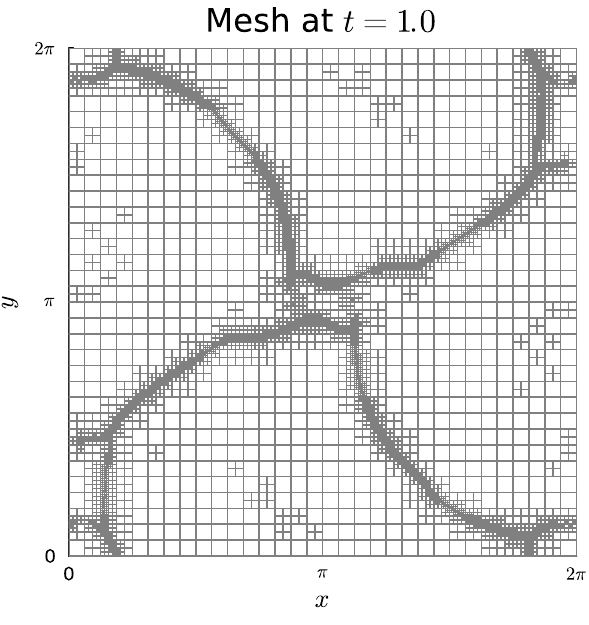}
		}
		\caption[Pressure $p$ at $t = 1.0$ for the visco-resistive Orszag-Tang vortex  and corresponding adaptively refined mesh.]
		{Pressure $p$ at $t = 1.0$ for the visco-resistive Orszag-Tang vortex (\cref{fig:OrszagTang_pressure}) and corresponding adaptively refined mesh (\cref{fig:OrszagTang_Mesh}).}
		\label{fig:OrszagTang_pressure_Mesh}
	\end{figure}

	To discretize the spatial derivatives a split-form \ac{DG} method \cite{gassner2016split} is used.
	The surface-flux itself is split up into conserved terms which are discretized using the Rusanov/Local Lax-Friedrichs flux and a non-conserved part for which a specialized flux is employed \cite{bohm2020entropy, powell1999solution}.
	Also for the volume flux/integral the conservative and non-conserved terms are treated differently.
	For this, we use the central flux for the conserved terms alongside the flux for the nonconservative Powell term \cite{bohm2020entropy}.
	As \ac{AMR} indicator we employ the shock-capturing indicator developed in \cite{hennemann2021provably} based on density-pressure $\rho \cdot p$.
	This particular example benefits heavily from \ac{AMR} since the configuration leads to initially sharp structures cf. \cref{fig:OrszagTang_pressure_Mesh} requiring accurate resolution which are then, due to the viscous effects, damped out over time with less need for high resolution.

	For time discretization, we employ the three-level $p=3, E=\{3, 4, 6\}$ optimized \ac{PERK} family and compare it against the standalone $S=6$ \ac{PERK} scheme, the $\text{SSP}_{3;3}$ method, $\text{DGLDD}_{3;7}$, and  $\text{PKD3S}_{3;5}$ with runtimes provided in \cref{tab:RunTimes_VRMHD_Orszag_Tang}.
	As the ideal GLM-MHD equations require a CFL-based timestep for the divergence cleaning mechanism, we also supply a CFL number to $\text{RDPKFSAL}_{3;5}$ instead of an adaptive error-controlled timestep.
	As for the other cases, the \ac{PERK} family is the fastest method among the considered ones.
	We observe that linear/absolute stability governs in this case the maximum admissible timestep as both the family as well as the standalone optimized \ac{PERK} schemes outperform the other schemes.
	Furthermore, the $\text{SSP}_{3;3}$ scheme is the slowest Runge-Kutta method, which matches the expectations as this scheme is not optimized at all for a special class of semidiscretizations, in contrast to the other methods.
	As for the previous example, we obtain for a four times larger stable timestep for the $\text{P-ERK}_{3;\{3,4,6\}}$ scheme compared to the $\text{SSP}_{3;3}$ a speedup by a factor of $2.005$ which is again the same ratio of speedup per timestep increase demonstrated in the examples in \cite{nasab2022third}.
	\begin{table}
		\def\arraystretch{1.2}
		\centering
		\begin{tabular}{l?{2}c|c|c}
			Method & $\sfrac{\bar{\tau}}{\bar{\tau}_{\text{P-ERK}_{3;\{3,4,6\}}}}$ & $\sfrac{\tau_\text{Opt}}{\bar{\tau}}$ & $N_\text{RHS}^\text{Actual}$ \\
			\Xhline{5\arrayrulewidth}
			$\text{P-ERK}_{3;\{3,4,6\}}$ 												& $1.0$   & $58.2 \%$ & $1.44 \cdot 10^{10}$ \\
			$\text{P-ERK}_{3;6}$ 																& $1.248$ & $46.6 \%$ & $2.07 \cdot 10^{10}$ \\
			$\text{SSP}_{3;3}$ 																	& $2.005$ & $29.0 \%$ & $3.28 \cdot 10^{10}$ \\
			$\text{DGLDD}_{3;7}$ \cite{toulorge2012optimal} 		& $1.575$ & $37.0 \%$ & $2.55 \cdot 10^{10}$ \\
			$\text{PKD3S}_{3;5}$ \cite{parsani2013optimized} 	  & $1.717$ & $33.9 \%$ & $2.72 \cdot 10^{10}$ \\
			$\text{RDPKFSAL}_{3;5}$ \cite{ranocha2022optimized} & $1.649$ & $35.3 \%$ & $2.72 \cdot 10^{10}$ 
		\end{tabular}
		\caption[Runtimes, optimality fraction \cref{eq:OptimalityMeasure}, and number of scalar \ac{RHS} evaluations \cref{eq:ActualRHSEvals} of different third-order integrators compared to the optimized $p=3, S= \{3, 4, 6\}$ \ac{PERK} integrator for the visco-resistive Orszag-Tang vortex.]
		{Runtimes, optimality fraction \cref{eq:OptimalityMeasure}, and number of scalar \ac{RHS} evaluations \cref{eq:ActualRHSEvals} of different third-order integrators compared to the optimized $p=3, S= \{3, 4, 6\}$ \ac{PERK} integrator for the visco-resistive Orszag-Tang vortex.}
		\label{tab:RunTimes_VRMHD_Orszag_Tang}
	\end{table}
	In terms of performance, there are $3.92 \cdot 10^9$ excess scalar \ac{RHS} evaluations which lead to an overall $1.718$ slower-than-optimal \ac{PERK} scheme, i.e., for this case a significant speed-up would be possible if a full range of optimized methods without monotonicity violation could be employed.
	The update of the partition identifiers amounts to $1.0\%$ of the total runtime, that is, the overhead due to the partitioning is again negligible.
	The memory requirements of the $\text{P-ERK}_{3;\{3,4,6\}}$ scheme are on average $28\%$ larger than for the standalone $\text{P-ERK}_{3;6}$ scheme.
	As a sidenote, we note that despite having the same number of scalar \ac{RHS} evaluations, the $\text{RDPKFSAL}_{3;5}$ scheme is slightly faster than the $\text{PKD3S}_{3;5}$ scheme which showcases the relevance of implementation of the timestepping algorithm itself.
	\subsubsection{3D Taylor-Green Vortex}
	The three-dimensional Taylor-Green vortex is a well-known reference problem showcasing emerging turbulence from a simple initial condition which is followed by the decay of the turbulent structures accompanied by kinetic energy dissipation \cite{debonis2013solutions, bull2014simulation}.
	The initial condition in primitive variables is given by
	\begin{equation}
		\label{eq:IC_TaylorGreen}
		\boldsymbol u_\text{prim}(t_0=0, x, y, z) = \begin{pmatrix} \rho \\ v_x \\ v_y \\ v_z \\ p \end{pmatrix} = \begin{pmatrix} 1 \\ \sin(x) \cos(y) \cos(z) \\ -\cos(x) \sin(y) \cos(z) \\ 0 \\ p_0 + \frac{1}{16} \rho \Big(\big(\cos(2x) + \cos(2y)\big) \big( 2+ \cos(2z) \big) \Big) \end{pmatrix} ,
	\end{equation}
	where $p_0 = \frac{\rho}{M^2 \gamma}$ with Mach number $M = 0.1$ as in \cite{jacobs2017opensbli}.
	The Prandtl and Reynolds number are $\text{Pr} = 0.72$ and $\text{Re} = 1600$, and the compressible Navier-Stokes equations are simulated on $\Omega = [0, 2 \pi]^3$ equipped with periodic boundaries until $t_f = 20.0$.
	
	In terms of spatial discretization, we employ solution polynomials of degree $k=2$, HLLE surface flux \cite{einfeldt1988godunov}, and an entropy-conservative volume flux \cite{ranocha2020entropy, gassner2016split}.
	We employ a simple \ac{AMR} indicator based on the enstrophy 
	\begin{equation}
		\label{eq:Enstrophy}
		\epsilon \coloneqq 0.5 \, \rho \, \boldsymbol{\omega} \cdot \boldsymbol{\omega} \, ,
	\end{equation}
	which refines cells once a specified threshold is exceeded.
	Doing so, we expect the mesh to be refined in regions where the vorticity $\boldsymbol{\omega} \coloneqq \nabla \times \boldsymbol{v}$ is high, i.e., where turbulent structures are emerging.
	The initial condition is discretized using $16^3$ cells which are allowed to be two times refined corresponding to a smallest cell size of $h = \frac{2 \pi}{64}$.
	While the Taylor-Green vortex is certainly not a classic use case for \ac{AMR} it allows us to examine the cost of the dynamic re-assignment of cells to partitions for a three-dimensional problem.

	In terms of the temporal integrators, we employ the three-level $p=3, E=\{3, 4, 6\}$ optimized \ac{PERK} family and compare it against the standalone $S=6$ \ac{PERK} scheme, the $\text{SSP}_{3;3}$ method, $\text{DGLDD}_{3;7}$ \cite{toulorge2012optimal}, $\text{PKD3S}_{3;5}$ \cite{parsani2013optimized}, and $\text{RDPKFSAL}_{3;5}$ \cite{ranocha2022optimized} with runtimes provided in \cref{tab:RunTimes_TaylorGreen}.

	When comparing the runtimes of the different schemes in \cref{tab:RunTimes_TaylorGreen} we observe that the $\text{P-ERK}_{3;\{3,4,6\}}$ scheme is again the fastest method among the considered ones, although the standalone $S=6$ \ac{PERK} scheme is only about $4\%$ slower.
	This is not surprising, as for the majority of timesteps most cells are completely refined and thus only very little speed-up may be gained.
	This is also reflected by the speedup that could be gained from a hypothetical ideally scaling set of methods.
	The $5.44 \cdot 10^9$ excess scalar \ac{RHS} evaluations lead only to an overall $1.034$ slower-than-optimal \ac{PERK} scheme.
	We included this example for instructive purposes since it showcases that also (at least for small) 3D problems the overhead of the re-computation of the partitioning datastructures does not exceed speedup gained.
	In particular, only $0.3\%$ of the total runtime is spent on the constuction of the flagging variables after mesh changes.
	In this case, the memory footprint of the $\text{P-ERK}_{3;\{3,4,6\}}$ scheme is on average only $7\%$ larger than the standalone $\text{P-ERK}_{3;6}$ scheme.
	\begin{table}
		\def\arraystretch{1.2}
		\centering
		\begin{tabular}{l?{2}c|c|c}
			Method & $\sfrac{\bar{\tau}}{\bar{\tau}_{\text{P-ERK}_{3;\{3,4,6\}}}}$ & $\sfrac{\tau_\text{Opt}}{\bar{\tau}}$ & $N_\text{RHS}^\text{Actual}$ \\
			\Xhline{5\arrayrulewidth}
			$\text{P-ERK}_{3;\{3,4,6\}}$ 												& $1.0$   & $96.8 \%$ & $1.25 \cdot 10^{11}$\\
			$\text{P-ERK}_{3;6}$ 																& $1.068$ & $90.6 \%$ & $1.42 \cdot 10^{11}$\\
			$\text{SSP}_{3;3}$																	& $1.768$ & $54.7 \%$ & $2.29 \cdot 10^{11}$\\
			$\text{DGLDD}_{3;7}$ \cite{toulorge2012optimal} 		& $1.250$ & $77.4 \%$ & $1.62 \cdot 10^{11}$\\
			$\text{PKD3S}_{3;5}$ \cite{parsani2013optimized} 	  & $1.622$ & $59.7 \%$ & $1.91 \cdot 10^{11}$\\
			$\text{RDPKFSAL}_{3;5}$ \cite{ranocha2022optimized} & $1.425$ & $67.9 \%$ & $1.82 \cdot 10^{11}$
		\end{tabular}
		\caption[Runtimes, optimality fraction \cref{eq:OptimalityMeasure}, and number of scalar \ac{RHS} evaluations \cref{eq:ActualRHSEvals} of different third-order integrators compared to the optimized $p=3, S= \{3, 4, 6\}$ \ac{PERK} integrator for the Taylor-Green vortex.]
		{Runtimes, optimality fraction \cref{eq:OptimalityMeasure}, and number of scalar \ac{RHS} evaluations \cref{eq:ActualRHSEvals} of different third-order integrators compared to the optimized $p=3, S= \{3, 4, 6\}$ \ac{PERK} integrator for the Taylor-Green vortex.}
		\label{tab:RunTimes_TaylorGreen}
	\end{table}
	\subsection{Hyperbolic Equations}
	In addition to the isentropic vortex test case used for verification in \cref{sec:Validation} we present additional applications to inviscid, hyperbolic problems.
	In particular, this is the first time the \ac{PERK} schemes are applied to invsicid problems, apart from the isentropic vortex used in convergence tests.
	For inviscid problems the issues concerning nonlinear stability discussed in \cref{subsec:NonlinearStabPERK} resurface and a simulation with three-level ($R=3$) \ac{PERK} schemes demands in general a reduction of the timestep compared te the optimized standalone scheme.
	Recall that for hyperbolic-parabolic problems the \ac{CFL} number did not have to be reduced for simulations with \ac{PERK} schemes compared to the standalone scheme.
	This clearly impairs the performance of the \ac{PERK} schemes for inviscid problems and correspondingly, we do not see the same speedup as for the viscous problems.
	%
	\subsubsection{Kelvin-Helmholtz Instability}
	\label{subsec:KelvinHelmholtzInstability}
	The inviscid Kelvin-Helmholtz instability is a suitable testcase for the \ac{AMR}-based dynamic \ac{ode} partitioning developed in this work.
	The initial condition can be accurately represented with a quite coarse discretization while the emerging instability requires a significantly finer resolution.
	We use the setup from \cite{rueda2021subcell} where the initial condition is given by 
	\begin{equation}
		\def\arraystretch{1.2}
		\label{eq:KelvinHelmholtzIC}
		\boldsymbol u_\text{primal}(t_0=0, x, y) = \begin{pmatrix} \rho \\ v_x \\ v_y \\ p \end{pmatrix} = \begin{pmatrix} 0.5 + 0.75 B(y) \\ 0.5\big(B(y) - 1\big) \\ 0.1 \sin(2 \pi x) \\ 1 \end{pmatrix}
	\end{equation}
	with $B(y) = \tanh(15y + 7.5) - \tanh(15y - 7.5)$ which corresponds to Mach number $\text{Ma} \leq 0.6$.
	The compressible Euler equations are discretized using the \ac{dgsem} with $k=3$ local polynomials, \ac{hlle} surface-flux \cite{einfeldt1988godunov}, subcell shock-capturing \cite{hennemann2021provably} to avoid spurious oscillations and entropy-conserving volume flux \cite{ranocha2020entropy, fisher2013high}.
	Initially, the domain $\Omega = [-1,1]^2$ is discretized with $2^4 \times 2^4$ cells which are allowed to be refined five times, i.e., the finest cells may have edge length $h = 2^{-8}$.
	The indicator for \ac{AMR} coincides in this case with the shock-capturing indicator \cite{hennemann2021provably} and is based on the density $\rho$.
	For time-integration, the $p=3, E= \{4, 6, 11\}$ \ac{PERK} family is employed where the admissible timestep doubles with increasing stage evaluations $E$.
	In contrast to the viscous flows from the previous section, the \ac{CFL} number of the $\text{P-ERK}_{3;\{4,6,11\}}$ scheme has to be reduced to $70 \%$ of the \ac{CFL} number of the $\text{P-ERK}_{3;11}$ scheme.
	This reduction is a consequence of the troubling nonlinear stability properties of the \ac{PERK} schemes at scheme boundaries, cf. \cref{subsec:NonlinearStabPERK}.

\begin{figure}[!t]
		\centering
		\subfloat[{Density $\rho$ at final time $t_f = 3.0$ for the Kelvin-Helmholtz instability.}]{
			\label{fig:KelvinHelmholtz_density}
			\centering
			\includegraphics[height=.45\textwidth]{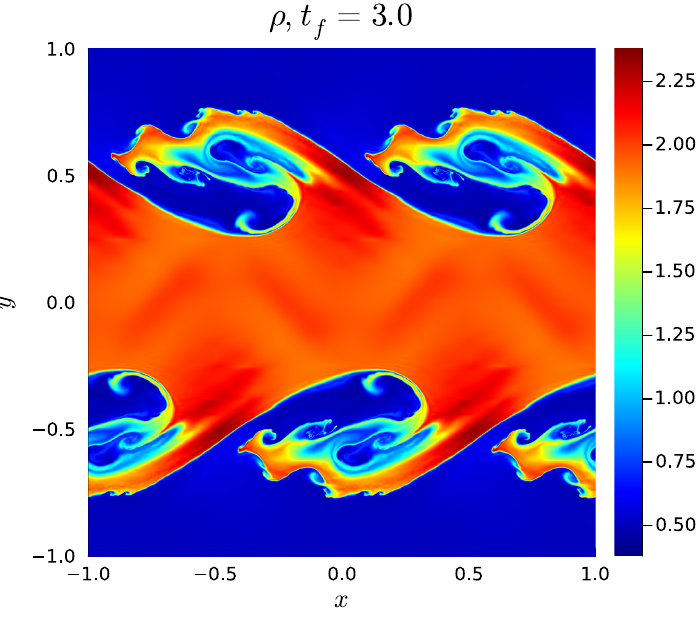}
		}
		\hfill
		\subfloat[{Five times refined mesh at $t_f = 3.0$.}]{
			\label{fig:KelvinHelmholtz_Mesh}
			\centering
			\includegraphics[height=.45\textwidth]{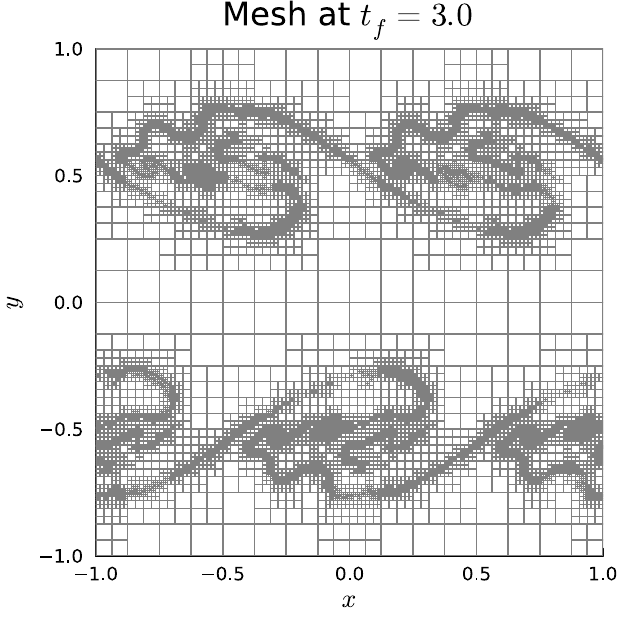}
		}
		\caption[Density and mesh at final time $t_f = 3.0$ for the Kelvin-Helmholtz instability.]
		{Density $\rho$ at final time $t_f = 3.0$ for the Kelvin-Helmholtz instability \cref{eq:KelvinHelmholtzIC} (\cref{fig:KelvinHelmholtz_density}) and corresponding adaptively refined mesh \cref{fig:KelvinHelmholtz_Mesh}.}
		\label{fig:KelvinHelmholtz_density_Mesh}
	\end{figure}

	The computational mesh and solution at final time $t_f=3.0$ are displayed in \cref{fig:KelvinHelmholtz_density_Mesh}.
	As for the viscous problems, we compare the \ac{PERK} family to the stand-alone $p=3, S=11$ scheme and the $\text{SSP}_{3;3}$, $\text{DGLDD}_{3;7}$, $\text{PKD3S}_{3;5}$, and $\text{RDPKFSAL}_{3;5}$ schemes.
	The AMR-interval $N_\text{AMR}$ (cf. \eqref{eq:AdditionalRHSEvals}) is adjusted according to the timesteps undertaken by the different methods to have an approximately equal AMR load.
	The runtimes $\tau$ are averaged over ten runs on a shared-memory multicore thread-parallelized machine and presented in \cref{tab:RunTimes_KelvinHelmholtz}.
	The $\text{P-ERK}_{3;\{4,6,11\}}$ is in this case the most efficient integrator with minimal runtime.
	Notably, the $\text{DGLDD}_{3;7}$ method is only slightly slower and even faster than the optimized $\text{P-ERK}_{3;11}$ scheme.
	This is attributed to favorable nonlinear stability properties (low-dispersion, low-dissipation), which are for this example (and for the inviscid problems in this section overall) more relevant than optimized absolute stability.
	Similarly, the $\text{SSP}_{3;3}$ scheme performs also better than the $\text{PKD3S}_{3;5}$ and $\text{RDPKFSAL}_{3;5}$ schemes despite that these are optimized for \ac{DG} methods and computational fluid dynamics, respectively.
	Note that this was not the case for the previous hyperbolic-parabolic examples, i.e., the $\text{SSP}_{3;3}$ scheme did not outperform any of the optimized schemes for viscous problems.
	\begin{table}[!ht]
		\def\arraystretch{1.2}
		\centering
		\begin{tabular}{l?{2}c|c|c}
			Method & $\sfrac{\bar{\tau}}{\bar{\tau}_{\text{P-ERK}_{3;\{4,6,11\}}}}$ & $\sfrac{\tau_\text{Opt}}{\bar{\tau}}$ & $N_\text{RHS}^\text{Actual}$ \\
			\Xhline{5\arrayrulewidth}
			$\text{P-ERK}_{3;\{4,6,11\}}$ 											& $1.0$   & $92.4 \%$ & $3.78 \cdot 10^9$ \\
			$\text{P-ERK}_{3;11}$ 															& $1.311$ & $70.5 \%$ & $5.30 \cdot 10^9$ \\
			$\text{SSP}_{3;3}$ 																	& $1.426$ & $64.8 \%$ & $6.62 \cdot 10^9$ \\
			$\text{DGLDD}_{3;7}$ \cite{toulorge2012optimal} 		& $1.169$ & $79.1 \%$ & $5.19 \cdot 10^9$ \\
			$\text{PKD3S}_{3;5}$ \cite{parsani2013optimized}		& $1.866$ & $49.4 \%$ & $8.36 \cdot 10^9$ \\
			$\text{RDPKFSAL}_{3;5}$ \cite{ranocha2022optimized} & $1.659$ & $55.7 \%$ & $6.54 \cdot 10^9$
		\end{tabular}
		\caption[Runtimes, optimality fraction \cref{eq:OptimalityMeasure}, and number of scalar \ac{RHS} evaluations \cref{eq:ActualRHSEvals} of different third-order integrators compared to the optimized $p=3, S= \{4, 6, 11\}$ \ac{PERK} integrator for the inviscid Kelvin-Helmholtz instability.]
		{Runtimes, optimality fraction \cref{eq:OptimalityMeasure}, and number of scalar \ac{RHS} evaluations \cref{eq:ActualRHSEvals} of different third-order integrators compared to the optimized $p=3, S= \{4, 6, 11\}$ \ac{PERK} integrator for the inviscid Kelvin-Helmholtz instability.}
		\label{tab:RunTimes_KelvinHelmholtz}
	\end{table}

	The update of the identifier variables required for the partitioning amounts in this case to about $0.4\%$ of the total runtime $\tau_\text{P-ERK}$.
	The memory footprint of the $\text{P-ERK}_{3;\{4,6,11\}}$ scheme is on average $23\%$ larger than the standalone $\text{P-ERK}_{3;11}$ scheme.
	To assess the optimality, or, in other words, the room for further improvements of using, e.g., a \ac{PERK} family with more members and higher stages, we compute the excess scalar \ac{RHS} evaluations according to \eqref{eq:AdditionalRHSEvals}.
	For this example, there are $3.80 \cdot 10^8$ additional scalar \ac{RHS} evaluations per conserved variable which increase the runtime of the employed $\text{P-ERK}_{3;\{4,6,11\}}$ scheme to $1.035$ of the optimal one, neglecting overhead due to additional levels.
	This value results from the fact that the vast majority of the cells is distributed among the finest three levels.
	In particular, at final time $t_f = 3.0$ the three non-optimal integrated levels form only $6\%$ of the entire number of cells, cf. \cref{fig:KelvinHelmholtz_Mesh}.
	In other words, the majority of the cells are integrated optimally with results in a speed-up of $30\%$ compared to the standalone $\text{P-ERK}_{3;11}$ scheme.

	When comparing the number of performed scalar \ac{RHS} evaluations $N_\text{RHS}^\text{Actual}$ to the runtime $\tau$ we observe not exactly a linear correspondence.
	This is due to two factors: First, the \ac{AMR} contributes a fixed-cost to the overall measured runtime $\tau$ which does not decrease with the number of scalar \ac{RHS} evaluations.
	Second, and more signficantly, the employed shock-capturing mechanism \cite{hennemann2021provably} leads to different computational costs for different cells, as for some a plain \ac{dgsem} reconstruction is employed while for others a blending of \ac{dgsem} and subcell Finite Volume reconstruction is necessary.
	The influence of this can clearly be seen by comparing the runtimes of the standalone $\text{P-ERK}_{3;11}$ scheme and the $\text{DGLDD}_{3;7}$ scheme which have roughly the same number of scalar \ac{RHS} evaluations.
	The $\text{DGLDD}_{3;7}$ scheme is, however, over-proportionally faster than the $\text{P-ERK}_{3;11}$ scheme, which is attributed to the lower computational cost of the employed shock-capturing mechanism.
	Similarly, the speedup of the $\text{P-ERK}_{3;\{4,6,11\}}$ scheme compared to the $\text{P-ERK}_{3;11}$ scheme is not exactly proportional to the number of scalar \ac{RHS} evaluations, again due to increased computational cost of the employed shock-capturing mechanism.
	This can be seen as a price of the poor nonlinear stability properties of the \ac{PERK} schemes which cause more cells to rely on a blend of \ac{dgsem} and subcell Finite Volume reconstruction.
	\subsubsection{Rayleigh-Taylor Instability}
	\label{subsec:RayleighTaylorInstability}
	Another classic example 
	is the Rayleigh-Taylor instability, where a heavier fluid gets pushed into a lighter one, leading to complicated flow structures beneath the characteristic mushroom cap.
	The configuration follows \cite{shi2003resolution} with reflective boundary conditions on the boundaries in $x$-direction and weakly enforced Dirichlet boundary conditions in $y$-direction.
	The domain $\Omega = [0, 0.25] \times [0, 1]$ is initially discretized with $N_x = 12, N_y = 4 \cdot N_x = 48$ third-order ($k=3$) elements.
	The spatial discretization is analogous to the Kelvin-Helmholtz instability from \cref{subsec:KelvinHelmholtzInstability} and the simulation is run up to $t_f = 1.95$ as in \cite{shi2003resolution}.
	Following \cite{shi2003resolution}, we set the ratio of specific heats $\gamma = \sfrac{5}{3}$ and the initial state to 
	\begin{equation}
		\label{eq:RTI_IC}
		\boldsymbol u_\text{prim}(t_0=0, x, y) = \begin{pmatrix} \rho \\ v_x \\ v_y \\ p \end{pmatrix} = \begin{pmatrix} \begin{cases} 2.0 & y \leq 0.5 \\ 1.0 & y > 0.5 \end{cases} \\ 0.0 \\ c k \cos(8 \pi x) \\ \begin{cases} 2.0 \cdot y + 1.0 & y \leq 0.5 \\ y + 1.5 & y > 0.5 \end{cases}  \end{pmatrix} 
	\end{equation}
	where $ c = \sqrt{\frac{\gamma p}{\rho}}$ denotes the speed of sound and $k= -0.025$ quantifies the strength of the initial perturbation.
	
	Each cell is allowed to be refined at most six times, i.e., the theoretically maximum resolution is $768 \times 3072$ with $h = \sfrac{1}{3072}$.
	For time integration we employ a $p=3, E = \{3, 4, 6\}$ scheme where the maximum admissible timestep doubles for each family member.
	Thus, the finest three levels can be integrated optimally and only the remaining levels are not efficiently integrated.
	Compared to the standalone $\text{P-ERK}_{3;6}$ scheme the \ac{CFL} number of the $\text{P-ERK}_{3;\{3,4,6\}}$ is reduced to $96 \%$ of the optimal value.
	We observe, however, that the $\text{P-ERK}_{3;4}$ has a higher stage-normalized \ac{CFL} number $\text{CFL}/E$, rendering it more efficient.
	Consequently, we employ $\text{P-ERK}_{3;4}$ for the runtime comparisons in \cref{tab:RunTimes_RayleighTaylor}.
	The computation of the partitioning variables consumes about $4.4\%$ of the total runtime.
	This increase compared to the other examples is due to the usage of the in principle unstructured \texttt{p4est} mesh which requires the computation of the minimum edge length for each cell.
	This is additional effort compared to the necessarily square/cubic cells forming the \texttt{TreeMesh}, where each cells' characteristic length is fully determined by the position of the associated node in the underlying tree datastructure.
	In terms of additional memory consumption, the $\text{P-ERK}_{3;\{3,4,6\}}$ scheme requires on average $34\%$ more memory than the standalone $\text{P-ERK}_{3;4}$ scheme.
	\begin{figure}[!t]
		\centering
		\subfloat[{Density $\rho$ at final time $t_f = 1.95$ for the Rayleigh-Taylor instability.}]{
			\label{fig:RayleighTaylor_density}
			\centering
			\includegraphics[height=.4\textheight]{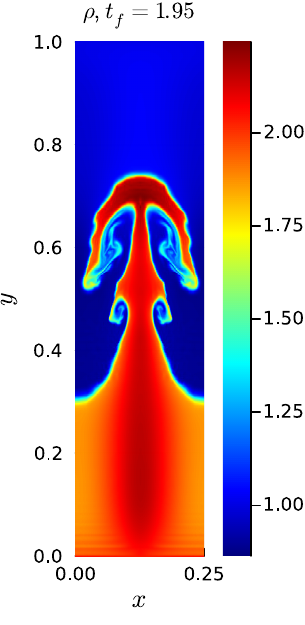}
		}
		\hspace{5em}
		\subfloat[{Five times refined mesh at $t_f = 1.95$.}]{
			\label{fig:RayleighTaylor_Mesh}
			\centering
			\includegraphics[height=.4\textheight]{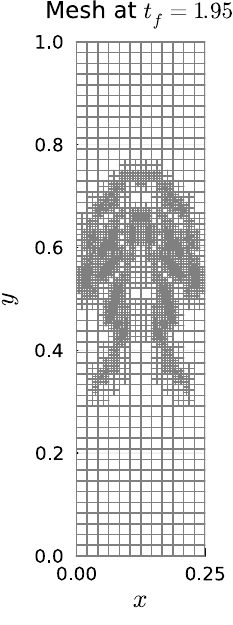}
		}
		\caption[Density and mesh at final time $t_f = 1.95$ for the Rayleigh-Taylor instability.]
		{Density $\rho$ at final time $t_f = 1.95$ for the Rayleigh-Taylor instability \cref{eq:RTI_IC} (\cref{fig:RayleighTaylor_density}) and corresponding adaptively refined mesh \cref{fig:RayleighTaylor_Mesh}.}
		\label{fig:RayleighTaylor_density_Mesh}
	\end{figure}

	As for the Kelvin-Helmholtz instability, the $\text{P-ERK}_{3;\{3,4,6\}}$ scheme is the fastest method among the considered ones, see \cref{tab:RunTimes_RayleighTaylor} followed by the standalone optimized $\text{P-ERK}_{3;4}$ scheme.
	Notably, the usually well-performing $\text{DGLDD}_{3;7}$ \cite{toulorge2012optimal} scheme is in this case relatively slow.
	As for the Kelvin-Helmholtz instability, we observe a discrepancy between the number of scalar \ac{RHS} evaluations $N_\text{RHS}^\text{Actual}$ and the runtime $\tau$.
	For instance, the $\text{P-ERK}_{3;\{3,4,6\}}$ scheme has roughly half the number of scalar \ac{RHS} evaluations of the $\text{PKD3S}_{3;5}$ scheme.
	However, the $\text{P-ERK}_{3;\{3,4,6\}}$ scheme is only about $1.4$ times faster than the $\text{PKD3S}_{3;5}$ scheme.
	This is attributed to more cells relying on the expensive subcell shock-capturing mechanism \cite{hennemann2021provably} for the $\text{P-ERK}_{3;\{3,4,6\}}$ scheme than for the $\text{PKD3S}_{3;5}$ scheme.
	Additionally, the \ac{AMR} fixed-costs which are for the $\text{P-ERK}_{3;\{3,4,6\}}$ scheme about $23 \%$ of the total runtime, do not decrease if $N_\text{RHS}^\text{Actual}$ is reduced.
	\begin{table}
		\def\arraystretch{1.2}
		\centering
		\begin{tabular}{l?{2}c|c|c}
			Method & $\sfrac{\bar{\tau}}{\bar{\tau}_{\text{P-ERK}_{3;\{3,4,6\}}}}$ & $\sfrac{\tau_\text{Opt}}{\bar{\tau}}$ & $N_\text{RHS}^\text{Actual}$ \\
			\Xhline{5\arrayrulewidth}
			$\text{P-ERK}_{3;\{3,4,6\}}$ 												& $1.0$   & $82.5 \%$ & $1.25 \cdot 10^{10}$ \\
			$\text{P-ERK}_{3;4}$ 																& $1.126$ & $73.3 \%$ & $1.84 \cdot 10^{10}$ \\
			$\text{SSP}_{3;3}$ 																	& $1.205$ & $68.5 \%$ & $2.04 \cdot 10^{10}$ \\
			$\text{DGLDD}_{3;7}$ \cite{toulorge2012optimal} 		& $1.751$ & $47.1 \%$ & $3.32 \cdot 10^{10}$ \\
			$\text{PKD3S}_{3;5}$ \cite{parsani2013optimized} 		& $1.402$ & $58.9 \%$ & $2.40 \cdot 10^{10}$ \\
			$\text{RDPKFSAL}_{3;5}$ \cite{ranocha2022optimized} & $1.292$ & $63.9 \%$ & $1.44 \cdot 10^{10}$
		\end{tabular}
		\caption[Runtimes, optimality fraction \cref{eq:OptimalityMeasure}, and number of scalar \ac{RHS} evaluations \cref{eq:ActualRHSEvals} of different third-order integrators compared to the optimized $p=3, S= \{3, 4, 6\}$ \ac{PERK} integrator for the inviscid Rayleigh-Taylor instability.]
		{Runtimes, optimality fraction \cref{eq:OptimalityMeasure}, and number of scalar \ac{RHS} evaluations \cref{eq:ActualRHSEvals} of different third-order integrators compared to the optimized $p=3, S= \{3, 4, 6\}$ \ac{PERK} integrator for the inviscid Rayleigh-Taylor instability.}
		\label{tab:RunTimes_RayleighTaylor}
	\end{table}
	\subsubsection{MHD Rotor}
	\label{sec:MHDRotor}
	In order to go beyond the Euler equations of gas dynamics we consider the ideal \ac{mhd} equations with divergence cleaning through a \ac{GLM} approach \cite{derigs2018ideal, munz2000divergence}.
	A classic testcase thereof is the \ac{mhd} rotor \cite{derigs2018entropy} which we also study here.
	The somewhat lenghty initial condition is provided in \ref{sec:MHDRotorIC}.
	
\begin{figure}[!t]
		\centering
		\subfloat[{Pressure $p$ at final time $t_f = 0.15$ for the ideal MHD rotor.}]{
			\label{fig:MHD_Rotor_pressure}
			\centering
			\includegraphics[height=.45\textwidth]{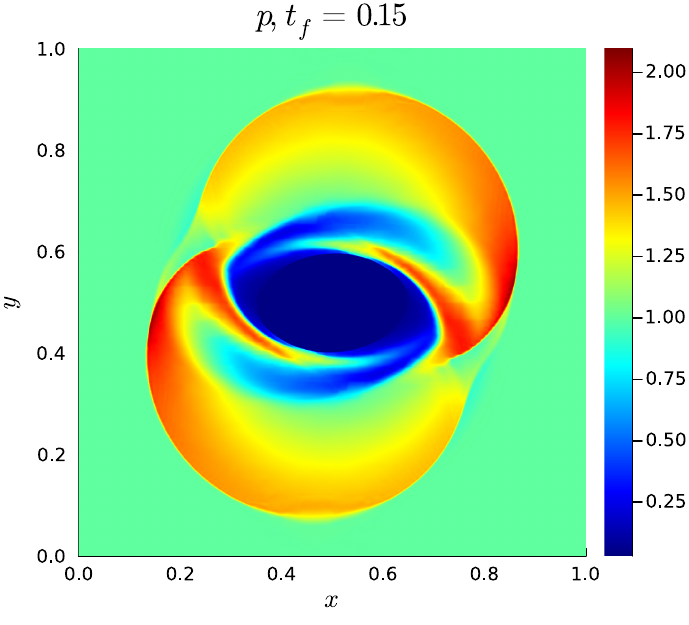}
		}
		\hfill
		\subfloat[{Six times refined mesh at $t_f = 0.15$.}]{
			\label{fig:MHD_Rotor_Mesh}
			\centering
			\includegraphics[height=.45\textwidth]{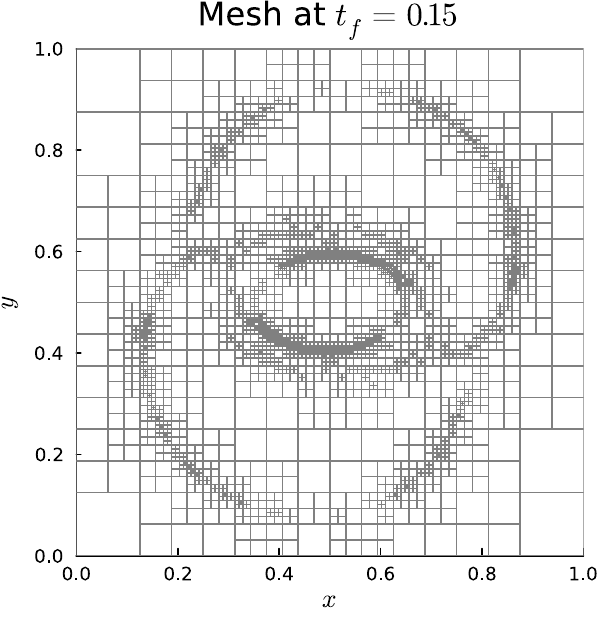}
		}
		\caption[Pressure $p$ and mesh at final time $t_f = 0.15$ for the MHD rotor.]
		{Pressure $p$ at final time $t_f = 0.15$ for the MHD rotor \cite{derigs2018entropy} (\cref{fig:MHD_Rotor_pressure}) and corresponding adaptively refined mesh \cref{fig:MHD_Rotor_Mesh}.}
		\label{fig:MHD_Rotor_pressure_Mesh}
	\end{figure}

	In terms of spatial discretization, we employ $k=4$ local polynomials with the fluxes being discretized analogous to the inviscid fluxes from the visco-resistive Orszag-Tang vortex studied in \cref{subsec:ViscoResistiveOrszagTang}.
	Again, for \ac{AMR} indication the shock-capturing mechanism from \cite{hennemann2021provably} based on density-pressure $\rho \cdot p$ is employed.
	The pressure distribution at final time $t_f = 0.15$ is displayed in \cref{fig:MHD_Rotor_pressure_Mesh} together with the corresponding mesh.
	For this example the $p=3$, $E = \{4, 6, 10\}$ scheme is a suitable choice for the \ac{PERK} family where the maximum admissible timestep doubles for each family member.
	Similar to the Rayleigh-Taylor instability, the \ac{CFL} number of the $\text{P-ERK}_{3;\{4,6,10\}}$ scheme has not to be significantly reduced compared to the standalone $\text{P-ERK}_{3;10}$ scheme.
	We compare the \ac{PERK} family to the best performing stand-alone $\text{P-ERK}_{3;6}$ optimized method, $\text{SSP}_{3;3}$, $\text{DGLDD}_{3;7}$, $\text{PKD3S}_{3;5}$, and $\text{RDPKFSAL}_{3;5}$ (with CFL-determined timestep).
	The results are tabulated in \cref{tab:RunTimes_MHDRotor}.
	As for the previous examples, the \ac{PERK} family is the fastest method among the considered ones with the standalone $S=6$ \ac{PERK} scheme being about $10\%$ slower.
	The runtimes for the $\text{PKD3S}_{3;5}$ and $\text{RDPKFSAL}_{3;5}$ scheme are almost identical which is not surprising since the scheme $\text{RDPKFSAL}_{3;5}$ developed in \cite{ranocha2022optimized} is based on the $\text{PKD3S}_{3;5}$ scheme \cite{parsani2013optimized}.
	\begin{table}
		\def\arraystretch{1.2}
		\centering
		\begin{tabular}{l?{2}c|c|c}
			Method & $\sfrac{\bar{\tau}}{\bar{\tau}_{\text{P-ERK}_{3;\{4,6,10\}}}}$ & $\sfrac{\tau_\text{Opt}}{\bar{\tau}}$ & $N_\text{RHS}^\text{Actual}$ \\
			\Xhline{5\arrayrulewidth}
			$\text{P-ERK}_{3;\{4,6, 10\}}$ 											& $1.0$   & $88.1 \%$ & $2.35 \cdot 10^9$ \\
			$\text{P-ERK}_{3;6}$ 																& $1.097$ & $80.6 \%$ & $3.23 \cdot 10^9$ \\
			$\text{SSP}_{3;3}$ 																	& $1.492$ & $57.9 \%$ & $4.43 \cdot 10^9$\\
			$\text{DGLDD}_{3;7}$ \cite{toulorge2012optimal} 		& $1.148$ & $76.7 \%$ & $3.51 \cdot 10^9$ \\
			$\text{PKD3S}_{3;5}$ \cite{parsani2013optimized} 		& $1.195$ & $73.7 \%$ & $3.56 \cdot 10^9$ \\
			$\text{RDPKFSAL}_{3;5}$ \cite{ranocha2022optimized} & $1.192$ & $73.9 \%$ & $3.56 \cdot 10^9$
		\end{tabular}
		\caption[Runtimes, optimality fraction \cref{eq:OptimalityMeasure}, and number of scalar \ac{RHS} evaluations \cref{eq:ActualRHSEvals} of different third-order integrators compared to the optimized $p=3, S= \{4, 6, 10\}$ \ac{PERK} integrator for the ideal MHD rotor.]
		{Runtimes, optimality fraction \cref{eq:OptimalityMeasure}, and number of scalar \ac{RHS} evaluations \cref{eq:ActualRHSEvals} of different third-order integrators compared to the optimized $p=3, S= \{4, 6, 10\}$ \ac{PERK} integrator for the ideal MHD rotor.}
		\label{tab:RunTimes_MHDRotor}
	\end{table}
	For this case, a fully optimal \ac{PERK} scheme could save $7.71 \cdot 10^7$ scalar \ac{RHS} calls per field.
	Correspondingly, the here employed three-level \ac{PERK} scheme is $1.135$ times slower than a hypothetical optimal \ac{PERK} family, again neglecting potential overhead due to the usage of additional levels.
	The dynamic updating of the partitioning indicator variables amounts in this case to $0.8\%$ of the total runtime.
	The memory footprint of the $\text{P-ERK}_{3;\{4,6,10\}}$ scheme is on average $11\%$ larger than the standalone $\text{P-ERK}_{3;6}$ scheme.
	As for the other two inviscid examples, the decrease in runtime $\tau$ of the \ac{PERK} schemes is not proportional to the number of scalar \ac{RHS} evaluations $N_\text{RHS}^\text{Actual}$.
	The reason behind this is that for the $\text{P-ERK}_{3;\{4,6,10\}}$ scheme a larger number of cells require the subcell stabilizing mechanism \cite{hennemann2021provably} compared to the standalone schemes.
	This is due to the poor nonlinear stability properties of the \ac{PERK} schemes and thus an increased computational cost per \ac{RHS} evaluation.
	\section{Conclusions}
	\label{sec:Conclusions}
	In this work, we applied the Paired-Explicit Runge-Kutta (P-ERK) schemes by Vermeire et. al. \cite{vermeire2019paired,nasab2022third} to convection-dominated, compressible flows with adaptive mesh refinement.
	For a range of hyperbolic and hyperbolic-parabolic examples we demonstrate speed-up compared to state-of-the-art optimized Runge-Kutta methods.
	In particular, the \ac{PERK} outperform for every considered case the non-multirate schemes by factors up to $2.4$ for viscous problems and $1.8$ for inviscid problems.

	Additional data structures are required for the partitioning, which can for all considered examples be set up in less of $5\%$ of the total runtime.
	The memory footprint of partitioning variables leads to $7-34\%$ additional memory consumption compared to the standalone schemes.
	For the shown applications a public reproducibility repository is provided \cite{doehring2024multirateRepro}.

	In general, the \ac{PERK} schemes perform better for \acp{PDE} involving viscous fluxes compared to purely inviscid problems.
	This is due to the fact that the \ac{PERK} schemes are in general not monotonicity-preserving and thus require the stabilizing shock-capturing strategy \cite{hennemann2021provably} on a larger number of cells than the standalone schemes.
	We investigated the mechanism behind this behaviour by conducting a thorough nonlinear stability analysis of the fully discrete, partitioned system.
	In particular, we show that the usage of the \ac{PERK} schemes leads in general to a loss of the monotonicity properties of the spatial discretization.
	More precise, the monotonicity properties are violated at the intersection of different schemes, i.e., where the partitioning is performed.
	This restricts the \ac{PERK} schemes to configurations with limited difference in stage evaluations among partitions for an efficient time integration.

	Future work will focus on the constuction of fourth-order consistent \ac{PERK} schemes and the application to unstructured meshes with adaptive mesh refinement.
	\section*{Data Availability}
	All data generated or analyzed during this study are included in this published article and its supplementary information files.
	\section*{Code Availability \& Reproducibility}
	We provide a reproducibility repository \cite{doehring2024multirateRepro} containing the source code and scripts and data to reproduce the results from \cref{sec:Validation} and \cref{sec:Applications}.
	\section*{Acknowledgments}
	Funding by German Research Foundation (DFG) under Research Unit FOR5409: 
	\\ "Structure-Preserving Numerical Methods for Bulk- and Interface-Coupling of Heterogeneous Models ~(SNuBIC)~" (grant \#463312734).
	\\
	The authors thank Lambert Theisen for fruitful discussion concerning the iterative spectrum computation.
	%
	%
	\section*{Declaration of competing interest}
	The authors declare the following financial interests/personal relationships which may be considered as potential competing interests:
	Daniel Doehring financial support was provided by German Research Foundation.

	\section*{CRediT authorship contribution statement}
	\noindent
	\textbf{Daniel Doehring}: Formal analysis, Investigation, Methodology, Software, Writing - original draft. \\
	\textbf{Michael Schlottke-Lakemper}: Conceptualization, Software, Funding acquisition, Writing - review \& editing. \\
	\textbf{Gregor J. Gassner}: Conceptualization, Funding acquisition, Writing - review \& editing. \\
	\textbf{Manuel Torrilhon}: Conceptualization, Funding acquisition, Supervision, Writing - review \& editing.

	\appendix

	\section{P-ERK Schemes on a Smoothly Refined, High-Resolution Mesh}
	\label{sec:AppSmoothHighResMesh}
\begin{figure}[!t]
		\centering
		\includegraphics[width=.7\textwidth]{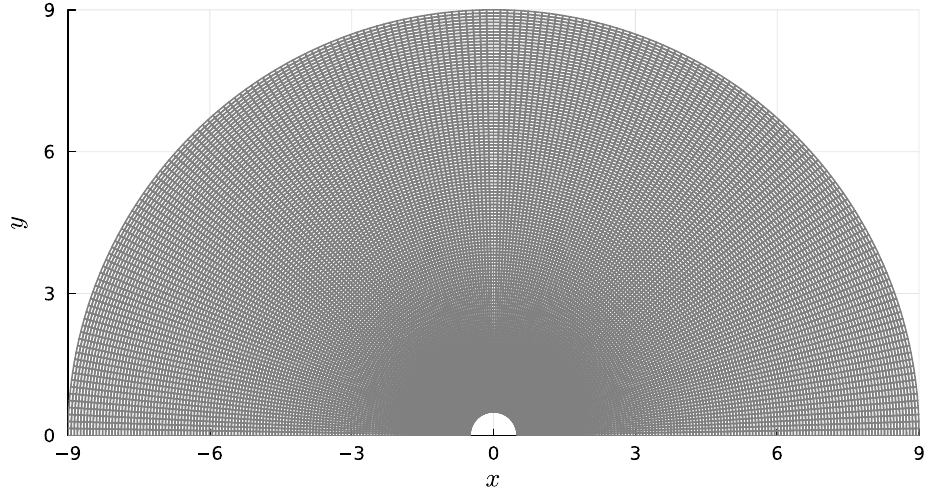}
		\caption[High-resolution, structured mesh for demonstration of working \ac{PERK} schemes on smooth meshes.]{High-resolution, structured mesh for demonstration of working \ac{PERK} schemes on smooth meshes.}
		\label{fig:StructuedTest_Mesh}
	\end{figure}

	In order to give further evidence for the nonlinear stability analysis conducted in \cref{subsec:NonlinearStabPERK}, we present a testcase where the \ac{PERK} schemes are applied to a smooth, high-resolution mesh.
	The refinement factor $\alpha$ of neighboring cells is for this example close to $1$, i.e., the discontinuities in the mesh are substantially less pronounced than for the quad/octree-based \ac{AMR} meshes.
	To this end, we construct a structured curved mesh by defining the mapping of the edges of the cartesian square $(\xi, \eta) \in [-1, 1]^2$ onto curved edges $\boldsymbol e_i, i = 1, \dots, 4$ of the computational domain.
	In particular, we employ the mapping
	\begin{subequations}
		\begin{align}
			\boldsymbol e_1 (\xi, \eta) &= \begin{pmatrix}
				r_0 + 0.5 (r_1 - r_0) (\xi + 1) \\ 0
			\end{pmatrix} \\
			\boldsymbol e_2 (\xi, \eta) &= \begin{pmatrix}
				r_1 \cos\big(0.5 \pi (\eta + 1)\big) \\ r_1  \sin\big(0.5 \pi (\eta + 1)\big)
			\end{pmatrix} \\
			\boldsymbol e_3 (\xi, \eta) &= \begin{pmatrix}
				-r_0 - 0.5 (r_1 - r_0) (\xi + 1) \\ 0
			\end{pmatrix} \\
			\boldsymbol e_4 (\xi, \eta) &= \begin{pmatrix}
				r_0 \cos\big(0.5 \pi (\eta + 1)\big) \\ r_0 \sin\big(0.5 \pi (\eta + 1)\big)
			\end{pmatrix}
		\end{align}
	\end{subequations}
	where the edges $\boldsymbol e_i$ of the $[-1, 1]^2$ cube are labeled counter-clockwise starting from the bottom $(\xi,\eta = 0)$ edge.
	This mapping corresponds to a half-circle with radius $r_1$ from which a half-circle with radius $r_0$ is subtracted.
	The mesh is discretized using $N_\xi = 192 \times N_\eta = 128$ cells and is displayed in \cref{fig:StructuedTest_Mesh}.

\begin{figure}[!t]
		\centering
		\subfloat[{Initial pressure perturbation according to \cref{eq:StructuedTest_IC}.}]{
			\label{fig:StructuedTest_p0}
			\centering
			\includegraphics[width=.45\textwidth]{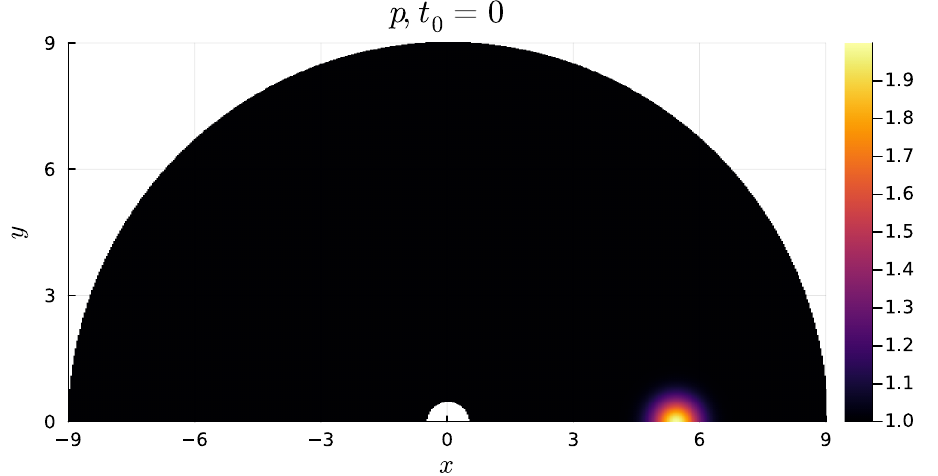}
		}
		\hfill
		\subfloat[{Pressure $p$ at final time $t_f = 5.0$.}]{
			\label{fig:StructuedTest_pf}
			\centering
			\includegraphics[width=.45\textwidth]{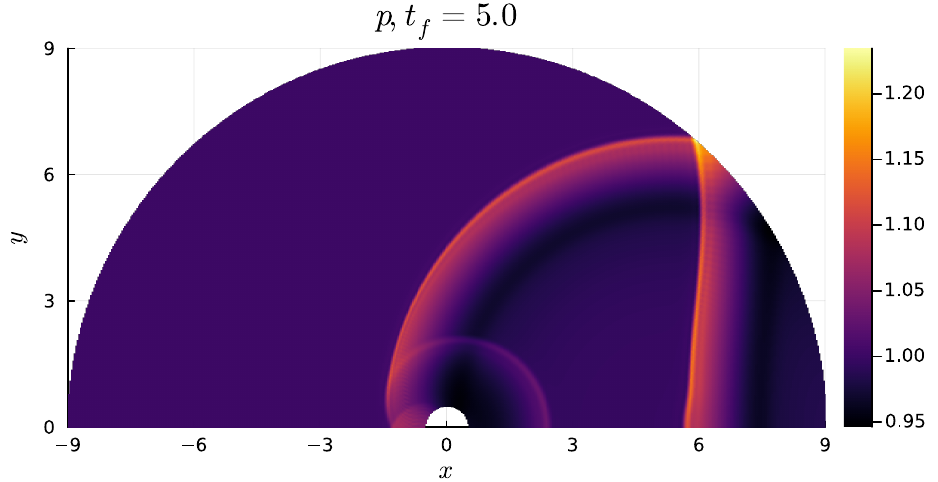}
		}
		\caption[Initial and final pressure distribution on the smoothly refined mesh.]{Initial pressure perturbation according to \cref{eq:StructuedTest_IC} (\cref{fig:StructuedTest_p0}) and pressure $p$ at final time $t_f = 5.0$ (\cref{fig:StructuedTest_pf}).}
		\label{fig:StructuedTest_p}
	\end{figure}

	For this testcase we simulate the inviscid Euler equations of gas dynamics with $\gamma = 1.4$ with initial condition
	\begin{equation}
		\label{eq:StructuedTest_IC}
		\boldsymbol u_\text{prim}(t_0=0, x, y) = \begin{pmatrix} \rho \\ v_x \\ v_y \\ p \end{pmatrix} = \begin{pmatrix} 1.0 \\ 0.0 \\ 0.0 \\ 1.0 + 2^{-\frac{(x-0.6 r_1)^2 + y^2}{0.45^2}} \end{pmatrix}	
	\end{equation}
	which corresponds to a pressure perturbation with magnitude $p_\text{max} = 2$ located at the right-bottom of the half-circle, cf. \cref{fig:StructuedTest_p0}.

	In terms of spatial discretization we employ the \ac{dgsem} with $k=2$ local polynomials, \ac{hll} surface-flux and an entropy-conservative volume flux \cite{ranocha2020entropy, gassner2016split} with reflective boundaires.
	In particular, no techniques like shock capturing or positivity-preserving limiters are employed to suppress spurious oscillations.

	For time-integration, second-order accurate stability polynomials with degrees
	$S= \allowbreak \{4, \allowbreak 6, \allowbreak 8, \allowbreak \dots, \allowbreak 16\}$ 
	are optimized, based on which we construct a range of \ac{PERK} schemes.
	The gridcells are sorted into $R$ (number of \ac{PERK} family members) bins based on their minimum edge length which govern the associated Runge-Kutta method.
	We test all combinations of \ac{PERK} schemes with $E^{(1)}=4$ minimum stage evaluations, i.e., from the standalone $S=4$ scheme up to the $S=16, E=\{4, 6, 8, \dots, 16\}$ family.
	For all combinations, the maximum admissible timestep is increased according to highest number of stage evaluations $E^{(R)}$.
	For this setup, across different \ac{PERK} schemes, the stage-evaluation normalized CFL-number for stable simulation is essentially constant across \ac{PERK} families, as shown in \cref{tab:CFL_EMax_StructuredTest}.
	For the sake of completeness we also provide the speedup compared to the Shu-Osher $p=3, S=3$ scheme, the reference method from \cite{vermeire2019paired,nasab2022third} in \cref{tab:CFL_EMax_StructuredTest}.
	\begin{table}[!ht]
		\def\arraystretch{1.5}
		\centering
		\begin{tabular}{c?{1.5}c|c|c|c|c|c|c}
			$E^{(R)}$ & 4 & 6 & 8 & 10 & 12 & 14 & 16 \\
			\Xhline{3.5\arrayrulewidth}
			$\text{CFL} \cdot \frac{E^{(1)}}{E^{(R)}}$ & $0.95$ & $0.97$ & $0.97$ & $0.94$ & $0.94$ & $0.94$ & $0.93$ \\
			\hline 
			$\frac{\tau_{\text{SSP}_{3;3}}}{\tau_{\text{P-ERK}_{2;4, \dots, E^{(R)}}}}$ & $1.26$ & $1.64$ & $1.90$ & $2.06$ & $2.22$ & $2.37$ & $2.45$
		\end{tabular}
		\caption[]{Normalized CFL numbers and speedup across different \ac{PERK} schemes for the structured mesh testcase. Here, $E^{(1)} = 4$ is the base scheme utilized in all cases.}
		\label{tab:CFL_EMax_StructuredTest}
	\end{table}
	\\
	For illustration, the final pressure distribution $p$ at $t_f = 5.0$ is displayed in \cref{fig:StructuedTest_p} alongside the initial perturbation according to \eqref{eq:StructuedTest_IC}.

	We repeated this testcase for a significantly coarser mesh with $N_\xi = 48 \times N_\eta = 32$ cells and obtained qualitatively similar results.
	In particular, the normalized CFL number is again of similar value, with the largest value $\text{CFL} = 1.04$ for the \ac{PERK} scheme with $E^{(R)} = 10$ and smallest value for $E^{(R)} = 16$ with $\text{CFL} = 0.95$.

	This example demonstrates that the \ac{PERK} schemes are a suitable multirate scheme, even for inviscid problems without stabilization.
	In particular, due to the smooth meshes, the discontinuity in the mesh is effectively eliminated and does not coincide with the discontinuity in the time integration method.
	This is the main difference to quad/octree-based\ac{AMR} meshes where the discontinuity in the mesh coincides with the discontinuity in the time integration method.
	In this case, the mesh refinement factor is always $\alpha = 2$ and the difference in stage evaluations scales quadratically with the number of levels, leading to a discontinuity in methods for the finest levels.
	\section{Ideal MHD Rotor: Initial Condition}
	\label{sec:MHDRotorIC}
	For the simulation of the ideal MHD rotor we employ the initial condition from \cite{derigs2018entropy}.
	The radius $r_0$ of the spinning cylinder is set to $r_0 = 0.1$ while the medium is initially at rest for $r \geq r_1 = 0.115$ with radius $r(x, y) = \sqrt{(x-0.5)^2 + (y-0.5)^2}$.
	The magnetic field is set uniformly to $\boldsymbol B(t_0 = 0, x, y) = \begin{pmatrix} \frac{5}{2 \sqrt{\pi} } & 0 & 0\end{pmatrix}$ as well as pressure $p(0, x, y) = 1$ and $v_z(0, x, y) = 0$.
	Inside the spinning cylinder, i.e., $r \leq r_0$ we set the density to $\rho(0, x, y) = 10$ and the velocities to $v_x(0, x, y) = -20 (y - 0.5)$, $v_y(0, x, y) = 20 (x - 0.5)$.
	In the transition from $r_0$ to $r_1$ we set the density to $\rho(0, x, y) = 1 + 9 f\big(r(x,y) \big)$ and the velocities to $v_x(0, x, y) = -20 (y - 0.5) f\big(r(x,y) \big)$, $v_y(0, x, y) = 20 (x - 0.5) f\big(r(x,y) \big)$ with $f\big(r(x,y) \big) \coloneqq \frac{r_1 - r(x,y)}{r_1}$.
	Outside of $r_1$, i.e., $r > r_1$ the medium is at rest $v_x(t_0, x, y) = 0$, $v_y(t_0, x, y) = 0$ and density is set to $\rho(0, x, y) = 1$.
	
	\bibliographystyle{elsarticle-num-names} 
	\bibliography{references.bib}
\end{document}